\documentclass{amsart}
\usepackage{tikz}
\usetikzlibrary{arrows,matrix}
\usepackage[all]{xy}

\begin{document}
\title[Arithmetic differential geometry in the arithmetic PDE setting, II]{Arithmetic differential geometry in the arithmetic PDE setting, II: curvature and cohomology}
\author{Alexandru Buium}
\address{Department of Mathematics and Statistics,
University of New Mexico, Albuquerque, NM 87131, USA}
\email{buium@math.unm.edu} 

\author{Lance Edward Miller}
\address{Department of Mathematical Sciences,  309 SCEN,
University of Arkansas, 
Fayetteville, AR 72701}
\email{lem016@uark.edu}

\def \cH{\mathcal H}
\def \cB{\mathcal B}
\def \d{\delta}
\def \ra{\rightarrow}
\def \bZ{{\mathbb Z}}
\def \cO{{\mathcal O}}
\newcommand{\Hom}{\operatorname{Hom}}

\def\uy{\underline{y}}
\def\uT{\underline{T}}

\newcommand{\OCp}{ {\mathbb C}_p^{\circ}  }

\newcommand{\sbt}{\,\begin{picture}(-1,1)(-1,-2)\circle*{2}\end{picture}\ }

\newtheorem{THM}{{\!}}[section]
\newtheorem{THMX}{{\!}}
\renewcommand{\theTHMX}{}
\newtheorem{theorem}{Theorem}[section]
\newtheorem{corollary}[theorem]{Corollary}
\newtheorem{lemma}[theorem]{Lemma}
\newtheorem{proposition}[theorem]{Proposition}
\theoremstyle{definition}
\newtheorem{definition}[theorem]{Definition}
\theoremstyle{definition}
\newtheorem{convention}[theorem]{Convention}
\newtheorem{notation}[theorem]{Notation}
\newtheorem{remark}[theorem]{Remark}
\newtheorem{example}[theorem]{\bf Example}
\numberwithin{equation}{section}
\subjclass[2000]{11S31, 17B56, 53B05, 53C21}

\begin{abstract}
This is the second paper in a series devoted to developing an arithmetic PDE analogue of Riemannian geometry. In Part 1 arithmetic PDE analogues of Levi-Civita and Chern connections were introduced and studied.
In this paper  arithmetic analogues of curvature and characteristic classes are developed.
\end{abstract}

\maketitle

\tableofcontents

\section{Introduction}

The present paper is Part 2 in a series;  the first part of this series \cite{BM22} will be referred to as Part 1. The aim of this series is to develop an arithmetic PDE version of the arithmetic ODE version of Riemannian geometry presented in \cite{BD, Bu17, Bu19}. Briefly, the ODE theory introduces analogous notions of familiar concepts from Riemann geometry where $p$-derivations, for $p$ a prime integer, are used instead of derivations. Essential concepts such as Chern and Levi-Civita connections are studied in detail and the associated theory of curvature for such connections surprisingly guards arithmetic information. In particular, the integers have some inherent curvature which can be seen as related to class field theory. This version of arithmetic geometry however is limited in that one works primarily over the completed maximal unramified extension of $\mathbb{Q}_p$ and as such there is a unique $p$-derivation to be had. Thus, this version of arithmetic Riemannian geometry is inherently of cohomogeneity one. 

\

Part I in detail adapts from the unramified context the basic structures of connections. In contrast, the focus of this article is curvature and associated cohomology theories which deepen the links between curvature and class field theory suggested in the unramified case. For the convenience of the reader the present paper is written so as to be logically independent of Part 1 in the sense that the relevant concepts and results of Part 1 will be (quickly) reviewed here as needed. However, we refer to the Introduction of Part 1
 for motivation and  context; we also refer to the body of Part 1 for more details. 
We include an Appendix which summarizes and explains the growing analogy with classical differential geometry from both Part 1 and the present paper.

\

\subsection{Quick review of Part 1}  \label{generall}
Following Part 1 we adopt the following 
 general conventions and notation {\it which will be enforced throughout the present paper}.
 
 \

Unless otherwise stated,  monoids will not necessarily have an identity. We denote by $\mathbb N$ the additive monoid of positive integers. We denote by $\mathbb M_n$ the free monoid with identity $0$ generated by the set $\{1,\ldots,n\}$; so the elements of the submonoid $\mathbb M_n^+:=\mathbb M_n\setminus\{0\}$ are words $\mu=i_1\ldots i_l$ with $i_1,\ldots,i_l\in \{1,\ldots,n\}$ and we denote by $\mathbb M_n^r$ the set of words of length $l\in \{1,\ldots,r\}$. For elements $\phi_1,\ldots,\phi_n$ in a monoid and for $\mu=i_1\ldots i_l\in \mathbb M_n^+$ we write $\phi_{\mu}:=\phi_{i_1}\ldots\phi_{i_l}$.

\

Associative $\mathbb Z$-algebras will   not be assumed commutative or with identity. By a ring we will mean a commutative associative $\mathbb Z$-algebra with identity. 
  For a (not necessarily algebraic) field extension $L/F$, denote by $\mathfrak G(L/F)$ the group of $F$-automorphisms of $L$.  For a square matrix $a$  with entries $a_{ij}$ 
in a ring we denote by $a^{(p^s)}$ the matrix whose entries are $a_{ij}^{p^s}$; we denote by $a^t$ the transpose of $a$. We adopt a similar notation for vectors instead of matrices.

\

We fix throughout an odd prime $p\in \mathbb Z$. For any ring $S$ and any Noetherian scheme $X$ we always denote by   $\widehat{S}$ and $\widehat{X}$ the respective $p$-adic completions.     
We denote by $R$  the $p$-adic 
completion of the maximum unramified extension of $\mathbb Z_p$
 and we consider the fields $K:=K[1/p]$, $k:=R/pR$.
 Let  $\pi$  be a prime element in a finite Galois extension of $\mathbb Q_p$; in Part 1 we eventually allowed $\pi$ to vary but in this paper we will fix $\pi$ once and for all.
Let $R_{\pi}:=R[\pi]$,  $K_{\pi}=K(\pi)$; we have $R_{\pi}/\pi R_{\pi}=k$.
  The degree $[K_{\pi}:K]$ of the field extension $K_{\pi}/K$ equals the ramification index $e=e(K_{\pi}/K)$ of this extension.  For every $\alpha\in k$ and $a\in R_{\pi}$ we write 
$\alpha=a$ mod $\pi$ if the image of $a$ in $k$ is $\alpha$.

\

  Fix an integer $s\in \mathbb N$.
 By  a {\bf higher $\pi$-Frobenius lift} of degree $s$ on a  ring $S$ containing $R_{\pi}$ we understand  a ring endomorphism $\phi$ of $S$ reducing to the $p^s$-power Frobenius modulo $\pi$. If $S$ is a $p$-torsion free ring containing $R_{\pi}$ then a  map $\delta:S\rightarrow S$ is called a {\bf higher $\pi$-derivation} of degree $s$ if the map $\phi:S\rightarrow S$ defined by $\phi(a):=a^{p^s}+\pi\delta(a)$ is a ring homomorphism (in which case it is a higher $\pi$-Frobenius lift of degree $s$); we say that $\phi$ and $\delta$ are {\bf attached} to each other. By a {\bf partial $\delta$-ring} of degree $s$ we understand a ring containing $R_{\pi}$ equipped with $n$ higher $\pi$-derivations of degree $s$. If we want to emphasize the degree $s$ in the notation we write $\phi^{(s)}$ and $\delta^{(s)}$ in place of $\phi$ and $\delta$, respectively.
 
 \

 Let $K^{\textup{alg}}$ be an algebraic closure of $K$ containing in $K_{\pi}$.
If $L$ be a subfield of $K^{\textup{alg}}$ containing $\mathbb Q_p$ by   a {\bf higher Frobenius automorphism} of degree $s\geq 1$ on  $L$  we understand  a continuous automorphism  $\phi\in \mathfrak G(L/\mathbb Q_p)$  that  induces the $p^s$-power automorphism on the residue field of the valuation ring  of $L$. For any subfield $L_0$ of $L$ we denote by $\mathfrak F^{(s)}(L/L_0)$ the set of higher Frobenius automorphisms on $L$ of degree $s$ that are the identity on $L_0$. Every element of $\mathfrak F^{(s)}(K^{\textup{alg}}/\mathbb Q_p)$ induces a higher $\pi$-Frobenius lift on $R_{\pi}$.

\

 We consider an $n$-tuple of $\pi$-Frobenius lifts, $(\phi_1^{(s)},\ldots,\phi_n^{(s)})$, of degree $s$ on $R_{\pi}$ and we consider the  higher $\pi$-derivations of degree $s$, $\delta_i^{(s)}:R_{\pi}\rightarrow R_{\pi}$, 
 attached to $\phi_i^{(s)}$. 
 We next fix an integer $N\geq 1$ and we consider the group scheme 
 $G:=\textup{GL}_{N/R_{\pi}}:=\textup{Spec}
 \ R_{\pi}[x,\det(x)^{-1}]$
  and the ring $\mathcal A:=\widehat{\mathcal O(G)}=R_{\pi}[x,\det(x)^{-1}]^{\widehat{\ }}$  
   where we denote by $x=(x_{ij})$  an $N\times N$ matrix of indeterminates. We define a {\bf $\pi$-connection} on $G$ of degree $s$ to be an $n$-tuple of higher $\pi$-Frobenius lifts 
   $\Phi^{(s)G}=((\phi_1^{(s)})^G,\ldots,(\phi_n^{(s)})^G)$ of degree $s$ on 
 $\mathcal A$ extending those on $R_{\pi}$. To give a $\pi$-connection is equivalent to giving the attached $n$-tuple of higher $\pi$-derivations
 $\Delta^{(s)G}=((\delta_1^{(s)})^G,\ldots,(\delta_n^{(s)})^G)$
 of degree $s$ on $\mathcal A$.  The numbers
$$e,s,n,N$$
 introduced above are not assumed to be related in general.
 For a $\pi$-connection as above we define the 
 {\bf Christoffel symbols} of the second kind to be the $n$-tuple of matrices $\Gamma^{(s)}=(\Gamma_1^{(s)},\ldots,\Gamma_n^{(s)}))$ whose transposed matrices are given by 
  \begin{equation}
\label{penu}\Gamma_i^{(s)t}:=(x^{(p^s)})^{-1} (\delta_i^{(s)})^G x
 \in \textup{Mat}_N(\mathcal A).\end{equation}
   Also we consider the matrices
    \begin{equation}
\Lambda_i^{(s)}:=(x^{(p^s)})^{-1}(\phi_i^{(s)})^G (x)=1+\pi\Gamma_i^{(s)t}\in \textup{GL}_N(\mathcal A)\end{equation}
and we  set $\Lambda^{(s)}=(\Lambda_1^{(s)},\ldots,\Lambda_n^{(s)})$.

\

 By a {\bf metric} we  understand a
 symmetric matrix $q\in \textup{GL}_N(R_{\pi})^{\textup{sym}}$; in the context of degree $s$ we write $q^{(s)}$ instead of $q$. (When $s=p^d$ for some $d\geq 1$ there is a conflict between our use of the superscript $(s)$ in $\Gamma_i^{(s)}$, $\Lambda_i^{(s)}$, $q^{(s)}$ where $(s)$ refers  to the degree  and the previously introduced definition $a^{(p^d)}=(a_{ij}^{p^d})$ for a matrix $a=(a_{ij})$. To resolve this conflict we make the convention that  $a^{(p^d)}$ will mean $(a_{ij}^{p^d})$   when the exponent $(p^d)$ is explicitly written as a $p$-power. On the other hand the exponent $(s)$ will be used to denote the degree  when $s$ is {\it not} explicitly written as a $p$-power.   We will say that  a $\pi$-connection as above is {\bf metric} with respect to $q^{(s)}$
   if the following diagrams are commutative for all $i$:
 $$
 \begin{array}{rcl}
\mathcal A& \stackrel{(\phi_{i,0}^{(s)})^G}{\longrightarrow} & \mathcal A \\
 \mathcal H_{q^{(s)}} \downarrow &\ &\downarrow \mathcal H_{q^{(s)}}\\
 \mathcal A & \stackrel{(\phi_i^{(s)})^G}{\longrightarrow} & \mathcal A\end{array}
$$
where $(\phi_{i,0}^{(s)})^G$ are the $\pi$-Frobenius lifts that send the matrix $x=(x_{jk})$ into the matrix $x^{(p^s)}:=(x_{jk}^{p^s})$ and $\mathcal H_{q^{(s)}}$ is the $R_{\pi}$-algebra homomorphism defined by $\mathcal H_{q^{(s)}}(x):=x^t q^{(s)} x$.  

\

Furthermore, in case $N=n$,
we define a {\bf torsion symbol} to be 
  an $n$-tuple $L^{(s)}=(L^{k(s)})$ of $n\times n$ 
  antisymmetric matrices $L^{k(s)}=(L^{k(s)}_{ij})$, 
  with entries in the ring $\widehat{R_{\pi}[y]}$ where $y=(y_1,\ldots,y_n)$ is an $n$-tuple of $n\times n$ matrix indeterminates.
    We say  that a $\pi$-connection 
   of degree $s$ is {\bf symmetric} with respect to $L^{(s)}$
   if the following equalities hold:
    \begin{equation}
\label{uf2}
\Gamma^{k(s)}_{ij}-\Gamma^{k(s)}_{ji}= L^{k(s)}_{ij}(\Lambda^{(s)}).\end{equation}
    One of our main results in Part 1 was the following arithmetic analogue of the Fundamental Theorem of Riemannian geometry (see \cite[Thm. 1.1]{BM22}):
    
     \begin{theorem}\label{boooster}
 For every metric $q^{(s)}$ and every torsion symbol $L^{(s)}$ there exists a unique $\pi$-connection of degree $s$ that is symmetric with respect to $L^{(s)}$ and metric with respect to $q^{(s)}$.
 \end{theorem}

The $\pi$-connection in the theorem above is referred to as the {\bf arithmetic Levi-Civita connection} of degree $s$ attached to $L^{(s)}$ and $q^{(s)}$.
In Part 1 we also introduced {\bf arithmetic Chern connections} of degree $s$ attached to metrics $q^{(s)}$ which are defined for $N$ not necessarily equal to $n$. Other topics were covered in Part 1 (geodesics, $\delta$-convergence, etc.) which we will not review in this Introduction.  

\

\subsection{Main results of the present paper}
Our main concern in this paper is to introduce and study {\bf curvature} and {\bf characteristic classes}. 
 Assume the setting of Subsection \ref{generall} and assume, until further notice, that $n=N$.
 Also  it will be convenient to  fix a positive integer $c$; its significance comes from the local class field theory examples for our theory, see Remark \ref{lcftc}.

\

We start by considering the following sequence of sets
\begin{equation}
\label{dda}
(\mathfrak D^{(s)})_{s\in c\mathbb N}\end{equation}
where $\mathfrak D^{(s)}=\{\phi_1^{(s)},\ldots,\phi_n^{(s)}\}$ and $\phi_i^{(s)}$ are higher $\pi$-Frobenius lifts 
 of degree $s$ 
on $R_{\pi}$ such that the set 
$\mathfrak D:=\cup_{s\in c\mathbb N}\mathfrak D^{(s)}$ 
is closed under composition; so $\mathfrak D$ is a monoid which  will be referred to
 as the {\bf Weil monoid} 
by analogy with the Weil group in local class field theory; see \cite[p. 69]{N80}.
We view $\mathfrak D$ as an analogue of an involutive distribution in the tangent bundle of a manifold; see Section 9. Assume in addition that we are given a sequence of metrics 
\begin{equation}
\label{ddaa}
(q^{(s)})_{s\in c\mathbb N}\end{equation}
and a sequence of torsion symbols 
\begin{equation}
\label{ddaaa}
(L^{(s)})_{s\in c\mathbb N}.\end{equation} 

\

\noindent The corresponding sequence $$(\Delta^{(s)G})_{s\in c\mathbb N}$$
of arithmetic Levi-Civita connections
$\Delta^{(s)G}=((\phi_1^{(s)})^G,\ldots, (\phi_n^{(s)})^G)$   attached to the data (\ref{dda}), (\ref{ddaa}), (\ref{ddaaa})
 will be referred to as the {\bf graded arithmetic Levi-Civita connection}.  To the latter one can attach 
 the family of elements
 $$R_{ijl}^{k(s,r)}\in \mathcal A, \ i,j,k,l\in \{1,\ldots,n\},\ \ s,r\in c\mathbb N$$ defined by
$$R_{ijl}^{k(s,r)}:=
\frac{1}{\pi}((\phi_i^{(s)})^G(\phi_j^{(r)})^G-(\phi_j^{(r)})^G(\phi_i^{(s)})^G
-(\phi_i^{(s)}\phi_j^{(r)})^G+(\phi_j^{(r)}\phi_i^{(s)})^G)(x_{kl}).$$

\

One defines the {\bf arithmetic Riemann curvature tensor} (see Definition \ref{arct}) as the family of elements 
  $$R^{(s)}_{ijkl}\in \mathcal A,\ \ \ i,j,k,l\in \{1,\ldots,n\},\ \ s\in c\mathbb N,$$
 obtained by making $s=r$ and ``lowering the indices":
$$R^{(s)}_{ijkl}:=\sum_m R_{ijk}^{m(s,s)} (q_{ml}^{(s)})^{p^{2s}}\in \mathcal A.$$
For $s=c$ we write $R_{ijkl}:=R_{ijkl}^{(c)}$.  Denoting by 
  $\mathcal M$ the maximal ideal of $\mathcal A$ corresponding
 to the identity element of the $k$-group scheme $G\otimes_{R_{\pi}} k$
 and assuming the Weil monoid is abelian
   we will 
  derive formulae for the images
  $\overline{R}_{ijkl}\in k=\mathcal A/\mathcal M$ of $R_{ijkl}$
   that are strikingly similar to the classical formulae in Riemannian geometry.  
   Specifically consider the higher $\pi$-Frobenius lift on $R_{\pi}$ of degree $2c$ given by
    $\phi_{ij}^{(2c)}:=\phi_i^{(c)}\phi_j^{(c)}$ and let $\delta_{ij}^{(2c)}$ be the higher $\pi$-derivation of degree $2c$ on $R_{\pi}$ attached to  $\phi_{ij}^{(2c)}$. Also let $L^{(c)}(1)$ be obtained from $L^{(c)}\in (\widehat{R_{\pi}[y]})^n$ by setting all matrices $y_1,\ldots,y_n$ equal to the identity matrix. 
   We will prove the following result (see Theorem \ref{tery}):
     
  \begin{theorem}
  Assume the Weil monoid is abelian and $L^{(c)}(1)=0$.
  Then
  $$
\overline{R}_{ijkl} =  -\frac{1}{2} \cdot (\delta_{jl}^{(2c)}q_{ik}-
\delta_{jk}^{(2c)}q_{il}-\delta_{il}^{(2c)}q_{jk}+\delta_{ik}^{(2c)}q_{jl}
) \ \ \  \textup{mod}\ \ \pi.$$ 
 \end{theorem}

 In particular, the family $(\overline{R}_{ijkl})$ satisfies the symmetries of the classical Riemann curvature tensor; see Corollary \ref{antisymm}. The formula in the above theorem is, of course, analogous to the formula for the ``order $2$ part" of the classical Riemann curvature tensor; cf. Equation \ref{zzoo} in Section 9. 

 \

Given a metric $q\in \textup{GL}_n(R_{\pi})^{\textup{sym}}$ with invertible diagonal entries
 (a condition which one can view as an analogue of positivity), there are  canonical ways to 
attach to  $q$  a sequence of metrics (\ref{dda}) with $q^{(c)}=q$; see Definition \ref{canmet}. Also given sets of $\pi$-Frobenius lifts  (\ref{ddaa}) there are canonical ways to define  a sequence of torsion symbols (\ref{ddaaa}) which will be referred to as the (additive, respectively multiplicative) {\bf canonical torsion symbols}; see Definition \ref{canonicaltorsion}.
 Note that if the Weil monoid is abelian then the canonical torsion symbols  will automatically satisfy $L^{(c)}(1)=0$ (see Remark \ref{agencyy}); so for these symbols the hypothesis $L^{(c)}(1)=0$ can be dropped from our theorem.

  
  \

We next turn our attention to characteristic classes. It is these classes which hold promise of deeper connections to class field theory. To keep the introduction brief, we giver here only a very rough, non-rigorous description to convey the ideas in summary. Consider in what follows arbitrary  {\bf graded $\pi$-connections} by which we mean  families $(\Delta^{(s)G})_{s\in c\mathbb N}$ of $\pi$-connections. We will attach to any such graded $\pi$-connection a class belonging to a certain second  Hochschild cohomology group $H^2_{\textup{Hoch}}$ which we call the {\bf Hochschild characteristic class}; this will induce, in particular, a class belonging to a  second Lie cohomology group $H^2_{\textup{Lie}}$ which we call the {\bf Lie characteristic class}.  

  \

  The above cohomology classes will be images of a purely group theoretic cohomology class which is independent of the choice of a graded $\pi$-connection and can be viewed as  analogous to a ``topological" class. The latter  class is attached to purely Galois theoretic data extracted from $K_{\pi}/K$ so the situation is  somewhat similar  to that in Chern-Weil theory where topological classes are represented by forms constructed from connections. 
   We will prove (see Corollary \ref{miarevenit} and Remark \ref{sidecenu} for precise statements) the following:
  
  \begin{theorem}\label{treillea}
  Assume $n\geq 2$ and the Weil monoid is abelian. Then, for  ``sufficiently general metrics,"  the Lie characteristic class (and hence the Hochschild characteristic class) of the graded arithmetic Levi-Civita connection  is non-trivial.
  \end{theorem}
  
  The ordering of the higher $\pi$-Frobenius lifts $\phi_1^{(s)},\ldots,\phi_n^{(s)}$
  in each degree $s$ will be determined by what we shall call {\bf labeling symbols}.
 The action of   the symmetric group $\Sigma_n$ on labeling symbols will induce an action on   curvature whose invariants we will consider in  Section \ref{invpol}; this will lead to analogues of scalar curvature. The symmetric group plays the role of group of ``coordinate changes"; this is consistent with the fact that  the geometry implicitly underlying our theory can be viewed as a    geometry over the field with one element; cf. \cite[Introduction]{Bu05} and \cite{Bor}.

 \

 A similar theory can be developed for the {\bf Chern connection} in which $\textup{GL}_n$ is replaced by $\text{GL}_N$ with $N$ not necessarily equal to $n$. The curvature of the arithmetic  Chern connection is generally non-trivial but,  in case the Frobenius lifts commute on $R_{\pi}$,  it  vanishes  mod $\pi$ (cf. Theorem \ref{wqq} for a precise statement). 
 
 \

 We will also introduce a ``multiplicative" version of curvature  and we will show that this multiplicative version
and the previous ``additive" version ``agree mod $\pi$ at the origin" (see Theorem \ref{ursuh} for a precise statement).
 The multiplicative curvature is a feature of our arithmetic perspective and doesn't have a direct classical counterpart.
 
 \

 The theory above  corresponds to the group $\textup{GL}_N$ viewed as a trivial torsor under itself. One can ask for a variant of the theory that involves non-trivial torsors.
 Since no non-trivial torsors for the algebraic group $\textup{GL}_N$ over a field exist
 (Hilbert's Theorem 90) we will replace the ``gauge group"  $\textup{GL}_N$ by a smaller group that has a 
 non-trivial first Galois cohomology group $H^1_{\textup{gr}}$. Inspired again by the link between geometry of $p$-derivations and the philosophy of the field with one element we will consider as ``gauge group" the semidirect product  $\mathfrak W$ of the  Weyl group $\Sigma_N$ of $\textup{GL}_N$ with the group of diagonal matrices whose entries are roots of unity in $R$.  We will then discuss  conditions under which the arithmetic Levi-Civita and Chern connections ``descend" to non-trivial torsors. The Lie characteristic class in the Lie algebra cohomology group $H^2_{\textup{Lie}}$ attached to  a connection (cf. Theorem \ref{treillea})  and the cohomology class in the Galois  cohomology group $H^1_{\textup{gr}}$ attached to  torsors (that we just mentioned) will not be related in our theory; intuitively  they will measure ``twists" in different directions.
 
 \

 Finally, we suspect these cohomology classes to be the start of a much deeper study relating concepts in arithmetic Riemannian geometry and class field theory. It would be beyond scope here to fully flesh out an arithmetic Chern-Weil theory however to justify these perspectives we will explain how the Christoffel symbols in our arithmetic PDE context are related to, and hence can be viewed as matrix analogues of, the classical  Legendre symbols in number theory; see Theorems \ref{legendre1} and \ref{recc}.
 
 \
 
\subsection{Leitfaden} Section 2  contains a series of preliminaries on monoids, algebras, and cohomology; in particular we introduce here  
a series of ``symbols" (Lie, associative, and labeling symbols)  which are derived from purely Galois theoretic data. 
Section 3
develops the formalism of curvature in the arithmetic PDE setting. Section 4 shows
 shows how to canonically attach sequences of metrics $(q^{(s)})$ and sequences $(L^{(s)})$
 of torsion symbols to a given metric $q$ and to given labeling symbols.
Section 5
 is concerned with invariant polynomials under the action of $\Sigma_n$ on the Riemannian curvature tensor. 
Section 6 discusses the multiplicative version of curvature. Section 7 discusses the gauge group and  torsors. Section 8 discusses the connection with Legendre symbols.
 The Appendix then offers a brief discussion of the analogies between our theory and classical Riemannian concepts.

\subsection{Acknowledgements} The first author was partially supported by the Simons Foundation through award 615356.

\section{Weil monoids}

In order to  introduce and study arithmetic curvature and characteristic classes we need a series of preliminaries on monoids, algebras, and their cohomology.

\subsection{Associative, Lie, and labeling symbols} \label{labelll}
 We start, in this Subsection, with an abstract preparation involving monoids; recall that, according to our general conventions, our monoids     are not assumed commutative and are not  assumed to have an identity element.

\begin{definition}\label{gm}
By a {\bf graded monoid} we understand a monoid $\mathfrak D$ together with a  
 monoid homomorphism $\deg:\mathfrak D\rightarrow \mathbb N$ whose image is of the form $c\mathbb N$ for some $c\in \mathbb N$.  We call $c$ the {\bf minimal degree} of $\mathfrak D$. We write
$\mathfrak D^{(s)}:=\deg^{-1}(s)$ for all $s\in c\mathbb N$. (Note that such a monoid does not have an identity.) \end{definition}

Graded monoids form a category where the morphisms are monoid homomorphisms that are compatible with the degree maps. Recall that a monoid is called {\bf cancellative} if it has the left and right cancellation property.
If a graded monoid as above  is cancellative  then the cardinalities of the sets $\mathfrak D^{(s)}$ satisfy
$$|\mathfrak D^{(c)}|\leq |\mathfrak D^{(2c)}|\leq |\mathfrak D^{(3c)}|\leq \cdots$$

In what follows if $\mathfrak S$ is a finite group, $c\in \mathbb N$,  and $\theta\in \textup{Aut}_{\textup{gr}}(\mathfrak S)$ is a group automorphism of $\mathfrak S$  we denote by $c\mathbb N\times_{\theta} \mathfrak S$ the semidirect product monoid whose multiplication  is defined by
\begin{equation}
\label{mamm}
(s_1,\sigma_1)(s_2,\sigma_2):=(s_1+s_2,\theta^{s/c}(\sigma_1)\sigma_2),\ \ s_1,s_2\in c\mathbb N,\  \ \sigma_1,\sigma_2\in \mathfrak S.
\end{equation}

\

\begin{proposition} \label{tdceeee}
Let $\mathfrak D$ be a graded monoid of minimal degree $c$. The following conditions are equivalent:

1) $\mathfrak D$  is cancellative and for all $s\in c\mathbb N$ the sets $\mathfrak D^{(s)}$ are finite of the same cardinality.

2)  
$\mathfrak D$ is isomorphic as a graded monoid with a semidirect product $c\mathbb N\times_{\theta}\mathfrak S$ where $\mathfrak S$ is a finite group and $\theta \in \textup{Aut}_{\textup{gr}}(\mathfrak S)$. 

3) 
There exists a group  homomorphism $\deg:\mathfrak G\rightarrow \mathbb Z$ with finite kernel such that   $\mathfrak D$ and  $\deg^{-1}(c\mathbb N)$ (with their degree maps
$\deg:\mathfrak D\rightarrow \mathbb Z$ and $\deg:\deg^{-1}(c\mathbb N)\rightarrow \mathbb Z$) are isomorphic as graded monoids. \end{proposition}

The types of monoids that appear most frequently in this article are those that meet the conditions Proposition~\ref{tdceeee} and we codify this condition as follows. 

\begin{definition}\label{rgm}
A  graded monoid $\mathfrak D$ of minimal degree $c$ is called a {\bf Weil monoid} if it  satisfies the equivalent conditions in Proposition \ref{tdceeee}. The cardinality $|\mathfrak D^{(c)}|$ is called the {\bf size} of $\mathfrak D$. If $\mathfrak D\simeq c\mathbb N\times_{\theta}\mathfrak S$ as in assertion 2 of 
Proposition \ref{tdceeee} we say that $\mathfrak D$ is {\bf represented by}  $(\mathfrak S,\theta)$.
\end{definition}

\begin{definition}
A subset of a graded monoid $\mathfrak D$ is {\bf abelian} if its elements are pairwise commuting. 
\end{definition} 

Before giving the proof of the proposition above we state the impact as the following corollaries. 

\begin{corollary}\label{ttoo}
Assume a Weil monoid $\mathfrak D$ of minimal degree $c$ is represented by $(\mathfrak S,\theta)$.
Then the following are equivalent:

1) $\mathfrak D^{(c)}$ is abelian;

2) $\mathfrak D$ is an abelian monoid;

3) $\mathfrak S$ is an abelian group and  $\theta=\textup{id}$.\end{corollary}

{\it Proof}. To see that assertion 2 implies assertion 3 take $s_1=s_2=s$ with $\theta^{s/c}=\textup{id}$  in Equation \ref{mamm} to get that $\mathfrak S$ is abelian; then take $s_1=s_2=c$ and $\sigma_2=1$ to see that $\theta=\textup{id}$. The rest of the implications follow.
\qed

\begin{corollary} \label{tdce}
Let $\mathfrak D$ be a Weil monoid of minimal degree $c$. Then
for every $\phi\in \mathfrak D^{(c)}$ there exists $t\in \mathbb N$  such that $\phi^t$ commutes with all the elements of $\mathfrak D$.
\end{corollary}

\

In order to prove Proposition \ref{tdceeee} (and for later reference) we start with the following lemma.

\begin{lemma}\label{prelop}
Assume $\mathfrak D$ is a graded monoid of minimal degree $c$ which is a submonoid of a group $\mathfrak H$ and assume that for all $s\in c\mathbb N$ the sets $\mathfrak D^{(s)}$ are finite of the same cardinality $n$. Let $\mathfrak D^{(c)}=\{\phi_1,\ldots,\phi_n\}$, let $\sigma_i:=\phi_1^{-1}\phi_i\in \mathfrak H$ for $i\in \{1,\ldots,n\}$ and let $\mathfrak S:=\{\sigma_1,\ldots,\sigma_n\}\subset \mathfrak H$. Then the following hold:

1) $\mathfrak S$ is a subgroup of $\mathfrak H$.

2) $\mathfrak S$ is normalized by $\phi_1$.

3) For all $s\in c\mathbb N$ we have $\mathfrak D^{(s)}=\{(\phi_1^{(c)})^{s/c}\sigma_i\ |\ 1\leq i\leq n\}$.
\end{lemma}

{\it Proof}.
In order to check assertion 1  it is sufficient to check that for  $i,j\in \{1,\ldots,n\}$
there exists $k\in \{1,\ldots,n\}$ such that 
\begin{equation}
\label{uunu}
\sigma_i\sigma_j=\sigma_k.\end{equation}
 Fix $i$ and $j$. The condition 
$|\mathfrak D^{(2c)}|=n$ implies that  we have
\begin{equation}
\label{forsure}
\mathfrak D^{(2c)}=\{\phi_1\phi_1,\phi_1\phi_2,\ldots,\phi_1\phi_n\}=\{\phi_1\phi_1,\phi_2\phi_1,\ldots,\phi_n\phi_1\}.\end{equation}
We get that there exists $l\in \{1,\ldots,n\}$ such that
\begin{equation}
\label{ddoi}
\phi_j\phi_1=\phi_1\phi_l.\end{equation}
 Now since 
$$\mathfrak D^{(2c)}=\{\phi_1\phi_l,\phi_2\phi_l,\ldots,\phi_n\phi_l\}$$
 there exists $k\in \{1,\ldots,n\}$ such that \begin{equation}
 \label{ttrei}
 \phi_i\phi_l=\phi_k\phi_1.\end{equation}
 We get
 $$\begin{array}{rcl}
 \sigma_i\sigma_j & = & \phi_1^{-1}\phi_i\phi_1^{-1}\phi_j\\
 \ & \ & \ \\
 \ & = & \phi_1^{-1}\phi_k\phi_1\phi_l^{-1}
 \phi_1^{-1}\phi_j \ \ \ \textup{by}\ \ (\ref{ttrei})\\
 \ & \ & \ \\
 \ & = & \phi_1^{-1}
 \phi_k\phi_1(\phi_1\phi_l)^{-1}\phi_j\\
 \ & \ & \ \\
 \ & = &  \phi_1^{-1}
 \phi_k\phi_1(\phi_j\phi_1)^{-1}\phi_j  \ \ \ \textup{by}\ \ (\ref{ddoi})\\
 \ & \ & \ \\
 \ & = & \phi_1^{-1}\phi_k\\
 \ & \ & \ \\
 \ & = & \sigma_k.
 \end{array}
 $$
  which proves (\ref{uunu}). 

  \

In order to check assertion 2 it is sufficient to check that  for all $i\in \{1,\ldots,n\}$ there exists  $j\in \{1,\ldots,n\}$ such that $\phi_1^{-1}\sigma_i\phi_1=\sigma_j$, equivalently $\phi_i\phi_1=\phi_1\phi_j$. This is, however clear from the equality (\ref{forsure}). Assertion 3 follows directly from the fact that $|\mathfrak D|=n$.
\qed

\bigskip

{\it Proof of Proposition \ref{tdceeee}}. The implication $3 \Rightarrow 1$ is trivial.

\

To check $1\Rightarrow 2$ note first that condition $1$  implies the following property: 
\begin{equation}
\label{reversible}
\mathfrak D \phi \cap \mathfrak D \phi'\neq \emptyset\ \ \ \textup{for all}\ \ \ \phi,\phi'\in \mathfrak D.\end{equation}
Indeed  if $\phi\in \mathfrak D^{(s)}$ and $\phi'\in \mathfrak D^{(r)}$ then, using the equality
$|\mathfrak D^{(s+r)}|=|\mathfrak D^{(s)}|$ and cancellativity,   there exists $\phi''\in \mathfrak D^{(s)}$ such that $\phi'\phi=\phi''\phi'$.
 On the other hand, by a classical result in monoid theory (\cite{CP61}, p. 35) every cancellative monoid with  property (\ref{reversible}) is embeddable into a group. Consider  an embedding of $\mathfrak D$ into a group $\mathfrak H$. We use now the notation in Lemma \ref{prelop} and we consider 
 the semidirect product  $c\mathbb N\times_{\theta} \mathfrak S$  where
 $\theta(\sigma_i):=\phi_1^{-1}\sigma_i\phi_1$.
  Then we have an isomorphism of graded monoids
 $$c\mathbb N\times_{\theta} \mathfrak S\rightarrow \mathfrak D,\ \ \ (s,\sigma_i)\mapsto \phi_1^{s/c}\sigma_i=
 \phi_1^{(s/c)-1}\phi_i\in \mathfrak D^{(s)}.$$
 To check $2\Rightarrow 3$ note that if we have an isomorphism of graded monoids $\mathfrak D\simeq c\mathbb N\times_{\theta} \mathfrak S$ then 
  $c\mathbb N\times_{\theta} \mathfrak S$ is the preimage of $c\mathbb N$
 under the group homomorphism $c\mathbb Z\times_{\theta} \mathfrak S\rightarrow \mathbb Z$ given by the first projection. \qed
 
 \begin{remark}
 If $\mathfrak D$ is a Weil monoid then for all $x\in \mathfrak D^{(r)}$ and $y\in \mathfrak D^{(s)}$ with $r> s$ there exist unique elements denoted by $x^{-1}y\in \mathfrak D^{(r-s)}$ and $yx^{-1}\in \mathfrak D^{(r-s)}$, respectively, such that 
 $$(yx^{-1})x=y\ \ \ \textup{and}\ \ \ x(x^{-1}y)=x.$$
 \end{remark}

\

We have need of careful bookkeeping to handle products in a Weil monoid. As such, we introduce a labeling scheme based on the introduction of various symbols. These play an important role in the results and we introduce canonical and natural choices for such symbols throughout. 

 \begin{definition}\label{ccaannabstract}
 Assume $\mathfrak D\rightarrow \mathbb N$ is a Weil monoid of minimal degree $c$ and size $n$.  A {\bf labeling symbol} for $\mathfrak D^{(s)}$ (where $s\in c\mathbb N$)
 is a bijection 
 $$\omega^{(s)}:\{1,\ldots,n\}\rightarrow \mathfrak D^{(s)}.$$
 A {\bf labeling symbol} for $\mathfrak D$ is a sequence $(\omega^{(s)})_{s\in c\mathbb N}$
 where $\omega^{(s)}$ is labeling symbol for $\mathfrak D^{(s)}$.
 Fix, in what follows, a labeling symbol  for $\mathfrak D$.
 We write $$\phi_i^{(s)}:=\omega^{(s)}(i).$$
 For $i,j\in \{1,\ldots,n\}$ and $s,r\in c\mathbb N$ we  define 
 $(i\star j)_{s,r}\in  \{1,\ldots,n\}$ to be the index with the property that
 $$\phi_i^{(s)}\phi_j^{(r)}=\phi_{(i\star j)_{s,r}}^{(s+r)}.$$
 For $s=r=c$ we simply write
 $$i\star j:=(i\star j)_{c,c}.$$
 We define  the {\bf associative symbol} (corresponding to our labeling symbol for $\mathfrak D$) to be 
 the $n$-tuple of matrices
 $$\alpha^{(s,r)}=(\alpha^{1(s,r)},\ldots,\alpha^{n(s,r)})$$ with entries in $\{0,1\}$ defined by letting
 $$\alpha_{ji}^{k(s,r)}=1\ \ \ \Longleftrightarrow\ \ \ k=(i\star j)_{s,r}.$$
  Note that for each $k$ and each $s$ the matrix $\alpha^{k(s,r)}$ is a permutation matrix, i.e.,
 a matrix obtained from the identity matrix by permuting its columns.
 Set
  \begin{equation} \ell_{ij}^{k(s,r)}:=\alpha_{ji}^{k(s,r)}-\alpha_{ij}^{k(r,s)} \in \{-1,0,1\}.\end{equation}

  \

 \noindent Set also $$\alpha_{ij}^{k(s)}:=\alpha_{ij}^{k(s,s)}, \alpha^{k(s)}=\alpha^{k(s,s)}, \alpha^{(s)}=\alpha^{(s,s)},$$ and
  $\ell_{ij}^{k(s)}:=\ell_{ij}^{k(s,s)}$. We then have $\ell^{k(s)}
_{ij}+\ell^{k(s)}_{ji}=0$; in particular for each $s$ we may consider the tuple of antisymmetric matrices $\ell^{(s)}:=(\ell^{1(s)},\ldots,\ell^{n(s)})$, with
 $\ell^{k(s)}:=(\ell_{ij}^{k(s)})$. We call $\ell^{(s)}$ the {\bf Lie symbol}.  
   \end{definition}

   \begin{remark}Note that the only labeling symbols that $\alpha_{ij}^{k(s,r)}$ and $\ell_{ij}^{k(s,r)}$ depend on are $\omega^{(s)},\omega^{(r)},\omega^{(s+r)}$.
   \end{remark}
   
   \

  One would like to make the labeling symbols $\omega^{(s)}$ ``compatible" as $s$ varies.
  A way to achieve that is by considering the following notion of coherence.
  
  \begin{definition}\label{cfff}
  A labeling symbol $(\omega^{(s)})_{s\in c\mathbb N}$ is called {\bf coherent} if for all $s\in c\mathbb N$ and all $i,j\in \{1,\ldots,n\}$ we have
  $$\phi_i^{(s+c)}(\phi_i^{(s)})^{-1}=\phi_j^{(s+c)}(\phi_j^{(s)})^{-1}.$$
   \end{definition}
  
We give now a classification of coherent symbols and show how each labeling symbol can be canonically modified to give a coherent symbol.

  \begin{lemma}
  Note that to every function $\gamma:\mathbb N\rightarrow \{1,\ldots,n\}$ and every labeling symbol
  $\omega:=\omega^{(c)}:\{1,\ldots,n\}\rightarrow \mathfrak D^{(c)}$ on $\mathfrak D^{(c)}$  one can attach a coherent labeling symbol $(\omega^{(s)})_{s\in c\mathbb N}$ on $\mathfrak D$ defined inductively by:
  \begin{equation}
  \label{cocomm}
  \phi^{(s+c)}_i:=\phi^{(c)}_{\gamma(s)}\phi^{(s)}_i\ \ \textup{for all}\ \ i\in \{1,\ldots,n\}.
  \end{equation}
   Every coherent labeling symbol on $\mathfrak D$ arises in this way for a unique $\omega$ and a unique function $\gamma$ as above.
   \end{lemma}

   \

   So, in the situation above, if one defines the words
  $$\mu^{(s)}_i:=\gamma(s-c)\gamma(s-2c)\ldots\gamma(2c)\gamma(c)i\in \mathbb M_n^{(s/c)}$$
  then
  $$\phi_i^{(s)}=\phi_{\mu_i^{(s)}}^{(c)}=\phi^{(c)}_{\gamma(s-c)}\phi^{(c)}_{\gamma(s-2c)}\ldots\phi^{(c)}_{\gamma(2c)}\phi^{(c)}_{\gamma(c)}\phi^{(c)}_i.$$
  If $\gamma$ is the constant function $\gamma=\gamma_h:c\mathbb N\rightarrow \{1,\ldots,n\}$, $\gamma_h(s)=h$ with $h\in \{1,\ldots,n\}$ we call $\omega^{(s)}=:\omega^{(s)}_h$ the {\bf canonical labeling symbol} attached to $\omega$ and $h$ in which case
  $$\phi^{(s)}_i=(\phi^{(c)}_h)^{s/c-1}\phi^{(c)}_i.$$
   
   \

 \subsection{Review of algebras and cohomology}
 
In this subsection we fix our conventions and notation for algebras, modules, and their cohomology and we review some basic facts needed in what follows. 

\

\subsubsection{Algebras}
Associative $\mathbb Z$-algebras will   not be assumed commutative or with identity. Recall that by a ring we mean a commutative associative $\mathbb Z$-algebra with identity. For the general context see \cite[Sect. 7 and 9]{W97}. Although, unlike \cite{W97}, our associative $\mathbb Z$-algebras are {\it not} assumed to have a unit. 

\

If $\mathfrak G$ is a monoid acting on a ring $S$ by ring endomorphisms  the {\bf skew monoid algebra} $S[\mathfrak G]$ attached to this action is the associative $\mathbb Z$-algebra   defined as the additive group of all formal finite $S$-linear combinations of elements of $\mathfrak G$  with   multiplication defined by  $(s_1g_1)(s_2g_2) =s_1g_1(s_2)g_1g_2$ for $g_1,g_2\in \mathfrak G$, $s_1,s_2\in S$. If $S=\mathbb Z$ we refer to $\mathbb Z[\mathfrak G]$ as the {\bf monoid $\mathbb Z$-algebra} attached to $\mathfrak G$.

\

For an associative $\mathbb Z$-algebra $\mathfrak a$ an {\bf $\mathfrak a$-$\mathfrak a$-bimodule} is an abelian group with a structure of left $\mathfrak a$-module and a structure of right $\mathfrak a$-module that commute with each other. By an ideal in an associative $\mathbb Z$-algebra we understand a  two-sided ideal; hence every ideal has a structure of bimodule.  
For an associative $\mathbb Z$-algebra $\mathfrak a$ left $\mathfrak a$-modules will be simply referred to as {\bf $\mathfrak a$-modules}.
If $\mathfrak a$ is an associative $\mathbb Z$-algebra and also an $S$-module (where $S$ is a ring) we say the two structures are {\bf compatible} if $(sx)y=s(xy)$ for all $s\in S$ and $x,y\in \mathfrak a$. The skew monoid algebra attached to an action of a monoid on a ring $S$ has a natural  $S$-module structure compatible with the associative $\mathbb Z$-algebra structure.
For an $\mathfrak a$-$\mathfrak a$-bimodule we have an induced $\mathfrak a$-module structure and an induced $\mathfrak a^{\textup{op}}$-module structure where $\mathfrak a^{\textup{op}}$ is the opposite associative $\mathbb Z$-algebra. If $\mathfrak a$ has an identity then an $\mathfrak a$-$\mathfrak a$-bimodule is the same as a $\mathfrak a\otimes_{\mathbb Z} \mathfrak a^{\textup{op}}$-module.

\

Modules over a Lie $\mathbb Z$-algebra $\mathfrak b$ will be referred to as {\bf Lie $\mathfrak b$-modules}; ideals in a  Lie $\mathbb Z$-algebra $\mathfrak b$ will be called {\bf Lie ideals}. For an associative $\mathbb Z$-algebra $\mathfrak a$ we will continue to denote by $\mathfrak a$ the Lie $\mathbb Z$-algebra whose abelian group is $\mathfrak a$ and whose bracket is the commutator; in this case if $M$ is an $\mathfrak a$-module then the multiplication map $\mathfrak a\times M\rightarrow M$ also makes $M$ a Lie $\mathfrak a$-module and the $\mathfrak a$-submodules of $M$ are the same as the Lie $\mathfrak a$-submodules of $M$. 
 Note that if $\mathfrak a$ is an associative $\mathbb Z$-algebra viewed as a Lie $\mathbb Z$-algebra then the structure of Lie $\mathfrak a$-module on $\mathfrak a$ induced by the structure of $\mathfrak a$-module on $\mathfrak a$ is {\it not} the adjoint Lie $\mathfrak a$-module structure on $\mathfrak a$.

 \

 Every ideal in an associative $\mathbb Z$-algebra $\mathfrak a$ is a Lie ideal in the Lie $\mathbb Z$-algebra $\mathfrak a$ but, of course, not every Lie ideal in $\mathfrak a$ is an ideal and not every Lie submodule of $\mathfrak a$ is a Lie ideal. If $\mathfrak a$ is an associative $\mathbb Z$-algebra and $M$ is an $\mathfrak a$-$\mathfrak a$-bimodule then $M$ has a natural structure of Lie $\mathfrak a$-module given by the commutator; this Lie $\mathfrak a$-structure  coincides with the Lie $\mathfrak a$-module structure defined by the $\mathfrak a$-module structure on $M$  if and only if the $\mathfrak a^{\textup{op}}$-module structure of $M$ is trivial.
 
 \

A {\bf graded} associative (respectively Lie) $\mathbb Z$-algebra is an associative (respectively  Lie) $\mathbb Z$-algebra $\mathfrak u$
together with an abelian group decomposition
$\mathfrak u=\bigoplus_{s\in \mathbb Z}\mathfrak u^{(s)}$
such that $\mathfrak u^{(r)}\mathfrak u^{(s)}\subset \mathfrak u^{(r+s)}$ (respectively 
$[\mathfrak u^{(s)},\mathfrak u^{(r)}]\subset \mathfrak u^{(s+r)}$) for all $s,r\in \mathbb Z$.

  \

\subsubsection{Lie algebra cohomology}
Let $\mathfrak c$ be a  Lie $\mathbb Z$-algebra and let $\mathfrak a$ be a  Lie $\mathfrak c$-module.
We denote by $Z^2_{\textup{Lie}}(\mathfrak c,\mathfrak a)$ the group of all  Lie $2$-cocycles $f: \mathfrak c\times \mathfrak c\rightarrow \mathfrak a$, i.e., alternating $\mathbb Z$-bilinear maps such that for all $x,y,z\in \mathfrak c$ we have
$$f([x,y],z)+f([y,z],x)+f([z,x],y)=x f(y,z)+y f(z,x)+z f(x,y).
$$
 We denote by $B^2_{\textup{Lie}}(\mathfrak c,\mathfrak a)$ the group of all Lie $2$-coboundaries $f:\mathfrak c\times \mathfrak c\rightarrow \mathfrak a$, i.e., maps of the form
 $$f(x,y)=xg(y)-yg(x)-g([x,y]),\ \ x,y\in \mathfrak c,$$
 where $g:\mathfrak c\rightarrow \mathfrak a$ is a $\mathbb Z$-module homomorphism.
  We consider the $2$nd Chevalley-Eilenberg cohomology group
$$H^2_{\textup{Lie}}(\mathfrak c,\mathfrak a):=Z^2_{\textup{Lie}}(\mathfrak c,\mathfrak a)/
B^2_{\textup{Lie}}(\mathfrak c,\mathfrak a).$$
We also recall that a Lie $\mathfrak c$-module $\mathfrak a$ is called {\bf trivial} if the multiplication map $\mathfrak c\times \mathfrak a\rightarrow \mathfrak a$ is the zero map.

\

The $H^2_{\textup{Lie}}$ group defined above is not a priori
 isomorphic to the second cohomology group defined via derived functors 
on the category of $\mathfrak c$-modules; however such an isomorphism exists if 
$\mathfrak c$ is a free $\mathbb Z$-module  
 or if 
both $\mathfrak c$ and $\mathfrak a$ are $\mathbb F_p$-linear spaces; cf. \cite[Cor. 7.7.3]{W97}.
 
 \

  Consider an exact sequence of Lie $\mathbb Z$-algebras
 $$0  \rightarrow  \mathfrak a  \rightarrow  \mathfrak b  \rightarrow  \mathfrak c  \rightarrow  0
 $$
 with $\mathfrak a$ abelian such that $\mathfrak b  \rightarrow  \mathfrak c$  has a section $\sigma$ in the category of  $\mathbb Z$-modules.
 Then $\mathfrak a$ has a natural induced structure of Lie  $\mathfrak c$-module:
 for $x\in \mathfrak c$ and $y\in \mathfrak a$ one defines $x\cdot y=[\tilde{x},y]$ where
 $\tilde{x}\in \mathfrak b$ is a lift of $x$.
 We may attach to our extension a class in the group $H^2_{\textup{Lie}}(\mathfrak c,\mathfrak a)$ by  considering the class of the Lie $2$-cocycle
 $$(x,y)\mapsto [\sigma(x),\sigma(y)]-\sigma([x,y]).$$
 This class does not depend on the choice of $\sigma$. We say that  the exact sequence is {\bf split} in the category of  Lie $\mathbb Z$-algebras if $\mathfrak b  \rightarrow  \mathfrak c$  has a section in that category. Then the sequence is split in the category of  Lie $\mathbb Z$-algebras if and only if its class in $H^2_{\textup{Lie}}(\mathfrak c,\mathfrak a)$ is trivial.

\

 We sometimes consider exact sequences as above in the category of graded Lie $\mathbb Z$-algebras;
 such a sequence is called {\bf split} in the category of graded Lie $\mathbb Z$-algebras
 if $\mathfrak b  \rightarrow  \mathfrak c$  has a section in that category.

\

 \subsubsection{Hochschild cohomology}.
Let $\mathfrak c$ be an associative $\mathbb Z$-algebra and let $\mathfrak a$ be a $\mathfrak c$-$\mathfrak c$-bimodule.
We denote by $Z^2_{\textup{Hoch}}(\mathfrak c,\mathfrak a)$ the group of all Hochschild $2$-cocycles $f: \mathfrak c\times \mathfrak c\rightarrow \mathfrak a$, i.e.,  $\mathbb Z$-bilinear maps such that for all $x,y,z\in \mathfrak c$ we have
$$x f(y,z)-f(xy,z)+f(x,yz)-f(x,y)z=0.
$$
 We denote by $B^2_{\textup{Hoch}}(\mathfrak c,\mathfrak a)$ the group of all Hochschild $2$-coboundaries $f:\mathfrak c\times \mathfrak c\rightarrow \mathfrak a$, i.e., maps of the form
 $$f(x,y)=xg(y)- g(xy)+g(x)y,\ \ x,y\in \mathfrak c,$$
 where $g:\mathfrak c\rightarrow \mathfrak a$ is a $\mathbb Z$-module homomorphism.
  We define the {\bf $2$nd Hochschild  cohomology group}
$$H^2_{\textup{Hoch}}(\mathfrak c,\mathfrak a):=Z^2_{\textup{Hoch}}(\mathfrak c,\mathfrak a)/
B^2_{\textup{Hoch}}(\mathfrak c,\mathfrak a).$$

\

  Consider an exact sequence of associative $\mathbb Z$-algebras
 $$0  \rightarrow  \mathfrak a  \rightarrow  \mathfrak b  \rightarrow  \mathfrak c  \rightarrow  0
 $$
 with $\mathfrak a$ an ideal in $\mathfrak b$ of square zero such that $\mathfrak b  \rightarrow  \mathfrak c$  has a section $\sigma$ in the category of  $\mathbb Z$-modules.
 Then $\mathfrak a$ has a natural induced structure of   $\mathfrak c$-$\mathfrak c$-bimodule:
 for $x\in \mathfrak c$ and $y\in \mathfrak a$ one defines $x\cdot y=\tilde{x}y$ and 
 $y\cdot x=y\tilde{x}$
 where
 $\tilde{x}\in \mathfrak b$ is a lift of $x$.
  
  \

 We may attach to our extension a class in the group $H^2_{\textup{Hoch}}(\mathfrak c,\mathfrak a)$ by  considering the class of the Hochschild $2$-cocycle
 $$(x,y)\mapsto \sigma(x)\sigma(y)-\sigma(xy).$$
 This class does not depend on the choice of $\sigma$. We say that  the exact sequence is {\bf split} in the category of  associative $\mathbb Z$-algebras if $\mathfrak b  \rightarrow  \mathfrak c$  has a section in that category. Then the sequence is split in the category of  associative $\mathbb Z$-algebras if and only if its class in $H^2_{\textup{Hoch}}(\mathfrak c,\mathfrak a)$ is trivial. 
 
\

 We sometimes consider exact sequences as above in the category of graded associative $\mathbb Z$-algebras;
 such a sequence is called {\bf split} in the category of associative $\mathbb Z$-algebras
 if $\mathfrak b  \rightarrow  \mathfrak c$  has a section in that category.
 
 \
 
 Note that if $\mathfrak c$ is an associative $\mathfrak a$-algebra and $\mathfrak a$ is a $\mathfrak c$-$\mathfrak c$-bimodule then, viewing $\mathfrak a$ as a Lie module over the Lie $\mathbb Z$-algebra $\mathfrak c$ we have a natural group homomorphism
 \begin{equation}
 \label{hochtolie}
 H^2_{\textup{Hoch}}(\mathfrak c,\mathfrak a)\rightarrow H^2_{\textup{Lie}}(\mathfrak c,\mathfrak a)\end{equation}
 defined by sending the class of any Hochschild cocycle $f:\mathfrak c\times \mathfrak c\rightarrow a$ into the class of the Lie cocycle $f^{\textup{Lie}}:\mathfrak c\times \mathfrak c\rightarrow \mathfrak a$ given by
 $$f^{\textup{Lie}}(x,y):=f(x,y)-f(y,x).$$
 
 \

 Finally note that in the classical treatment of Hochschild cohomology (e.g. in \cite[Sect. 9]{W97})
 all associative algebras are assumed to have a unit (and all bimodules are unitary) whereas in our setting we do not assume our associative $\mathbb Z$-algebras have  a unit. If we assume $\mathfrak c$ has a unit (and $\mathfrak a$ is  unitary)   then the  Hochschild cohomology groups $H^2_{\textup{Hoch}}(\mathfrak c,\mathfrak a)$ defined above coincide with the Hochschild cohomology groups $H^2(\mathfrak c,\mathfrak a)$ in \cite[Sect. 9]{W97}.  On the other hand if $\mathfrak c$ is an arbitrary (not necessarily unital)  associative $\mathbb Z$-algebra  (which is the case in our applications) and $\mathfrak a$ is a $\mathfrak c$-$\mathfrak c$-bimodule one can consider the  associative algebra  $\mathfrak c^+:=\mathbb Z\oplus \mathfrak c$
 in which 
 $$(n,x)(m,y)=(nm,ny+mx+xy)$$
 for $(n,x),(m,y)\in \mathfrak c^+$
 and we may consider $\mathfrak a$ as a $\mathfrak c^+$-$\mathfrak c^+$-bimodule
 via the rules
 $$(n,x) a=na+xa,\ \ a(n,x)=na+ax$$
 for all $a\in \mathfrak a$. 
 Then $\mathfrak c^+$ has unit $(1,0)$ (which is, of course,  different from the unit of $\mathfrak c$ in case $\mathfrak c$ has a unit). With these definitions
 there is a natural  group isomorphism
 $$H^2_{\textup{Hoch}}(\mathfrak c,\mathfrak a)\simeq H^2_{\textup{Hoch}}(\mathfrak c^+,\mathfrak a)$$
 defined by sending the class of any Hochschild cocycle $f:\mathfrak c\times \mathfrak c\rightarrow a$
 into the class of the Hochschild  cocycle $f^+:\mathfrak c^+\times \mathfrak c^+\rightarrow a$
 given by
 $$f^+((n,x),(m,y)):=f(x,y).$$
 
\

\subsection{Cohomology classes of Weil monoids} \label{ccwm}
Fix an integer $c\in \mathbb N$ and let 
$\mathbb M_{n,c}^+$ be the monoid
 $\mathbb M_n^+$ viewed as a  graded monoid with components of degree $s\in c\mathbb N$ consisting of the words of length $s/c$ and all the components of degree $s\in \mathbb N\backslash c\mathbb N$ empty. Note that $\mathbb M_{n,c}^+$   is a Weil monoid if and only if $n=1$.
 
 \

 Let $\mathfrak D$ be a Weil monoid of minimal degree $c$ and let 
  $\mathfrak D^{(c)}=\{\phi_1^{(c)},\ldots,\phi_n^{(c)}\}$. Consider the unique graded monoid homomorphism $\mathbb M_{n,c}^+\rightarrow \mathfrak D$ that sends every $i\in \{1,\ldots,n\}\subset \mathbb M_n^+$ into $\phi_i^{(c)}$.
   Denote by $\mathfrak d:=\mathbb Z[\mathfrak D]$ the monoid $\mathbb Z$-algebra of $\mathfrak D$ and by $\mathfrak l_{n,c}=\mathbb Z[\mathbb M_{n,c}^+]$  the  monoid algebra attached to $\mathbb M_{n,c}$; we view these as graded associative $\mathbb Z$-algebras 
   with gradings induced by those of the corresponding monoids in the natural way.
      Then we have an induced surjective graded associative $\mathbb Z$-algebra  homomorphism $\mathfrak l_{n,c}\rightarrow \mathfrak d$ whose kernel we denote by
    $\mathfrak n_{\mathfrak D}$.    
    
    \

   \begin{definition} For a Weil monoid $\mathfrak D$
   the {\bf Lie class} and the {\bf Hochschild class}    are the classes
    \begin{equation}
   \label{maddr}
   \kappa_{\mathfrak D}\in H^2_{\textup{Lie}}(\mathfrak d,\mathfrak n_{\mathfrak D}/\mathfrak n_{\mathfrak D}^2),\ \ \ 
   h_{\mathfrak D}\in H^2_{\textup{Hoch}}(\mathfrak d,\mathfrak n_{\mathfrak D}/\mathfrak n_{\mathfrak D}^2)\end{equation}
    attached to the exact sequence
   \begin{equation}
   \label{connff}
   0  \rightarrow  \mathfrak n_{\mathfrak D}/\mathfrak n_{\mathfrak D}^2  \rightarrow  \mathfrak l_{n,c}/\mathfrak n_{\mathfrak D}^2 
   \rightarrow  \mathfrak d  \rightarrow  0.\end{equation}
    \end{definition}

     The above definition makes sense because  (\ref{connff})  is  an extension of Lie (respectively associative) $\mathbb Z$-algebras which is split in the category of $\mathbb Z$-modules and $\mathfrak n_{\mathfrak D}/\mathfrak n_{\mathfrak D}^2$ is an ideal of square zero (respectively an abelian ideal) in the associative (respectively Lie) $\mathbb Z$-algebra $\mathfrak l_{n,c}/\mathfrak n_{\mathfrak D}^2$.

    \

    The groups and the classes in Equation (\ref{maddr})  are objects of a purely group  theoretic nature and it would be interesting to have an ``explicit" description of these objects; see Example \ref{shegun} for a discussion of some special cases. Note that $\kappa_{\mathfrak D}$ is the image
   of $h_{\mathfrak D}$ under the natural map between the cohomology groups. 
   
   \

   Assume in what follows that $\mathfrak D$ is a Weil monoid of minimal degree $c$ and size $n$. Note that $\kappa_{\mathfrak D}=0$ if $n=1$. We expect that $\kappa_{\mathfrak D}\neq 0$ (and hence $h_{\mathfrak D}\neq 0$) for all $n\geq 2$. We will check 
   that $h_{\mathfrak D}\neq 0$ for $n\geq 2$; we will also check that $\kappa_{\mathfrak D}\neq 0$ 
    in case $\mathfrak D^{(c)}$ contains two distinct commuting elements; see Proposition \ref{knz} below. The condition on the two commuting elements can always be achieved by replacing the monoid $\mathfrak D$ by the monoid $\bigcup_{s\in ct\mathbb N}\mathfrak D^{(s)}$ for some $t\in \{1,\ldots,n\}$; see Corollary \ref{tdce}.

\begin{proposition}\label{knz} For a Weil monoid $\mathfrak D$ of minimal degree $c$ and size $n$ the following properties hold:

1) Assume $n\geq 2$. Then the  Hochschild class $h_{\mathfrak D}$ attached to $\mathfrak D$
is non-trivial.

2) Assume $\mathfrak D^{(c)}$ contains two distinct commuting elements. Then the  Lie class $\kappa_{\mathfrak D}$ attached to $\mathfrak D$  is non-trivial.
\end{proposition}

{\it Proof}. 
Choose any map of sets $S:\mathfrak D\rightarrow \mathbb M_{n,c}^+$ which is a section of the surjection $\mathbb M_{n,c}^+\rightarrow \mathfrak D$. Consider the unique extension of $S$ to an $\mathbb Z$-linear  map (still denoted by) $S:\mathfrak d\rightarrow \mathfrak l_{n,c}$ and denote by $s:\mathfrak d\rightarrow \mathfrak l_{n,c}/\mathfrak n_{\mathfrak D}^2$ the induced $\mathbb Z$-module homomorphism.

\

We first prove assertion 2; assume  $\kappa_{\mathfrak D}$ is trivial and
 we will derive a contradiction.
Assume $X,Y\in \mathfrak D^{(c)}$ are distinct and commuting. Since $\kappa_{\mathfrak D}$ is trivial there exists a $\mathbb Z$-module homomorphism $g:\mathfrak d\rightarrow \mathfrak n_{\mathfrak D}/\mathfrak n_{\mathfrak D}^2$ such that 
$$s(X)s(Y)-s(Y)s(X)-s(XY)+s(YX)=X g(Y)-Yg(X)-g(XY)+g(YX).$$
Since $XY=YX$ we have
$$s(X)s(Y)-s(Y)s(X)=X g(Y)-Yg(X).$$
Note that the component $\mathfrak n_{\mathfrak D}^{(c)}$ of $\mathfrak n_{\mathfrak D}$ of degree $c$ vanishes because the component of degree $c$ of $\mathfrak l_{n,c}$ maps isomorphically onto $\mathfrak d^{(c)}$; hence the component of degree $c$ of $\mathfrak n_{\mathfrak D}/\mathfrak n_{\mathfrak D}^2$ vanishes. So the homogeneous components of
$X g(Y)-Yg(X)$  have degree $\geq 3c$. Since $s(X)s(Y)-s(Y)s(X)$ is homogeneous of degree $2c$ we get that $s(X)s(Y)-s(Y)s(X)=0$. Hence 
$S(X)S(Y)-S(Y)S(X)\in \mathfrak n_{\mathfrak D}^2$. Since $S(X)S(Y)-S(Y)S(X)$ is homogeneous of degree $2c$ and since all elements in $\mathfrak n_{\mathfrak D}^2$ have homogeneous components of degree $\geq 4c$ it follows that $S(X)S(Y)=S(Y)S(X)$  in $\mathbb M_n$ which is a contradiction because $S(X)$ and $S(Y)$ are two distinct elements
in the set $\{1,\ldots,n\}$ that freely generates  $\mathbb M_n$. This proves assertion 2.

\

We next prove assertion 1; assume  $h_{\mathfrak D}$ is trivial and we will derive a contradiction.
Since $n\geq 2$ the map
$$\mathfrak D^{(c)}\times \mathfrak D^{(c)}\rightarrow \mathfrak D^{(2c)},\ 
\ \ (X,Y)\mapsto XY$$
is not injective hence one can find two distinct pairs $(X_1,Y_1),(X_2,Y_2)\in \mathfrak D^{(c)}\times \mathfrak D^{(c)}$ such that $X_1Y_1=X_2Y_2$. In particular $X_1\neq X_2$ hence $S(X_1)\neq S(X_2)$. Since $h_{\mathfrak D}$ is trivial there exists a $\mathbb Z$-module homomorphism $g:\mathfrak d\rightarrow \mathfrak n_{\mathfrak D}/\mathfrak n_{\mathfrak D}^2$ such that 
$$s(X_1)s(Y_1)-s(X_1Y_1)= X_1g(Y_1)-g(X_1Y_1)+g(X_1)Y_1,$$
$$s(X_2)s(Y_2)-s(X_2Y_2)= X_2g(Y_2)-g(X_2Y_2)+g(X_2)Y_2.$$
Subtracting these equalities and using $X_1Y_1=X_2 Y_2$ we get
$$s(X_1)s(Y_1)-s(X_2)s(Y_2)=X_1g(Y_1)+g(X_1)Y_1-X_2g(Y_2)-g(X_2)Y_2.$$

\

Again, since the component of degree $c$ of $\mathfrak n_{\mathfrak D}/\mathfrak n_{\mathfrak D}^2$ vanishes, the homogeneous components of
$X_1g(Y_1)$, $g(X_1)Y_1$, $X_2g(Y_2)$, $g(X_2)Y_2$  have degree $\geq 3c$. 
Since $s(X_1)s(Y_1)-s(X_2)s(Y_2)$ is homogeneous of degree $2c$ we get that $s(X_1)s(Y_1)-s(X_2)s(Y_2)=0$. Hence 
$S(X_1)S(Y_1)-S(X_2)S(Y_2)\in \mathfrak n_{\mathfrak D}^2$. Since $S(X_1)S(Y_1)-S(X_2)S(Y_2)$ is homogeneous of degree $2c$ and since all elements in $\mathfrak n_{\mathfrak D}^2$ have homogeneous components of degree $\geq 4c$ it follows that $S(X_1)S(Y_1)=S(X_2)S(Y_2)$  in $\mathbb M_n$ which is a contradiction because $S(X_1)$ and $S(X_2)$ are two distinct elements
in the set $\{1,\ldots,n\}$ that freely generates  $\mathbb M_n$. This proves assertion 1.
\qed

\begin{remark}\label{ofofof}
There is an obvious way to express the algebra 
$\mathfrak l_{n,c}=\mathbb Z[\mathbb M_{n,c}^+]$ 
 above and as an ideal in the algebra of non-commutative polynomials. Indeed consider the associative $\mathbb Z$-algebra of {\bf  non-commutative polynomials}
 \begin{equation}
 \label{rrt}
\mathbb Z\langle T_1,\ldots,T_n\rangle\end{equation}
  whose elements are formal sums 
  $$\lambda_0+\sum_{\mu\in \mathbb M_n}\lambda_{\mu} T_{\mu}$$ with $\lambda_0,\lambda_{\mu}\in \mathbb Z$  and $T_{\mu}:=T_{i_1}\cdots T_{i_m}$ for $\mu=i_1\cdots i_m$. 

  \

Denote 
by $\mathfrak p_n=(T_1,\ldots,T_n)$ the  ideal  in (\ref{rrt}) generated by $T_1,\ldots,T_n$. We view $\mathfrak p_n$ 
 as a graded associative $\mathbb Z$-algebra with  component of degree $s\in \mathbb N$ 
 the $\mathbb Z$-submodule generated by the monomials of degree $s/c$ if $c$ divides $s$ and $0$ otherwise.
Then we have a graded associative $\mathbb Z$-algebra isomorphism 
 \begin{equation}
 \label{ffuutt}
 \mathfrak l_{n,c}\simeq \mathfrak p_n,\ \ \mu\mapsto T_{\mu}.\end{equation} 

 \

Now, let $\mathfrak D$ be a Weil monoid of minimal degree $c$ and size $n$. Under the 
identification (\ref{ffuutt}),
 the ideal $\mathfrak n_{\mathfrak D}$ is easily seen to be generated by the family of 
 $n(n-1)$ 
 non-commutative quadratic polynomials
\begin{equation}
\label{techno}
(T_iT_j-T_1T_{k_{ij}})_{(i,j)\in \{2,\ldots,n\}\times\{1,\ldots,n\}},\end{equation}
where for each $i,j$ as above $k_{ij}\in \{1,\ldots,n\}$ is the unique index such that
$$\phi_i^{(c)}\phi_j^{(c)}=\phi_1^{(c)}\phi^{(c)}_{k_{ij}}.$$
The component $\mathfrak n_{\mathfrak D}^{(c)}$ of degree $c$ 
in $\mathfrak n_{\mathfrak D}$ vanishes.
The component $\mathfrak n_{\mathfrak D}^{(2c)}$ of degree $2c$ in $\mathfrak n_{\mathfrak D}$ (which identifies with the  component $(\mathfrak n_{\mathfrak D}/\mathfrak n_{\mathfrak D}^2)^{(2c)}$ of degree $2c$ in $\mathfrak n_{\mathfrak D}/\mathfrak n_{\mathfrak D}^2$) is a free $\mathbb Z$-module with basis the family (\ref{techno}). \end{remark}

\begin{example} \label{shegun}
Let $\mathfrak D$ be a Weil monoid of minimal degree $c$ and size $n$.
Consider the labeling symbols $\omega^{(c)}$ and $\omega^{(2c)}$  defined by the equalities 
$$\mathfrak D^{(c)}=\{\phi_1^{(c)},\ldots,\phi_n^{(c)}\},\ \ \ \mathfrak D^{(2c)}=\{\phi_1^{(2c)},\ldots,\phi_n^{(2c)}\},\ \ \ \phi_i^{(2c)}=\phi_1^{(c)}\phi_i^{(c)}.$$ 
(So $\omega^{(2c)}$ is the canonical labeling symbol attached to $\omega^{(c)}$ and $h=1$.)
Write $\alpha^k=\alpha^{k(c)}$, $\ell^k=\ell^{k(c)}$ for the associative and Lie symbols
corresponding to the above labeling symbols.
Furthermore consider an embedding of $\mathfrak D$  in a group $\mathfrak H$. Then by Lemma \ref{prelop} there exists a subgroup $\mathfrak S=\{\sigma_1,\ldots,\sigma_n\}$ of $\mathfrak H$ with $\sigma_1=1$ such that $\phi_1^{(c)}$ normalizes $\mathfrak S$ and such that $\phi_i^{(c)}=\phi_1^{(c)}\sigma_i$ for $i\in \{1,\ldots,n\}$.

\

In what follows we illustrate our concepts in case $n=2$ and $n=3$.

\

1) If $n=2$ then $\mathfrak S=C_2$, a cyclic group of order $2$, and $\mathfrak D$ 
 is represented by the pair $(C_2,\textup{id})$. We have
\begin{equation}
\label{you0}
\alpha^1=\left(\begin{array}{cc} 1 & 0\\ 0 & 1\end{array}\right),\ \  \alpha^2=\left(\begin{array}{cc} 0 & 1\\ 1 & 0\end{array}\right),\ \ \ \ell^1=\ell^2=0.\end{equation} 
By Remark \ref{ofofof}   the ideal $\mathfrak n_{\mathfrak D}$ is generated in this case by the $2$ quadratic polynomials
$$T_2T_1-T_1T_2,\ \ \ T_2^2-T_1^2,$$
which also  form a $\mathbb Z$-module basis for the component $\mathfrak n_{\mathfrak D}^{(2c)}$.
To go one step further the component $\mathfrak n_{\mathfrak D}^{(3c)}$ of degree $3c$ in $\mathfrak n_{\mathfrak D}$ (which identifies with the component $(\mathfrak n_{\mathfrak D}/\mathfrak n_{\mathfrak D}^2)^{(3c)}$ of degree $3c$ in $\mathfrak n_{\mathfrak D}/\mathfrak n_{\mathfrak D}^2$) has a $\mathbb Z$-module basis consisting of 
the $6$ cubic polynomials
$$T_1^3-T_1T_2^2,\ \ T_2^3-T_2T_1^2,$$
 $$T_2^2T_1-T_1T_2^2,\ \ T_1^2T_2-T_2T_1^2,$$ 
$$T_1T_2T_1-T_1^2T_2,\ \ T_2T_1T_2-T_2^2T_1.$$

\

2) If $n=3$ and  $\phi_1^{(c)}$ centralizes $\mathfrak S$ then  $\mathfrak D$ is  represented by $(C_3,\textup{id})$ where $\mathfrak S=C_3$ is a cyclic group of order $3$ and we have
\begin{equation}
\label{you1}
\alpha^1=\left(\begin{array}{ccc} 1 & 0 & 0\\ 0 & 0 & 1\\ 0 & 1 & 0\end{array}\right),\ \  \alpha^2=\left(\begin{array}{ccc} 0 & 1 & 0\\ 1 & 0 & 0\\ 0 & 0 & 1\end{array}\right),\ \ 
\alpha^3=\left(\begin{array}{ccc} 0 & 0 & 1\\ 0 & 1 & 0\\ 1 & 0 & 0\end{array}\right),\end{equation}
hence
\begin{equation}
\label{you11}
\ell^1=\ell^2=\ell^3=0.\end{equation}
In this case
  the ideal $\mathfrak n_{\mathfrak D}$ is generated  by the $6$ quadratic polynomials
$$T_2T_1-T_1T_2,\ \ \ T_2T_3-T_3T_2,\ \ \ T_3T_1-T_1T_3,$$
$$T_2^2-T_1T_3,\ \ \ T_1^2-T_2T_3,\ \ \ T_3^2-T_1T_2,$$
which also form a $\mathbb Z$-module basis of $\mathfrak n_{\mathfrak D}^{(2c)}$.

\

3) If $\phi_1^{(c)}$ does not centralize $\mathfrak S$  then  $\mathfrak D$ is  represented by $(C_3,\theta)$ with $\theta\neq \textup{id}$. In this case we have
\begin{equation}
\label{you2}
\alpha^1=\left(\begin{array}{ccc} 1 & 0 & 0\\ 0 & 1 & 0\\ 0 & 0 & 1\end{array}\right),\ \  \alpha^2=\left(\begin{array}{ccc} 0 & 1 & 0\\ 0 & 0 & 1\\ 1 & 0 & 0\end{array}\right),\ \ 
\alpha^3=\left(\begin{array}{ccc} 0 & 0 & 1\\ 1 & 0 & 0\\ 0 & 1 & 0\end{array}\right)\end{equation}
hence 
\begin{equation}
\label{you22}
\ell^1=0,\ \ \ell^2=-\ell^3=\left(\begin{array}{rrr} 0 & 1 & -1\\ -1 & 0 & 1\\ 1 & -1 & 0\end{array}\right).
\end{equation}
Furthermore, in this case the ideal $\mathfrak n_{\mathfrak D}$ is generated  by the $6$ quadratic polynomials
$$T_2T_1-T_1T_3,\ \ \ T_2^2-T_1^2,\ \ \ T_2T_3-T_1T_2,$$
$$T_3T_1-T_1T_2,\ \ \ T_3T_2-T_1T_3,\ \ \ T_3^2-T_1^2,$$
which also form a $\mathbb Z$-module basis of $\mathfrak n_{\mathfrak D}^{(2c)}$.
Note that in this case every two distinct elements of $\mathfrak D^{(c)}$ are non-commuting. 

\

The explicit determination of all the homogeneous components of $\mathfrak n_{\mathfrak D}/\mathfrak n_{\mathfrak D}^2$ and of the corresponding cohomology groups and classes $h_{\mathfrak D}$ in the above example  seems to be an intricate combinatorics/linear algebra problem that we will not pursue here.\end{example}

  \section{Curvature and characteristic classes in the  PDE setting}

\subsection{Canonical  algebra} \label{cannn}
So far our discussion after the Introduction did not involve local fields. We now return to the local field setting in Subsection \ref{generall}. In this subsection we consider a canonical associative (and hence also Lie) $\mathbb Z$-algebra attached
 to $K_{\pi}$ and we implement the monoid abstract nonsense in Subsection \ref{labelll} to this concrete setting. 

\

Recall that  the extension $K_{\pi}/K$ is Galois and consider its Galois group  $\mathfrak G(K_{\pi}/K)$.
 For every integer $s\geq 1$ let 
 $\mathfrak F^{(s)}:=\mathfrak F^{(s)}(K_{\pi}/\mathbb Q_p)$ be the set of 
 higher Frobenius automorphisms of $K_{\pi}$ of degree $s$.  It is clear that  
 $\mathfrak F^{(s)}$ is a principal homogeneous space for $\mathfrak G(K_{\pi}/K)$ under the action
 $(\phi^{(s)},\sigma)\mapsto \phi^{(s)}\sigma$ for $\phi^{(s)}\in \mathfrak F^{(s)}$ and $\sigma \in \mathfrak G(K_{\pi}/K)$. By the ``linear independence of characters"  the union $\mathfrak F:=\bigcup_{s=1}^{\infty} \mathfrak F^{(s)}$ 
 is a disjoint union and 
 is a $K_{\pi}$-linearly independent subset of the $K_{\pi}$-linear space $\textup{End}_{\mathbb Z-\textup{mod}}(K_{\pi})$ of $\mathbb Z$-module endomorphisms of $K_{\pi}$. 
 We have  $\mathfrak F^{(s)} \mathfrak F^{(r)}\subset
 \mathfrak F^{(s+r)}$ in the associative $\mathbb Z$-algebra $\textup{End}_{\mathbb Z-\textup{mod}}(K_{\pi})$
 hence $\mathfrak F$ is a monoid which is clearly a Weil monoid with surjective degree map $\mathfrak F\rightarrow \mathbb N$  sending $\mathfrak F^{(s)}$ into $s$; cf Definition \ref{rgm}. (This is precisely the context in which Weil monoids appear in \cite[p. 69]{N80}.)

 \

 The $\mathbb Z$-linear span $\mathfrak f$ of $\mathfrak F$ in $\textup{End}_{\mathbb Z-\textup{mod}}(K_{\pi})$ is a graded associative $\mathbb Z$-subalgebra of $\textup{End}_{\mathbb Z-\textup{mod}}(K_{\pi})$ graded by the additive monoid $\mathbb N$ and whose component of degree $s$ is the $\mathbb Z$-linear span 
 $\mathfrak f^{(s)}$
 of $\mathfrak F^{(s)}$. Note that $\mathfrak f^{(s)}$ is a free $\mathbb Z$-module with basis $\mathfrak F^{(s)}$ and $\mathfrak f$ is a free $\mathbb Z$-module with basis $\mathfrak F$. 
 So $\mathfrak f$ is isomorphic to the   monoid algebra 
 $\mathbb Z[\mathfrak F]$
 attached to the monoid $\mathfrak F$.    We may  view $\textup{End}_{\mathbb Z-\textup{mod}}(R_{\pi})$ as a Lie $\mathbb Z$-algebra with respect to the commutator operator. So
  $\mathfrak f$ can (and will) be identified with a  Lie $\mathbb Z$-subalgebra of $\textup{End}_{\mathbb Z-\textup{mod}}(R_{\pi})$.
  This is simply because every $\mathbb Z$-module endomorphism of $K_{\pi}$ which is zero on $R_{\pi}$ must be zero.

 \begin{definition} The associative (and hence also Lie)  $\mathbb Z$-algebra  $\mathfrak f$ is called the {\bf canonical  algebra} of $R_{\pi}$.
 \end{definition}
 
 Let $\phi_K$ be the unique  element of $\mathfrak F^{(1)}(K/\mathbb Q_p)$. 
  
 \begin{proposition}
 The $K_{\pi}$-linear span of $\mathfrak F^{(s)}$ 
 in $\textup{End}_{\mathbb Z-\textup{mod}}(K_{\pi})$
 coincides with the  $K_{\pi}$-linear subspace,
 $\textup{End}_{\phi_K^s}(K_{\pi})$,
  of $\textup{End}_{\mathbb Z-\textup{mod}}(K_{\pi})$
 consisting of all $\phi_K^s$-linear  endomorphisms
 i.e., endomorphisms $\alpha$ such that $\alpha(cx)=\phi_K^s(c)\alpha(x)$ for $c\in K$ and $x\in K_{\pi}$.
 \end{proposition}
 
 {\it Proof}.
 Pick any element $\phi^{(s)}_0\in \mathfrak F^{(s)}$. Then the map
  $\alpha \mapsto (\phi^{(s)}_0)^{-1}\circ \alpha$ is a bijection from $\textup{End}_{\phi_K^s}(K_{\pi})$
  to the $K_{\pi}$-linear space, $\textup{End}_{K-\textup{mod}}(K_{\pi})$, of all $K$-linear endomorphisms of $K_{\pi}$. This bijection preserves $K_{\pi}$-linear independence. On the other hand
  $\textup{End}_{K-\textup{mod}}(K_{\pi})$ has dimension $n^2$ over $K$ for $n:=[K_{\pi}:K]$ hence it has dimension $n$ over $K_{\pi}$. We conclude by noting that the $K_{\pi}$-linear span of 
  $\mathfrak F^{(s)}$ is contained in $\textup{End}_{\phi_K^s}(K_{\pi})$ and has dimension $n$ over $K_{\pi}$. 
 \qed
 
 \begin{remark}
 The above statement cannot be refined into an ``integral" statement in the sense that for $\phi_R$ the Frobenius lift on $R$ 
 not every $\phi_R^s$-linear endomorphism of $R_{\pi}$ is an $R_{\pi}$-linear combination of restrictions to $R_{\pi}$ of the elements of $\mathfrak F^{(s)}$. Indeed if $\phi_1^{(s)}$ and $\phi_2^{(s)}$ are two distinct elements of $\mathfrak F^{(s)}$ then the linear combination
 $$\frac{1}{\pi}\phi_1^{(s)}-\frac{1}{\pi}\phi_2^{(s)}:R_{\pi}\rightarrow R_{\pi}$$
 is a counterexample.
 \end{remark}

 It is useful to consider certain Lie subalgebras 
 of the Lie $\mathbb Z$-algebra 
 $\mathfrak f$ as follows. First for $c\in \mathbb N$ and for a subset  $\mathfrak D^{(c)}\subset \mathfrak F^{(c)}$ we define a graded monoid inductively by $\mathfrak D^{(s)}:=\mathfrak D^{(s-c)}\mathfrak D^{(c)}$ if $s\in c\mathbb N\setminus \{c\}$ and setting $\mathfrak D = \bigcup_{s \in c \mathbb N} D^{(s)}$. We call $\mathfrak D$ the {\bf monoid associated to $\mathfrak D^{(c)}$}. We now define two important properties. 
 
 \begin{definition}\label{rregular} Fix $\mathfrak D^{(c)} \subset \mathfrak F^{(c)}$. 
 \begin{enumerate}
\item We say $\mathfrak D^{(c)}$ is {\bf abelian} if its elements are pairwise commuting in $\mathfrak D$, i.e., that $\mathfrak D$ is an abelian monoid. 
\item We say $\mathfrak D^{(c)}$ is a {\bf involutive} if $\mathfrak D^{(s)}$ has the same cardinality as $\mathfrak D^{(c)}$ for all $s\in c\mathbb N$, i.e., that $\mathfrak D$ is a Weil monoid. 
\end{enumerate} 
\end{definition}

 The terminology in the definition above is motivated by the analogy, which will become apparent, between involutive subsets of $\mathfrak F^{(c)}$ in the sense above and (Frobenius) involutive distributions in classical differential geometry. 
 
 \begin{remark}\

 \
 
 1) $\mathfrak F^{(c)}$ is involutive for all $c\in \mathbb N$.

 \

 2) If $\mathfrak D^{(c_0)}\subset \mathfrak F^{(c_0)}$ is an arbitrary subset then there exists an integer
 $c\in c_0\mathbb N$ such that $\mathfrak D^{(c)}$ is involutive.
 \end{remark}

\begin{remark}\label{ccoonn}
 Assume $\mathfrak D^{(c)}\subset \mathfrak F^{(c)}$ is involutive, so $\mathfrak D$ has a natural structure of Weil monoid with degree map surjecting onto $c\mathbb N$ (cf. Definition \ref{rgm}). In particular, as in Definition \ref{ccaannabstract}, one can consider 
 labeling symbols
 $$\omega^{(s)}:\{1,\ldots,n\}\rightarrow \mathfrak D^{(s)},\ \ \omega^{(s)}(i)=:\phi_i^{(s)};$$
 note that the tuple $(\phi_1^{(s)},\ldots,\phi_n^{(s)})$ defines a structure of partial $\delta$-ring 
   of degree $s$ on $R_{\pi}$.
  According to our conventions in Part 1 if $\phi_{ij}^{(2c)}:=\phi_i^{(c)}\phi_j^{(c)}$ and if $\delta_{ij}^{(2c)}$ is the higher $\pi$-derivation of degree $2c$ attached to the higher $\pi$-Frobenius lift $\phi_{ij}^{(2c)}$ of degree $2c$ then
  for every  $a\in R_{\pi}$  we have
 \begin{equation}\label{staru}
 \delta_{i\star j}^{(2c)}a=\frac{\phi^{(2c)}_{i\star j}a-a^{p^{2c}}}{\pi}=\frac{\phi_i^{(c)}\phi^{(c)}_j a-a^{p^{2c}}}{\pi}=
 \frac{\phi_{ij}^{(2c)} a-a^{p^{2c}}}{\pi}=\delta_{ij}^{(2c)}a.
 \end{equation}
 Furthermore, as in Definition \ref{ccaannabstract}, one can consider, for $s\in c\mathbb N$,  the associative symbols
 $$\alpha^{(s)}:=(\alpha^{1(s)},\ldots,\alpha^{n(s)})$$
 and the Lie symbols 
 $$\ell^{(s)}=(\ell^{1(s)},\ldots,\ell^{n(s)})$$
  attached to our labeling symbols.
  We can write
  \begin{equation}
  \phi_i^{(s)}\phi_j^{(r)}=\sum_k \alpha^{k(s,r)}_{ji} \phi_k^{(s+r)},\end{equation}
 \begin{equation}
 \label{nygov}
 [ \phi_i^{(s)},\phi_j^{(r)}]:=
 \phi_i^{(s)}\phi_j^{(r)}-\phi_j^{(r)}\phi_i^{(s)}=\phi_{(i\star j)_{s,r}}^{(s+r)}-\phi_{(j\star i)_{r,s}}^{(s+r)}=
 \sum_k \ell_{ij}^{k(s,r)} \phi_k^{(s+r)}.\end{equation}
  \end{remark}
  
  \begin{remark}
  Given a graded monoid $\mathfrak D$  it is an interesting Galois theoretic problem to determine for which $\pi$ there exists an embedding (i.e. injective graded monoid homomorphism) $\mathfrak D\rightarrow \mathfrak F$, where $\mathfrak F$ is the monoid of higher $\pi$-Frobenius lifts on $K_{\pi}$.   \end{remark}
  
 \begin{definition}\label{ccaann}
 Assume $\mathfrak D^{(c)}\subset \mathfrak F^{(c)}$ is involutive.
 The $\mathbb Z$-linear span in $\mathfrak f$ 
 of the set $\mathfrak D$ is denoted by 
$\mathfrak d$. It is an associative and hence Lie $\mathbb Z$-subalgebra  referred to as the {\bf canonical algebra} attached to $\mathfrak D^{(c)}$.
(Note that $\mathfrak d$ is isomorphic as an associative $\mathbb Z$-algebra  to the  monoid algebra 
$\mathbb Z[\mathfrak D]$
attached to the monoid $\mathfrak D$ and our notation matches the one in Subsection \ref{ccwm}.)
The numbers $\ell_{ij}^{k(s,r)}$ in Definition \ref{ccaannabstract} are called the  {\bf structure constants} of the Lie $\mathbb Z$-algebra $\mathfrak d$. \end{definition}

\begin{remark} If $\mathfrak D^{(1)}:=\mathfrak F^{(1)}$ we have $\mathfrak d=\mathfrak f$.
The reason why it is natural to consider more general involutive subsets  $\mathfrak D^{(c)}\subset \mathfrak F^{(c)}$ is twofold. First we want to consider the case $c\neq 1$ in order to accommodate higher Frobenius automorphisms in $\mathfrak D^{(c)}$ (for instance those arising from local class field theory as in Remark \ref{lcftc} below). Second,  we want to consider metrics $q$ of dimension $n$  where 
$n$ is the number of $\pi$-derivatives 
in our $\pi$-connections  of degree $s$; but in this case  it is not natural to insist that  $n$ be
 equal to the degree $e$ of $K_{\pi}$ over $K$; the degree $e$ should be allowed to be arbitrary as we want to be able to vary the fields $K_{\pi}$ and the matrices $q$ while keeping $n$ fixed.
\end{remark}

We can completely characterize involutive subsets via the following theorem.

\begin{theorem}\label{boooo} The following statements hold:

1) 
For every  pair 
 $(\mathfrak S,\phi)$ consisting of a subgroup
  $\mathfrak S=\{\sigma_1,\ldots,\sigma_n\}$ of $\mathfrak G(K_{\pi}/K)$  and an element 
 $\phi\in \mathfrak F^{(c)}$ that normalizes $\mathfrak S$ (i.e., $\phi^{-1}\mathfrak S\phi=\mathfrak S$) the subset 
 $$\mathfrak D^{(c)}:=\{\phi^{s/c}\sigma_1,\ldots,\phi^{s/c}\sigma_n\}\subset \mathfrak F^{(c)}$$ is involutive and hence $\mathfrak D:=\bigcup_{s\in c\mathbb N} \mathfrak D^{(s)}$ is a Weil monoid.
 
 2) Every involutive subset $\mathfrak D^{(c)}$ of $\mathfrak F^{(c)}$ 
 (and hence every Weil submonoid $\mathfrak D$ of $\mathfrak F$ of minimal degree $c$)
 arises from a pair $(\mathfrak S,\theta)$ as in 1) above. \end{theorem}

 {\it Proof}.  Assertion 1 is clear. Assertion 2 follows from Lemma \ref{prelop}.
\qed

\begin{remark}\label{incaa} In the notation of Theorem \ref{boooo} if $\sigma_1=1$ and 
$\phi^{(s)}_i:=\phi^{s/c}\sigma_i$ then 
$$\mathfrak D^{(s)}=\{\phi^{(s)}_i\ |\  i\in \{1,\ldots,n\}\}$$ 
and $\omega^{(s)}(i)=\phi_i^{(s)}$ define natural labeling symbols. 
 Of course,  the Weil monoid $\mathfrak D$
 is isomorphic as a graded monoid to the semidirect product $c\mathbb N
 \times_{\textup{inn}_{\phi}} \mathfrak S$ where $\textup{inn}_{\phi}(\sigma_i)=\phi^{-1}\sigma_i\phi$. So the monoid $\mathfrak D$ arising from $(\mathfrak S,\phi)$ as in Theorem \ref{boooo} coincides with the monoid represented by $(\mathfrak S,\textup{inn}_{\phi})$
 in the sense of Definition \ref{rgm}.
 If for $i,j\in \{1,\ldots,n\}$ one defines $\phi(i),i\circ j\in \{1,\ldots,n\}$ via the the equalities $\sigma_{\phi(i)}:=\phi^{-1}\sigma_i\phi$ and $\sigma_i\sigma_j=\sigma_{i\circ j}$ 
 (see Equation (\ref{defcircc})) then, for $s,r\in c\mathbb N$, one has 
 \begin{equation}
 \label{istarjsr}
 (i\star j)_{s,r}=\phi^{r/c}(i)\circ j,\end{equation}
  hence
 $$\phi_i^{(s)}\phi_j^{(r)}=\phi_{\phi^{r/c}(i)\circ j}^{(s+r)}.$$
So  $\alpha^{k(s,r)}_{ji}=1$ if and only if $\phi^{r/c}(i)\circ j=k$. 
\end{remark}

\begin{example}\label{ogarcenusiu} Consider an involutive subset $\mathfrak D^{(c)}\subset \mathfrak F^{(c)}$ arising from a pair $(\mathfrak S,\phi)$ as in Theorem~\ref{boooo} and consider the labeling symbols in Remark \ref{incaa}. 

\

1) If  $n=2$  the associative and Lie symbols are given by (\ref{you0}).  An example of this situation is given by $c=1$, $\pi^2=p$, $\mathfrak S=\mathfrak G(K_{\pi}/K)$.

\

2) If  $n=3$   
and the action of $\phi$ on $\mathfrak S$ by conjugation is trivial (e.g. in the case $c=1$, $\pi=p^{1/3}$, $\phi \pi=\pi$, $p\equiv 1$ mod $3$,
$\mathfrak S=\mathfrak G(K_{\pi}/K)$)
 then the associative and Lie symbols are given by (\ref{you1}) and (\ref{you11}).

\

 3) If $n=3$ and  the action of $\phi$ on $\mathfrak S$ by conjugation
   is non-trivial (e.g. in the case $c=1$, $\pi=p^{1/3}$, $\phi \pi=\pi$, $p\equiv 2$ mod $3$, $\mathfrak S=\mathfrak G(K_{\pi}/K)$) then the associative and Lie symbols are given by (\ref{you2}) and (\ref{you22}).
\end{example}

\subsection{Examples from local class field theory}

A special role will be played in this paper by sets of commuting higher Frobenius automorphisms on finite extensions of $K$.
There is a systematic way to construct such sets via local class field theory. 
 
 \

Start with a finite extension $L$ of $\mathbb Q_p$ and let $p^s$ be the cardinality of the residue field of $L$. Let $L^{\textup{ur}}$ be the maximum unramified extension of $L$ in $L^{\textup{alg}}$; so $L^{\textup{ur}}$ is generated over $L$ by the roots of unity of order prime to $p$. Let $\phi_{L^{\textup{ur}}}^{(s)}$ be the unique element of $\mathfrak F^{(s)}(L^{\textup{ur}}/L)$; we denote by $\langle  \phi_{L^{\textup{ur}}}^{(s)}\rangle \simeq \mathbb Z$ the cyclic subgroup of $\mathfrak G(L^{\textup{ur}}/L)$ generated by $\phi_{L^{\textup{ur}}}^{(s)}$.
The {\bf Weil group} $W_L$ is defined as the group of all elements in $\mathfrak G(L^{\textup{alg}}/L)$ whose image in $\mathfrak G(L^{\textup{ur}}/L)$ belongs to
$\langle  \phi_{L^{\textup{ur}}}^{(s)}\rangle$; cf. \cite[p.69]{N80}. 

\

Let $W_L^{\textup{ab}}$ be the abelianization of $W_L$ and let $L^{\textup{ab}}$ be the maximum abelian extension of $L$ contained in $L^{\textup{alg}}$. There is a natural surjective homomorphism
$\deg_L:W_L\rightarrow \mathbb Z$ inducing a homomorphism still denoted by $\deg_L:W_L^{\textup{ab}}\rightarrow \mathbb Z$. The norm residue symbol in local  class field theory provides an isomorphism
$$(-, L^{\textup{ab}}/L):L^{\times}\rightarrow W_L^{\textup{ab}}$$
whose composition with $\deg_L:W_L^{\textup{ab}}\rightarrow \mathbb Z$ equals the normalized (i.e., surjective) valuation map $v_L:L^{\times}\rightarrow \mathbb Z$; cf. \cite[pp. 69-70]{N80}. 

\

  Let $\pi_{L,1},\ldots,\pi_{L,n}\in L^{\times}$ be  prime elements, i.e., 
elements of valuation  $v_L(\pi_{L,i})=1$ and let $\phi^{(s)}_{i,L^{\textup{alg}}}\in W_L$ be  automorphisms whose images in $W_L^{\textup{ab}}$ equal $(\pi_{L,i},L^{\textup{ab}}/L)$, respectively, thus $\deg(\phi^{(s)}_{i,L^{\textup{alg}}})=1$
so $\phi^{(s)}_{\pi_{L,i},L^{\textup{alg}}}\in \mathfrak F^{(s)}(L^{\textup{alg}}/L)$.  Let $E/L$ be a finite abelian extension. The 
restrictions $\phi^{(s)}_{\pi_{L,i},E}\in  \mathfrak G(E/L)$ of the automorphisms $\phi^{(s)}_{i,L^{\textup{alg}}}$ are pairwise commuting, they belong to $\mathfrak F^{(s)}(E/L)$, and they only depend on the prime elements $\pi_{L,i}$.
Let $F$ be the maximum unramified extension of $\mathbb Q_p$ contained in $E$ and $\pi_E$  a prime element in $E$ so $E=F(\pi_E)$ and we have that $K_{\pi_E}:=K(\pi_E)$ is isomorphic to the tensor product $E\otimes_F K$, in particular equals the compositum $KE$; cf. \cite[Rmk. 2.3]{BM22}. Let $\phi_K$ be the Frobenius automorphism on $K$. We may consider the higher Frobenius automorphisms
 of degree $s$ on $K_{\pi}$ defined by 
\begin{equation}
\label{xio}
\phi^{(s)}_{\pi_{L,i},K_{\pi_E}}:=\phi^{(s)}_{\pi_{L,i},E}\otimes_F \phi_K^s\in 
\mathfrak F^{(s)}(K_{\pi_E}/L).
\end{equation}

\

\noindent To summarize the above discussion in the form of the following lemma.

\begin{lemma}\label{lorna}
Let $L$ be a finite extension of $\mathbb Q_p$ with residue field of cardinality $p^s$, let  $E/L$ be a finite abelian extension, let $\pi_E$ be a prime element of $E$ and let $K_{\pi_E}:=K(\pi_E)=KE$. To every family of prime elements $\pi_{L,1},\ldots,\pi_{L,n}\in L$ one can canonically attach a family of pairwise commuting higher Frobenius automorphisms,
$$\phi^{(s)}_{\pi_{L,1},K_{\pi_E}},\ldots,\phi^{(s)}_{\pi_{L,n},K_{\pi_E}}\in  \mathfrak F^{(s)}(K_{\pi_E}/L).
$$
This attachment is compatible, in the obvious sense, with enlarging the integer $n$ and enlarging the field $E$.
\end{lemma}

\begin{remark}
In order to apply the construction above to the main setting of our paper we need to have natural
examples of situations as in the Lemma 
\ref{lorna} such that 
$E/\mathbb Q_p$ is Galois.
This is achieved via the following trivial observation: 
for every finite Galois extension  $L/\mathbb Q_p$ and every abelian extension $E_0/L$ there is an abelian extension $E/L$ containing $E_0/L$ with the property that $E/\mathbb Q_p$ is Galois. Indeed one can
take $E$ to be the compositum of all conjugates of $E_0$ over $\mathbb Q_p$. In particular for every finite Galois extension $L/\mathbb Q_p$ we have that the maximum abelian extension $L^{\textup{ab}}$ of $L$ is Galois over $\mathbb Q_p$.
\end{remark}

\begin{remark}\label{lcftc}
Here is a  natural construction of  abelian involutive subsets $\mathfrak D^{(c)}\subset \mathfrak F^{(c)}$ based on the above discussion of local class field theory. Let $E/\mathbb Q_p$ be a finite Galois extension and $\pi$  a prime element in $E$.
Let $L$ be a field between $\mathbb Q_p$ and $E$ with $E/L$  abelian, let $\pi_{L,1},\ldots,\pi_{L,n}$ be prime elements of $L$,
and let $p^{c_0}$ be the cardinality of the residue field of $L$. By Lemma \ref{lorna} one can attach to this data an abelian subset
$$\mathfrak D^{(c_0)}=\{\phi^{(c_0)}_{\pi_{L,1},K_{\pi}},\ldots,\phi^{(c_0)}_{\pi_{L,n},K_{\pi}}\}\subset  \mathfrak F^{(c_0)}(K_{\pi}/L)\subset \mathfrak F^{(c_0)}(K_{\pi}/\mathbb Q_p).$$
There exists $c\in c_0\mathbb N$ such that $\mathfrak D^{(c)}$ is involutive and automatically abelian. \end{remark}

\

   \subsection{Complements on $\pi$-connections}
 In addition to the review in Subsection \ref{generall} we need to recall and complement some of the material in Part 1. 

 \

 First recall that for every ring $S$ and every integer $N \geq 1$ we denote by $\textup{Mat}_N(S)$ the ring of $N\times N$ matrices with coefficients in $S$ and by  $\textup{GL}_N(S)$ the group of invertible elements of $\textup{Mat}_N(S)$.
 We denote by $\textup{Mat}_N(S)^{\textup{sym}}$ and $\textup{GL}_N(S)^{\textup{sym}}$ the sets of symmetric matrices in $\textup{Mat}_N(S)$ and $\textup{GL}_N(S)$, respectively. The identity matrix  will  be denoted by $1=1_N=(\delta_{kl})$ where $\delta_{kl}$ is the Kronecker symbol. 

 \

 For every  $X\in \textup{Mat}_N(S)$  we denote by $X_{ij}$ its entries and we write $X=(X_{ij})$. Recall that for  every integer $s\geq 1$ we denote by $X^{(p^s)}$ the matrix $(X_{ij}^{p^s})$ and by $X^t$ the transpose of $X$. Similarly for $v=(v_i)$  a row or column vector with entries $v_i\in S$ we denote by $v^{(p^s)}$ the vector with entries $v_i^{p^s}$.  If $u:S\rightarrow S'$ is a map of sets we write $u(X):=(u(X_{ij}))$ and $u(v)=(u(v_i))$. In what follows we will be led to consider families 
 $(X_1,\ldots,X_n)$ of matrices in $\textup{Mat}_n(S)$; for such a family we denote by $(X_i)_{jk}$ the  $jk$-entries of $X_i$ so we write $X_i=((X_i)_{jk})$. 

 \

 However we will never use the notation $X_{ijk}$ to denote the $jk$-entry of $X_i$; this notation  will be reserved for objects obtained by  ``lowering the indices," in a sense that will be explained  later!)  On the other hand we will also consider matrices  $X^k$ indexed by superscripts $k$ (rather than subscripts);  as a rule these superscripts will never mean ``raising $X$ to the $k$-th power."
 For a matrix $X^k$ we denote by $X^k_{ij}:=(X^k)_{ij}$ its $ij$-entries. For every $n$-tuple of matrices $(X_1,\ldots,X_n)$ we  define the $n$-tuple of matrices $(X^1,\ldots,X^n)$ by the formula 
\begin{equation}\label{treimat}
(X^k)_{ij}=X^k_{ij}:=(X_i)_{jk}=(X_i^t)_{kj};\end{equation}
 and, vice versa, if one is given an $n$-tuple of matrices $(X^1,\ldots,X^n)$ the above formula defines an $n$-tuple of matrices
 $(X_1,\ldots,X_n)$.

\

As in Subsection \ref{generall} let $x=(x_{ij})$ be an $N\times N$  matrix of indeterminates with $n\in \mathbb N$. We consider the group scheme
 $G=\textup{GL}_{N,R_{\pi}}=\textup{Spec }R_{\pi}[x,\det(x)^{-1}]$,
 so for the ring of global functions we have 
$\mathcal O(G)=R_{\pi}[x,\det(x)^{-1}]$ and for the group of points we have $G(R_{\pi})=\textup{GL}_N(R_{\pi})$. 
We consider the ring
$\mathcal A:=\widehat{\mathcal O(G)}=R_{\pi}[x\det(x)^{-1}]^{\widehat{\ }}$
and we consider the prime ideal 
$\mathcal P:=(x-1)$
generated by  the entries $x_{kl}-\delta_{kl}$ of the matrix $x-1$ where  $1:=1_n=(\delta_{kl})$ is the identity matrix. We also consider the 
maximal ideal 
$\mathcal M:=(\pi,\mathcal P)=(\pi,x-1)$
 generated by $\pi$ and $\mathcal P$.
 The residue field $\mathcal A/\mathcal M$ equals $k$. 

 \

\noindent We have a canonical homomorphism 
$$\mathcal A\rightarrow \mathcal A/\mathcal P\simeq R_{\pi},\ \ \ F\mapsto F(1)=F_{|x=1}.$$
Assume we are given an $n$-tuple of $\pi$-Frobenius lifts on $R_{\pi}$ of degree $s$.
Recall that a {\bf $\pi$-connection} of degree $s$ on $G$  is an $n$-tuple 
 $$\Delta^{(s)G}=((\delta^{(s)}_1)^G,\ldots,(\delta^{(s)}_n)^G)$$ 
 of higher $\pi$-derivations of degree $s$ on $\mathcal A$ extending the given higher $\pi$-derivations of degree $s$ on $R_{\pi}$, respectively. To give a 
$\pi$-connection  of degree $s$  on $G$ is the same as to give a structure of  partial $\delta$-ring of degree $s$ 
on 
$\mathcal A$.
We denote by $$\Phi^{(s)G}=((\phi^{(s)}_1)^G,\ldots,(\phi^{(s)}_n)^G)$$ the attached  family of higher 
$\pi$-Frobenius lifts.

\

Recall from Subsection \ref{generall} the Christoffel symbol of the second kind 
 $\Gamma^{(s)}=
  (\Gamma_1^{(s)},\ldots,\Gamma_n^{(s)})$
 of a 
 $\pi$-connection  $((\delta_1^{(s)})^G,\ldots,(\delta_n^{(s)})^G)$  of degree $s$ on $G$. 
In view of  Equation (\ref{treimat})  we write 
\begin{equation}
\label{trematt}
\Gamma_{ij}^{k(s)}:=(\Gamma_i^{(s)})_{jk}=(\Gamma_i^{(s)t})_{kj}\in \mathcal A.\end{equation}
(We will {\it not} use the notation $\Gamma^{(s)}_{ijk}$ to denote the $jk$ entry of $\Gamma^{(s)}_i$; the notation $\Gamma^{(s)}_{ijk}$ will be reserved for the ``Christoffel symbols of the first kind" to be introduced later.)
We have the following formula
\begin{equation}
\label{la11}
((\phi_i^{(s)})^G(x_{kj}))(1)=\delta_{kj}+\pi\Gamma^{(s)k}_{ij}(1),
\end{equation}
where $\delta_{kj}$ is the Kronecker symbol.

\

For an ``intrinsic" description of Christoffel symbols we refer to \cite[Subsect. 4.3]{BM22}.
For a discussion of the analogy with classical differential geometry we refer to  Section \ref{diskus}.
Our Christoffel symbols here are slightly different from the ones in \cite{Bu17} and \cite{Bu19};
the change was necessary in order to accommodate our new, PDE, setting. 

\
 
 Explicitly for every 
 $q^{(s)}\in \textup{GL}_N(R_{\pi})^{\textup{sym}}$ and every  $\pi$-connection $\Delta^{(s)G}$  of degree $s$ as above
we consider the following $N\times N$ matrices 
 with entries in $\mathcal A$:
\begin{equation}
\label{1955}
\begin{array}{rcl}
A_i^{(s)} & := & x^{(p^s)t}\cdot \phi^{(s)}_{i}(q^{(s)})\cdot x^{(p^s)},\\
\ & \ & \ \\
 B^{(s)} & := & (x^tq^{(s)}x)^{(p^s)}.\end{array}\end{equation} 

\

Note that the reduction mod $\mathcal M^2$ of $(\phi_i^{(s)})^G$ is uniquely determined by 
the reduction mod $\mathcal M^2$ of $\Lambda_i^{(s)}$; the latter is completely determined by the reduction mod $\mathcal M$ of $\Gamma_i^{(s)}$.

\

The condition that a $\pi$-connection as above is metric with respect to $q^{(s)}$ (see Subsection \ref{generall}) is  equivalent to the matrix equalities:
\begin{equation}
\label{LALB}
\Lambda_{i}^{(s)t}A_i^{(s)} \Lambda_{i}^{(s)}=B^{(s)}\end{equation}
for $i\in \{1,\ldots,n\}$. 
The condition that a $\pi$-connection as above is symmetric with respect to a torsion symbol
$L^{(s)}=(L^{k(s)})$ (see Subsection \ref{generall}) is equivalent to the equalities:
\begin{equation}
\label{uf22}
(\Lambda_i^{(s)}-1)_{kj}-(\Lambda_j^{(s)}-1)_{ki}=\pi\cdot L_{ij}^{k(s)}(\Lambda^{(s)}).
\end{equation}
Finally we consider  the  ``lowering of the indices" operation by defining
\begin{equation}
\label{loweringg}
\begin{array}{rcl}
L_{ijk}^{(s)} & := & \sum_m L_{ij}^{m(s)}(q_{mk}^{(s)})^{p^s},\\
\ & \ & \ \\
\Gamma_{ijk}^{(s)} & := & \sum_m \Gamma_{ij}^{m(s)} (q_{mk}^{(s)})^{p^s}.\end{array}\end{equation}
Note that, consistent with our previously explained notation,
 $\Gamma_{ijk}^{(s)}$ is {\it not} the $jk$ entry of $\Gamma_i^{(s)}$ (the latter being denoted by $(\Gamma_i^{(s)})_{jk}=\Gamma^{k(s)}_{ij}$, see Equation \ref{trematt}); a similar remark holds 
 for $L_{ijk}^{(s)}$.

 \begin{definition}
The quantities
$\Gamma_{ijk}^{(s)}$ are the {\bf Christoffel symbols of the first kind} (relative to the metric $q^{(s)}$);
$L_{ijk}^{(s)}$ are the {\bf torsion symbols of the first kind} (relative to the metric $q^{(s)}$).
\end{definition}

\begin{remark}\label{nutzii}
 Since the ring endomorphisms $(\phi_i^{(s)})^G$   send any $f\in \mathcal M$ into the ideal $(\pi,f^{p^s})$ hence into the ideal $(\pi,\mathcal M^{p^s})$ we get that for all $\nu\geq 1$ the ring endomorphisms $(\phi_i^{(s)})^G$ send $\mathcal M^{\nu}$
 into the ideal
 $(\pi,\mathcal M^{p^s})^{\nu}$. 
 The latter is contained in the ideal $(\pi^{\nu},\mathcal M^{\nu+1})$ which in its turn is contained in the ideal $\mathcal M^{\nu}$.
 Hence, for every $\nu\geq 1$ the ring homomorphisms 
  $(\phi_i^{(s)})^G$ induce ring endomorphisms
  $(\phi_i^{(s)})^{\mathcal A/\mathcal M^{\nu}}$ of the ring 
  $\mathcal A/\mathcal M^{\nu}$.
     We also have induced $\phi_i^{(s)}$-linear endomorphisms
  $(\phi_i^{(s)})^{\mathcal M^{\nu}/\mathcal M^{\nu+1}}$ of the abelian group
  $\mathcal M^{\nu}/\mathcal M^{\nu+1}$.
\end{remark}
  \

  Note, in particular, that the ring endomorphism $(\phi_i^{(s)})^{\mathcal A/\mathcal M^2}$ 
   is completely determined by its restriction to $R_{\pi}$ and by the classes 
modulo $\mathcal M$
of the Christoffel symbols of the second kind, 
$\Gamma_i^{(s)}$ (i.e., by the images of $\Gamma_i^{(s)}$ in 
$\textup{Mat}_N(k)$) and hence by the classes modulo $\mathcal M$
of the Christoffel symbols of the first kind, $\Gamma_{ijk}^{(s)}$.
Explicitly we have the following formula. 

\begin{lemma}
Modulo $\mathcal M^2$ we have 
\begin{equation}
\label{cucucu} 
(\phi_i^{(s)})^G(x_{kl}-\delta_{kl})\equiv \pi \Gamma^{k(s)}_{il}\ \ \textup{mod}\ \ \mathcal M^2.
\end{equation}
\end{lemma}
{\it Proof}.
To check the formula note first that since $$(x_{kk}-1)^2, \ p(x_{kk}-1)\in \mathcal M^2$$ we have
\begin{equation}
\label{cucsc}
x_{kk}^p= ((x_{kk}-1)+1)^p=\sum_{i=0}^p \left(\begin{array}{c} p\\ i\end{array}\right)(x_{kk}-1)^i\equiv 1\ \ \textup{mod}\ \ \mathcal M^2
\end{equation}
and
\begin{equation}
\label{cucscc}
x_{kl}^p\equiv 0\ \ \textup{mod}\ \ \mathcal M^2\ \ \textup{for}\ \ k\neq l
\end{equation}
hence
\begin{equation}
\label{cucsccc}
x_{kl}^p\equiv \delta_{kl}\ \ \textup{mod}\ \ \mathcal M^2.
\end{equation}
Then using Equations  (\ref{1955}) and (\ref{cucsccc}) we have
$$
\begin{array}{rcl}
(\phi_i^{(s)})^G(x_{kl}) & = & (x^{(p^s)}\Lambda_i^{(s)})_{kl}\\
\ & \ & \ \\
\ & = & (x^{(p^s)}(1+\pi \Gamma_i^{(s)t}))_{kl}\\
\ & \ & \ \\
\ & = & \sum_m x_{km}^{p^s}(\delta_{ml}+\pi\Gamma_{il}^{m(s)})\\
\ & \ & \ \\
\ & \equiv & \delta_{kl}+\pi \Gamma_{il}^{k(s)}  \ \ \ \  \textup{mod}\ \ \mathcal M^2,
\end{array}
$$
which checks formula (\ref{cucucu}). \qed

\

Now, recall  the following result from Part 1; cf. \cite[Thm. 3.30]{BM22} which gives the arithmetic version of the fundamental theorem for existence of the Levi-Civita connection.

\begin{theorem}\label{LCC}
(Arithmetic Levi-Civita Theorem)
Assume $N=n$.
Assume $s\geq 1$ is an integer, $q^{(s)}\in G(R_{\pi})^{\textup{sym}}$ is a metric,   and $L^{(s)}=(L^{1(s)},\ldots,L^{n(s)})$ is a torsion symbol of the second kind. 
 There exists a unique   $\pi$-connection $$\Delta^{(s)\textup{LC}}=((\delta_1^{(s)})^{\textup{LC}}, \ldots, (\delta_n^{(s)})^{\textup{LC}})$$  of degree $s$  on $G$ that is metric with respect to $q^{(s)}$ and symmetric with respect to $L^{(s)}$.
Moreover if $(\Gamma_{ijk}^{\textup{(s)\textup{LC}}})$ and $(L_{ijk}^{(s)})$ are  the corresponding Christoffel symbols and the torsion symbols of the first kind
 then we have the following congruence  modulo $\mathcal M$ in the ring $\mathcal A$
 
$$\Gamma_{ijk}^{(s)\textup{LC}} \equiv -\frac{1}{2}(\delta_{i}^{(s)}q^{(s)}_{jk}+\delta_{j}^{(s)}q^{(s)}_{ik}-\delta_{k}^{(s)}q^{(s)}_{ij})+
\frac{1}{2}
(L^{(s)}_{kij}(1)+L^{(s)}_{ijk}(1)-L^{(s)}_{jki}(1)).$$
\end{theorem}

\medskip

The unique  $\pi$-connection  $\Delta^{(s)\textup{LC}}$  of degree $s$ in Theorem \ref{LCC} is called the {\bf arithmetic  Levi-Civita connection} of degree $s$ attached to $q^{(s)}$ and $L^{(s)}$.

\

For $N$ not necessarily equal to $n$ and  $\phi_i^{(s)}$  a higher $\pi$-Frobenius lift of degree $s$  on $\mathcal A$ extending $\phi_i^{(s)}$ on $R_{\pi}$ we say that $\phi_i^{(s)}$ is  {\bf $\cB_{q^{(s)}}$-symmetric} provided the following  diagram commutes
\begin{equation}
 \label{got-chern}
 \begin{array}{rcl}
\mathcal A& \stackrel{\cB_{q^{(s)}}}{\longrightarrow} & \mathcal A \otimes \mathcal A \\
\cB_{q^{(s)}}  \downarrow &\ &\downarrow \phi_i^{(s)} \otimes \phi_{i,0}^{(s)}\\
 \mathcal A \otimes \mathcal A & \stackrel{\phi_{i,0}^{(s)} \otimes \phi_i^{(s)}}{\longrightarrow} & \mathcal A\end{array}
\end{equation}
where $\cB_{q^{(s)}}$ is defined by $\cB_{q^{(s)}}(x) = x_1^tq^{(s)}x_2$ and $x$, $x_1:=x\otimes 1$, $x_2:=1\otimes x$ are the corresponding $N\times N$ matrices of indeterminates. 
As in \cite{BD} or \cite{Bu17} and with notation as in (\ref{1955})
the commutativity of the diagram (\ref{got-chern}) is equivalent to the equality
\begin{equation}\label{motuu}
A_i^{(s)} \Lambda_i^{(s)} = \Lambda_i^{(s)t}A_i^{(s)}.\end{equation}

\

We also make use of the Chern connection guaranteed by the following result from Part 1; cf. \cite[Thm. 3.38 and Cor. 3.41]{BM22}.

\begin{theorem}\label{ch} For $q^{(s)} \in \textup{GL}_N(R_\pi)^{\textup{sym}}$ and $i \in \{1,\ldots,n\}$ there is a 
unique higher $\pi$-Frobenius lift $\phi_i^{(s)\textup{Ch}}$ of degree $s$ on $\mathcal A$ that is $\cB_{q^{(s)}}$-symmetric  and metric with respect to $q^{(s)}$. Moreover we have the following congruences mod $\mathcal M$ in the ring $\mathcal A$:
$$\Gamma_i^{\textup{(s)\textup{Ch}}} \equiv -\frac{1}{2} \delta_i^{(s)} q^{(s)} \cdot ((q^{(s)})^{(p^s)})^{-1}.$$
\end{theorem}

The $\pi$-connection $\Delta^{(s)\textup{Ch}}=(\delta_1^{(s)\textup{Ch}},\ldots,\delta_n^{(s)\textup{Ch}})$  of degree $s$ 
defined by the higher Frobenius lifts $\phi_i^{(s)\textup{Ch}}$ is called the {\bf arithmetic Chern connection} of degree $s$ attached to $q^{(s)}$.

\

\subsection{Atiyah extensions} \label{nucah} We return now to the setting of Subsection \ref{cannn}.
 
 \begin{notation}
In what follows we fix an  involutive subset $\mathfrak D^{(c)}\subset \mathfrak F^{(c)}$ of cardinality $n$.  We consider the canonical Lie subalgebra (which is also an associative $\mathbb Z$-subalgebra)
$$\mathfrak d\subset \mathfrak f$$ attached
to $\mathfrak D^{(c)}$ and we consider the group scheme $G:=\textup{GL}_{N/R_{\pi}}$ with $N$ not necessarily equal to $n$.
Recall that we denoted  by  $\textup{End}_{\mathbb Z-\textup{mod}}(\mathcal A)$
 the associative $\mathbb Z$-algebra (also viewed as a Lie $\mathbb Z$-algebra)  of all $\mathbb Z$-module endomorphisms of $\mathcal A$ and recall   the ideal $\mathcal M:=(\pi,x-1)$ in $\mathcal A$. For $s\in c\mathbb N$ we denote by 
  $\mathfrak E^{(s)}$ the set of all the higher $\pi$-Frobenius lifts on $\mathcal A$
  of degree $s$ 
  that send $R_{\pi}$ into itself and whose restriction to $R_{\pi}$ is in $\mathfrak D^{(s)}$. We denote by
  $$\mathfrak e^{(s)}\subset \textup{End}_{\mathbb Z-\textup{mod}}(\mathcal A)$$ 
  the $\mathbb Z$-linear span of $\mathfrak E^{(s)}$. By the ``linear independence of characters" $\mathfrak e^{(s)}$ is a free $\mathbb Z$-module with basis $\mathfrak E^{(s)}$.
 The set 
 \begin{equation}
 \label{Esbt}
 \mathfrak E:=\bigcup_{s\in c\mathbb N} \mathfrak E^{(s)}\end{equation}
 has a natural structure of graded monoid of minimal degree $c$.
  We set
  $$\mathfrak e:=\bigoplus_{s\in c\mathbb N}\mathfrak e^{(s)}.$$

  \

  One easily sees that $\mathfrak e$ is a free $\mathbb Z$-submodule with basis $\mathfrak E$
 and also a graded associative
   $\mathbb Z$-subalgebra (hence a graded Lie $\mathbb Z$-sub\-algebra) of $\textup{End}_{\mathbb Z-\textup{mod}}(\mathcal A)$.    So $\mathfrak e$ is isomorphic to the  monoid algebra $\mathbb Z[\mathfrak E]$ attached to the monoid $\mathfrak E$.  We have  a surjective graded associative (and hence Lie) $\mathbb Z$-algebra homomorphism
  $\mathfrak e \rightarrow \mathfrak d$ defined by restriction and we denote by $\mathfrak r$ the kernel of this homomorphism which is therefore  an ideal in the associative $\mathbb Z$-algebra (and hence also a Lie ideal in the Lie $\mathbb Z$-algebra) $\mathfrak e$. So we have 
  an exact sequence  in the category of graded associative (and hence also Lie) $\mathbb Z$-algebras
  \begin{equation}\label{nutziiii}
  0\rightarrow \mathfrak r\rightarrow \mathfrak e\rightarrow \mathfrak d\rightarrow 0
  \end{equation}
  which is split in the category of graded associative (and hence also Lie) $\mathbb Z$-algebras (the splitting in both categories being given for each $s$ by the trivial  $\pi$-connection of degree $s$).  
  \end{notation}
 
\begin{remark}
The extension (\ref{nutziiii}) can be viewed as an arithmetic analogue   of the Atyiah  extension in \cite{A57}. However, the classical Atiyah extension is generally 
 non-split (in the category of vector bundles) and it is this non-splitting  that leads to the non-triviality of Chern classes. 

 \

 In our case the splitting (in the category of graded Lie $\mathbb Z$-algebras) of the extension (\ref{nutziiii}) shows that no analogue of Chern classes  is to be expected to emerge from its consideration. We will introduce below  more refined exact sequences that {\it will} lead to non-trivial  Lie algebra and Hochschild cohomology classes attached to a given ``graded $\pi$-connection" (see Definition \ref{morerefined} and Corollary \ref{cincisute} below). 
By the way the splitting of (\ref{nutziiii}) 
can be viewed as a reflection of the fact that $G$ is  a {\it trivial} torsor 
over itself. A discussion of ``non-trivial torsors" will be undertaken in Subsection \ref{torsors}.
\end{remark}

\begin{definition}
Fix $\mathfrak D^{(c)}$ and $N$ in $\textup{GL}_N$  as above and consider the attached sequence (\ref{nutziiii}).
An {\bf Atiyah extension} is a graded associative $\mathbb Z$-subalgebra $\mathfrak l$  of $\mathfrak e$ such that  $\mathfrak l$  surjects onto $\mathfrak d$. \end{definition}

Hence we have an exact sequence 
in the category of graded associative (hence also Lie) $\mathbb Z$-algebras, which is  also a split exact sequence of  $\mathbb Z$-modules, 
\begin{equation}
\label{nusha1}
0\rightarrow \mathfrak n \rightarrow \mathfrak l \stackrel{r}{\rightarrow} \mathfrak d\rightarrow 0,
\end{equation}
where $\mathfrak n:=\textup{Ker}(r:\mathfrak l\rightarrow \mathfrak d)=\mathfrak r\cap \mathfrak l$.
To (\ref{nusha1})  one can attach the exact sequence in the category of graded associative $\mathbb Z$-algebras (and hence also in the category of graded Lie $\mathbb Z$-algebras)
given by
\begin{equation}\label{nusha222}
0  \rightarrow  \mathfrak n/\mathfrak n^2  \rightarrow  \mathfrak l/\mathfrak n^2  \stackrel{r^{(2)}}{\rightarrow}  \mathfrak d  \rightarrow  0. \end{equation}
This is also a split exact sequence of $\mathbb Z$-modules  because $\mathfrak d$ is a free $\mathbb Z$-module. Note that $\mathfrak n/\mathfrak n^2$ is  an ideal of square zero in the associative $\mathbb Z$-algebra $\mathfrak l/\mathfrak n^2$
hence also an abelian Lie ideal in the Lie $\mathbb Z$-algebra $\mathfrak l/\mathfrak n^2$.
So we can make the following definition.

\begin{definition}\label{morerefined}
 The {\bf Lie characteristic class}
of $\mathfrak l$ is the cohomology class 
$$\kappa(\mathfrak l)\in H^2_{\textup{Lie}}(\mathfrak d,\mathfrak n/\mathfrak n^2)
$$
of the sequence (\ref{nusha222}).
The {\bf Hochschild characteristic class}
of $\mathfrak l$ is the cohomology class 
$$h(\mathfrak l)\in H^2_{\textup{Hoch}}(\mathfrak d,\mathfrak n/\mathfrak n^2)
$$
of the sequence (\ref{nusha222}). Clearly $\kappa(\mathfrak l)$ is the image of $h(\mathfrak l)$
via the natural map between the corresponding cohomology groups.
\end{definition}

 We next consider a $k$-linear analogue of the above construction. Fix $\nu\in \mathbb N$. Let $[\pi^{\nu}]$ be the image of $\pi^{\nu}$ in the $k$-linear space $\mathcal M^{\nu}/\mathcal M^{\nu+1}$ and let $k\cdot [\pi^{\nu}]\subset \mathcal M^{\nu}/\mathcal M^{\nu+1}$ be the one dimensional $k$-linear subspace spanned by $[\pi^{\nu}]$.
 Let $\textup{Fr}:k\rightarrow k$ be the $p$-power Frobenius automorphism.
  For all $s\in c\mathbb N$ we let 
  $$\mathfrak b^{(s)}\subset \textup{End}_{\mathbb Z-\textup{mod}}(\mathcal M^{\nu}/\mathcal M^{\nu+1})$$ be the $k$-linear subspace of all $\textup{Fr}^s$-linear maps $u\in \textup{End}_{\mathbb Z-\textup{mod}}(\mathcal M^{\nu}/\mathcal M^{\nu+1})$ such that
 $$u(\mathcal M^{\nu}/\mathcal M^{\nu+1})\subset k\cdot [\pi^{\nu}]$$
 and we let $\mathfrak a^{(s)}$
 be the $k$-linear subspace  of $\mathfrak b^{(s)}$ consisting of all  $u\in \mathfrak b^{(s)}$ such that
 $$u(k\cdot [\pi^{\nu}])=0.$$
 Finally let
 $$\mathfrak c^{(s)}\subset \textup{End}_{\mathbb Z-\textup{mod}}(k\cdot [\pi^{\nu}])$$
 be the $k$-linear subspace of all $\textup{Fr}^s$-linear maps in $\textup{End}_{\mathbb Z-\textup{mod}}(k\cdot [\pi^{\nu}])$. 

 \

 Define
 $$\mathfrak a:=\bigoplus_{s\in c\mathbb N}\mathfrak a^{(s)},\ \ \mathfrak b:=\bigoplus_{s\in c\mathbb N}\mathfrak b^{(s)},\ \ \mathfrak c:=\bigoplus_{s\in c\mathbb N}\mathfrak c^{(s)}.$$
 Then $\mathfrak b$ and $\mathfrak c$ are graded associative (hence Lie) $\mathbb Z$-subalgebras and also $k$-linear subspaces
 of 
 $$\textup{End}_{\mathbb Z-\textup{mod}}(\mathcal M^{\nu}/\mathcal M^{\nu+1})\ \ \textup{and}\ \ \textup{End}_{\mathbb Z-\textup{mod}}(k\cdot [\pi^{\nu}]),
 $$
 respectively. Also $\mathfrak a$ is a $k$-linear subspace and  a graded ideal of square zero in the associative $\mathbb Z$-algebra $\mathfrak b$ and  hence $\mathfrak a$ is an abelian graded Lie ideal in the Lie $\mathbb Z$-algebra $\mathfrak b$. 
 
 \

 By Remark \ref{nutzii} we have compatible $k$-linear associative (and hence Lie) $\mathbb Z$-algebra restriction maps $\mathfrak l\rightarrow \mathfrak b$ and $\mathfrak d\rightarrow \mathfrak c$, hence an induced map $\mathfrak n\rightarrow \mathfrak a$. Note that the $k$-linear span of the image of the  map $\mathfrak d\rightarrow \mathfrak c$ equals $\mathfrak c$.
 
 \

 For every Atiyah extension $\mathfrak l$ the various objects considered  above fit into  a commutative  diagram  of graded associative (and hence Lie) $\mathbb Z$-algebras with exact rows:
 \begin{equation}
\label{westw}
\begin{array}{ccccccccc}
0 & \rightarrow & \mathfrak n& \rightarrow  &  \mathfrak l & \stackrel{r}{\rightarrow} & \mathfrak d & \rightarrow & 0 \\
\ & \ & \downarrow& \ & \downarrow  & \ & ||  & \ & \ \\
0 & \rightarrow & \mathfrak n/\mathfrak n^2 & \rightarrow & \mathfrak l/\mathfrak n^2& \stackrel{r^{(2)}}{\rightarrow} & \mathfrak d & \rightarrow & 0\\
\ & \ & \downarrow& \ & \downarrow  & \ & \downarrow  & \ & \ \\
0 & \rightarrow & \mathfrak a & \rightarrow & \mathfrak b & \rightarrow & \mathfrak 
c & \rightarrow & 0
\end{array} \end{equation}

        The bottom row in (\ref{westw})  is, of course, independent of $\mathfrak l$ 
and its structure can be explicitly described as follows.
 Consider
a $k$-basis 
\begin{equation}
\label{kbasis}
\beta_0,\ldots,\beta_M\end{equation}
 of $\mathcal M^{\nu}/\mathcal M^{\nu+1}$
where $$\beta_0=[\pi^{\nu}]$$ 
and 
the identification
$$\textup{Mat}_{(M+1)\times 1} \simeq \mathcal M^{\nu}/\mathcal M^{\nu+1},\ \ (c_0,c_1,\ldots,c_M)^t\mapsto \sum_{i=0}^M c_i\beta_i,\ \ c_i\in k.$$
Then $\mathfrak b^{(s)}$ identifies with 
the $k$-linear space  of matrices
\begin{equation}
\label{oprottu}
\left(\begin{array}{cc}\lambda_0\cdot \textup{Fr}^s & \lambda\cdot \textup{Fr}^s \\ 0 & 0\end{array}\right),\ \ \lambda_0\in k,\ \ \lambda \in\textup{Mat}_{1\times M}(k)\end{equation}
acting via left multiplication on column vectors $(c_0,c_1,\ldots,c_M)^t$. 
Here the left bottom $0$ is the $M\times 1$ matrix with zero entries and the right bottom $0$ is the $M\times M$ matrix with zero entries; we will adopt this convention from now on.

\

\noindent Furthermore under the identification
$$k\simeq  k\cdot [\pi^{\nu}]=\pi^{\nu}R_{\pi}/\pi^{\nu+1} R_{\pi},\ \ \lambda\mapsto \lambda \pi^{\nu}\ \ \textup{mod}\ \ \pi^{\nu+1}$$
the space $\mathfrak c^{(s)}$ identifies with the $k$-linear space of maps $k\rightarrow k$ of the form
$\lambda_0\cdot \textup{Fr}^s$ with $\lambda\in k$. If, under the above identification, $u\in \mathfrak c^{(s)}$  corresponds to  $\lambda_0\cdot \textup{Fr}^s$ then we have
$$\lambda_0=\frac{u(\pi^{\nu})}{\pi^{\nu}}\ \ \textup{mod}\ \ \pi.$$
In particular, if $u$ is the image of some $\phi_i^{(s)}\in \mathfrak d$ then
$$\lambda_0=\left(\frac{\phi_i^{(s)}(\pi)}{\pi}\right)^{\nu}
\equiv (\delta_i^{(s)}\pi)^{\nu}
\ \ \textup{mod}\ \ \pi.$$

\

\

Under the identifications above 
the map $\mathfrak b\rightarrow \mathfrak c$ corresponds to the map  sending a matrix (\ref{oprottu})  to $\lambda_0\cdot \textup{Fr}^s$; moreover
$\mathfrak a^{(s)}$ identifies  with the $k$-linear space of matrices
\begin{equation}
\label{rr0}
\left(\begin{array}{cc}0 & 
\lambda\cdot \textup{Fr}^s
\\ 0 & 0\end{array}\right)\end{equation}
where $\lambda\in \textup{Mat}_{1\times M}(k)$.
The bottom row of (\ref{westw}) is  split in the category of graded associative (and also Lie) $\mathbb Z$-algebras: a splitting  $\mathfrak c\rightarrow \mathfrak b$ is given by
 $$\lambda_0\cdot \textup{Fr}^s\mapsto \left(\begin{array}{cc} \lambda_0\cdot \textup{Fr}^s & 0 \\ 0 & 0 \end{array}\right).$$

 \

 The natural structure of Lie
$\mathfrak c$-$\mathfrak c$-bimodule on 
 $\mathfrak a$ defined by the bottom row of (\ref{westw}), viewed as an associative $\mathbb Z$-algebra sequence, is given for $X:=\lambda_0\cdot \textup{Fr}^r\in \mathfrak c$ and for a matrix as in (\ref{rr0})
via the rules 
\begin{equation}
\label{rulibu}
\begin{array}{rcl}
X \cdot \left(\begin{array}{cc}0 & \lambda\cdot \textup{Fr}^s\\ 0 & 0\end{array}\right) & := & 
\left(\begin{array}{cc}\lambda_0\cdot \textup{Fr}^r & 0\\ 0 & 0\end{array}\right) \left(\begin{array}{cc}0 & \lambda\cdot \textup{Fr}^s\\ 0 & 0\end{array}\right)\\
\ & \ & \ \\
\ & = & \left(\begin{array}{cc}0 & \lambda_0\cdot \lambda^{p^r}\cdot \textup{Fr}^{s+r}\\ 0 & 0\end{array}\right),\end{array}\end{equation}
\begin{equation}
\left(\begin{array}{cc}0 & \lambda\cdot \textup{Fr}^s\\ 0 & 0\end{array}\right)\cdot  X=\left(\begin{array}{cc}0 & \lambda\cdot \textup{Fr}^s\\ 0 & 0\end{array}\right)\left(\begin{array}{cc}\lambda_0\cdot \textup{Fr}^r & 0\\ 0 & 0\end{array}\right)
=0.
\end{equation}

\

The natural   Lie
$\mathfrak c$-module structure on 
 $\mathfrak a$ defined by the bottom row of (\ref{westw}), viewed as a Lie algebra sequence, is given for $X:=\lambda_0\cdot \textup{Fr}^r\in \mathfrak c$ and for a matrix as in (\ref{rr0})
via the rule 
\begin{equation}
\label{ruli}
\begin{array}{rcl}
X \cdot \left(\begin{array}{cc}0 & \lambda\cdot \textup{Fr}^s\\ 0 & 0\end{array}\right) & := & \left[
\left(\begin{array}{cc}\lambda_0\cdot \textup{Fr}^r & 0\\ 0 & 0\end{array}\right),\left(\begin{array}{cc}0 & \lambda\cdot \textup{Fr}^s\\ 0 & 0\end{array}\right)\right]\\
\ & \ & \ \\
\ & = & \left(\begin{array}{cc}0 & \lambda_0\cdot \lambda^{p^r}\cdot \textup{Fr}^{s+r}\\ 0 & 0\end{array}\right).\end{array}\end{equation}
So the above Lie
$\mathfrak c$-module structure on 
 $\mathfrak a$  is the same as 
 the  Lie
$\mathfrak c$-module structure induced by
 the  $\mathfrak c$-module structure on $\mathfrak a$ coming from the above $\mathfrak c$-$\mathfrak c$-bimodule structure  on $\mathfrak a$. 
 On the other hand the $\mathfrak c^{\textup{op}}$-module structure on $\mathfrak a$ induced by the $\mathfrak c$-$\mathfrak c$-bimodule structure on $\mathfrak a$ is trivial.
 Finally 
 \begin{equation}
 \label{erer}
 \mathfrak a\simeq \mathfrak c^M\ \ \ \textup{as}\ \ \ \mathfrak c-\textup{modules}.\end{equation}

 \

 Whenever we talk about the $\mathfrak c$-module structure,  $\mathfrak c^{\textup{op}}$-module structure, $\mathfrak c$-$\mathfrak c$-bimodule structure, and Lie $\mathfrak c$-module structure on $\mathfrak a$ we shall be referring to the structures just described.
 They induce a $\mathfrak d$-module structure,  $\mathfrak d^{\textup{op}}$-module structure, $\mathfrak d$-$\mathfrak d$-bimodule structure, and Lie $\mathfrak d$-module structure on $\mathfrak a$ via the map $\mathfrak d\rightarrow \mathfrak c$. With this structures the left vertical maps
 in the diagram (\ref{westw}) are $\mathfrak d$-module maps.

\begin{remark}\label{pasulpaull}
There is a direct link between the algebra $\mathfrak c$ considered above and skew polynomial algebras which we record here. Indeed
for our fixed $c\in \mathbb N$ consider    the skew monoid algebra 
$k[\mathbb N\cup\{0\}]$
 attached to the additive monoid $\mathbb N\cup\{0\}$ acting on $k$ via 
the rule $m\cdot \lambda=\lambda^{p^{cm}}$, $\lambda\in k$, $m\in \mathbb N\cup\{0\}$. Then  $k[\mathbb N\cup\{0\}]$ identifies with the associative $\mathbb Z$-algebra of {\bf skew polynomials}
 $k[T,\textup{Fr}^c]$ whose elements are formal sums $\sum_{i\geq 0}\lambda_i T^i$ with $\lambda_i\in k$ 
 and multiplication  defined by $T\cdot \lambda=\lambda^{p^c}T$. 

 \

 Denote 
by $\mathfrak m$ the  ideal of $k[T,\textup{Fr}^c]$ generated by $T$. We view $\mathfrak m$ 
 as a graded associative $\mathbb Z$-algebra with  component of degree $s\in \mathbb N$ equal to $k\cdot T^{s/c}$ if $c$ divides $s$ and $0$ otherwise.
Then we have a graded associative $\mathbb Z$-algebra $k$-linear isomorphism 
 $$\mathfrak c\simeq \mathfrak m,\ \ \lambda\cdot \textup{Fr}^s\mapsto \lambda 
 \cdot T^{s/c},\ \ \text{for}\ \ \lambda\in k,\ s\in c\mathbb N.$$ \end{remark}
   
   \

   \subsection{Definition of curvature and characteristic classes}\label{gencur} 
In what follows we fix  an  involutive subset $\mathfrak D^{(c)}\subset \mathfrak F^{(c)}$ of cardinality $n$  and labeling symbols $\omega^{(s)}$ for $s\in c\mathbb N$  as in Subsection \ref{cannn}.  We consider the canonical Lie subalgebra 
$\mathfrak d\subset \mathfrak f$ attached
to $\mathfrak D^{(c)}$. For all $s\in c\mathbb N$
we have at our disposal 
the set
$$\mathfrak D^{(s)}=\{\phi^{(s)}_1,\ldots,\phi^{(s)}_n\}$$
of higher $\pi$-Frobenius lifts  of degree $s$ on $R_{\pi}$ and  the corresponding families
$$\Phi:=(\Phi^{(s)})_{s\in c\mathbb N}=(\phi^{(s)}_1,\ldots,\phi^{(s)}_n)_{s\in c\mathbb N},$$
$$\Delta:=(\Delta^{(s)})_{s\in c\mathbb N}=(\delta^{(s)}_1,\ldots,\delta^{(s)}_n)_{s\in c\mathbb N}$$
of higher $\pi$-Frobenius lifts and higher $\pi$-derivations on $R_{\pi}$, respectively.
So $\Delta$ gives a structure of partial $\delta$-ring of degree $s$ on $R_{\pi}$ for every $s\in c\mathbb N$.

\

 We have the
family of structure constants $(\ell_{ij}^{k(s,r)})$ of the Lie algebra 
$\mathfrak d$
and for each $s\in c\mathbb N$ we have the  
associative symbol $\alpha^{(s)}=(\alpha^{1(s)},\ldots,\alpha^{n(s)})$ and the 
Lie symbol $\ell^{(s)}:=(\ell^{1(s)},\ldots,\ell^{n(s)})$; recall that $\ell^{k(s)}_{ij}=\alpha^{k(s)}_{ji}-\alpha^{k(s)}_{ij}$. We freely use here and below the notation in Subsection \ref{nucah}.    In addition  we assume that $\nu$ in that subsection satisfies $\nu=1$.
In this case we have $M=n^2$ and the elements $\beta_1,\ldots,\beta_M$ in (\ref{kbasis}) will be taken to be the classes of the polynomials $x_{kl}-\delta_{kl}$ in the following order:
$$x_{11}-1,x_{12},\ldots, x_{n,n-1},x_{nn}-1.$$

\begin{definition}
Assume that for a fixed $c\in \mathbb N$ we are given a collection 
$$\Delta^{G}=(\Delta^{(s)G})_{s\in c\mathbb N}=((\delta_1^{(s)})^G,\ldots,(\delta_n^{(s)})^G)_{s\in c\mathbb N},$$ 
 of  $\pi$-connections of degree $s$ 
on $G$  with attached higher $\pi$-Frobenius lifts 
$$\Phi^{G}=(\Phi^{(s)G})_{s\in c\mathbb N}=((\phi_1^{(s)})^G,\ldots,(\phi_n^{(s)})^G)_{s\in c\mathbb N};$$
 we refer to such a collection as a {\bf graded $\pi$-connection} on $G$. 
 A  graded $\pi$-connection as above induces a unique $\mathbb Z$-linear map
\begin{equation}
\nabla:\mathfrak d\rightarrow \mathfrak e,\ \ \ X\mapsto \nabla_X\end{equation}
such that
$$\nabla_{\phi_i^{(s)}}:=(\phi_i^{(s)})^G.$$

\

To give a  graded $\pi$-connection as above is equivalent to giving $\nabla$ so  we also refer to $\nabla$ a {\bf  graded $\pi$-connection}. 
If $\mathfrak l$ is an Atiyah extension we say that $\nabla$ is {\bf $\mathfrak l$-valued} if
$\nabla(\mathfrak d)\subset \mathfrak l$. In this case
 $\nabla$ is a section of $r:\mathfrak l\rightarrow \mathfrak d$ in the category of graded $\mathbb Z$-modules. 

 \

 We may consider  the $\mathbb Z$-bilinear maps
   \begin{equation}
  \label{nua}\begin{array}{ll}
   \mathcal S:\mathfrak d\times \mathfrak d\rightarrow \mathfrak n, &  \mathcal S(X,Y):=
   \nabla_X\nabla_Y-\nabla_{XY}\\
   \ & \ \\
  \mathcal R:\mathfrak d\times \mathfrak d\rightarrow \mathfrak n, &  \mathcal R(X,Y):=\mathcal S(X,Y)-\mathcal S(Y,X)=[\nabla_X,\nabla_Y]-\nabla_{[X,Y]}\end{array}\end{equation}
  for $X,Y\in \mathfrak d$. The right hand sides of Equation  (\ref{nua}) are computed  in $\mathfrak l$ but in fact belong to $\mathfrak n$. 
  We refer to the map $\mathcal R$ above as the {\bf curvature} of the $\mathfrak l$-valued graded $\pi$-connection $\nabla$.\end{definition}

\begin{remark}
  In view of our definition of the bracket on $\mathfrak d$, 
  if we specialize $X=\phi^{(s)}_i,Y=\phi^{(r)}_j\in \mathfrak D$ 
  the following holds:
  $$\begin{array}{rcl}
  \mathcal R(\phi_i^{(s)},\phi_j^{(r)}) &= &
  (\phi_i^{(s)})^G(\phi_j^{(r)})^G -(\phi_j^{(r)})^G(\phi_i^{(s)})^G -
\sum_k \ell_{ij}^{k(s,r)}(\phi_k^{(s+r)})^G\\
\ & \ & \ \\
\ & = & (\phi_i^{(s)})^G(\phi_j^{(r)})^G -(\phi_j^{(r)})^G(\phi_i^{(s)})^G -
(\phi_{(i\star j)_{s,r}}^{(s+r)})^G+(\phi_{(j\star i)_{r,s}}^{(s+r)})^G.\end{array}$$
In particular  $\mathcal R(\phi_i^{(s)},\phi_j^{(r)})(a)\in \pi \mathcal A$ for all $a\in \mathcal A$.
\end{remark}

\

\begin{definition}\label{arct}
We set 
$$R_{ijl}^{k(s,r)}:=\frac{1}{\pi}\mathcal R(\phi_i^{(s)},\phi_j^{(r)})(x_{kl})\in \mathcal A.$$
Given a family of metrics $(q^{(s)})_{s\in c\mathbb N}$
the {\bf arithmetic Riemann curvature tensor} attached to $\nabla$ is the collection of elements
$$R_{ijkl}^{(s)}:=\sum_m R_{ijk}^{m(s,s)} (q_{ml}^{(s)})^{p^{2s}}\in \mathcal A,$$
with 
$i,j,k,l\in \{1,\ldots,n\},\ s,r\in c \mathbb N$. We write
$$R_{ijl}^k:=R_{ijl}^{k(c,c)},\ \ \ R_{ijkl}:=R_{ijkl}^{(c)}.$$
\end{definition}

\

As in \cite{Bu17, Bu19}  multiplication by the matrix $((q_{ij}^{(s)})^{p^{2s})}$ is an analogue of ``lowering the indices by a metric" in classical differential geometry. We introduce more notation as follows.

\begin{definition}
We consider the 
 compositions 
  \begin{equation}
 \label{brechtta00}
 \nabla^{(2)}:\mathfrak d\stackrel{\nabla}{\rightarrow}  \mathfrak l\rightarrow \mathfrak l/\mathfrak n^2,\ \ \ X\mapsto \nabla^{(2)}_X,\end{equation}
 \begin{equation}
 \label{brechtta1}
 \overline{\nabla}:\mathfrak d\stackrel{\nabla}{\rightarrow}  \mathfrak l\rightarrow \mathfrak b,\ \ \ X\mapsto \overline{\nabla}_X,\end{equation}
\begin{equation}
 \mathcal S^{(2)}:
 \mathfrak d\times \mathfrak d\stackrel{\mathcal S}{\rightarrow} \mathfrak n\rightarrow  \mathfrak n/\mathfrak n^2. \end{equation}
\begin{equation}
\label{babanicaadu}
 \overline{\mathcal S}:
 \mathfrak d\times \mathfrak d\stackrel{\mathcal S}{\rightarrow} \mathfrak n\rightarrow   \mathfrak a. \end{equation}
\begin{equation}
 \mathcal R^{(2)}:
 \mathfrak d\times \mathfrak d\stackrel{\mathcal R}{\rightarrow} \mathfrak n\rightarrow  \mathfrak n/\mathfrak n^2. \end{equation}
\begin{equation}
\label{babanicaa}
 \overline{\mathcal R}:
 \mathfrak d\times \mathfrak d\stackrel{\mathcal R}{\rightarrow} \mathfrak n\rightarrow   \mathfrak a. \end{equation}
The following equalities hold for $X,Y\in \mathfrak d$:
\begin{equation}
 \mathcal R^{(2)}(X,Y)=[\nabla^{(2)}_X,\nabla^{(2)}_Y]-\nabla^{(2)}_{[X,Y]},
 \end{equation}
 \begin{equation}
 \label{ppuu}
 \overline{\mathcal R}(X,Y)=[\overline{\nabla}_X,\overline{\nabla}_Y]-\overline{\nabla}_{[X,Y]}.
 \end{equation}
  \end{definition}
  
  \medskip
  
  Note that $\mathcal R^{(2)}$ is a Lie $2$-cocycle of $\mathfrak d$ with values in the Lie $\mathfrak d$-module $\mathfrak n/\mathfrak n^2$ and it defines 
 the class $\kappa(\mathfrak l)$. 
 Similarly $\mathcal S^{(2)}$ is a Hochschild $2$-cocycle of $\mathfrak d$ with values in the $\mathfrak d$-$\mathfrak d$-bimodule $\mathfrak n/\mathfrak n^2$ and it defines 
 the class $h(\mathfrak l)$.
 On the other hand $\overline{\mathcal R}$ is a Lie $2$-cocycle 
 of $\mathfrak d$ with values in the Lie $\mathfrak d$-module $\mathfrak a$ (where $\mathfrak a$ is viewed as a Lie $\mathfrak d$-module via its structure of Lie $\mathfrak c$-module). 
 Also $\overline{\mathcal S}$ is a Hochschild $2$-cocycle 
 of $\mathfrak d$ with values in the $\mathfrak d$-$\mathfrak d$-bimodule $\mathfrak a$ (where $\mathfrak a$ is viewed as a $\mathfrak d$-$\mathfrak d$-bimodule via its structure of $\mathfrak c$-$\mathfrak c$-bimodule).

\medskip

For a given  graded $\pi$-connection $\nabla$ there is a smallest Atiyah extension $\mathfrak l$ such that $\nabla$ is $\mathfrak l$-valued; cf. the Definition below.

\begin{definition}\label{liehulll}
Assume $\nabla:\mathfrak d\rightarrow \mathfrak e$ is a  graded $\pi$-connection.
For $s\in c\mathbb N$ let $\mathfrak L_{\nabla}^{(s)}\subset \mathfrak E^{(s)}$ be the set of all 
 products
\begin{equation}
\label{furt}
\nabla_{\phi_{i_1}^{(s_1)}}\nabla_{\phi_{i_2}^{(s_2)}}\cdots
\nabla_{\phi_{i_m}^{(s_m)}}
\end{equation}
where $m\in \mathbb N$, $i_1,\ldots,i_m\in \{1,\ldots,n\}$, 
  $s_1,\ldots,s_m\in c\mathbb N$, $s_1+\cdots+s_m=s$. (These products are not necessarily distinct for distinct tuples $(i_1,\ldots,i_m,s_1,\ldots,s_m)$.) Consider the submonoid 
  $$\mathfrak L_{\nabla}:=\bigcup_{s\in c\mathbb N} \mathfrak L_{\nabla}^{(s)}$$
  of $\mathfrak E$.
The {\bf Atiyah extension attached to  $\nabla$} is the $\mathbb Z$-submodule $\mathfrak l_{\nabla}$ of $\mathfrak e$ generated by $\mathfrak L_{\nabla}$. \end{definition}
  
\begin{remark}  Note that  $\mathfrak l_{\nabla}$ has a natural structure of graded 
associative $\mathbb Z$-algebra
$$\mathfrak l_{\nabla}=\bigoplus_{s\in c\mathbb N} \mathfrak l_{\nabla}^{(s)}$$ with component 
$\mathfrak l_{\nabla}^{(s)}$ of degree $s$ spanned by $\mathfrak L_{\nabla}^{(s)}$.
This shows that $\mathfrak l _{\nabla}$
is indeed an Atiyah extension.
 By the ``linear independence of characters" each $\mathfrak l_{\nabla}^{(s)}$ is a free $\mathbb Z$-module with basis $\mathfrak L_{\nabla}^{(s)}$. In particular  $\mathfrak l^{(c)}$ is a free $\mathbb Z$-module
with basis $$\nabla_{\phi_1^{(c)}},\ldots,\nabla_{\phi_n^{(c)}}$$ 
and $\mathfrak l_{\nabla}$ is isomorphic to the monoid algebra
$\mathbb Z[\mathfrak L_{\nabla}]$ attached to  $\mathfrak L_{\nabla}$. The graded component $\mathfrak n_{\nabla}^{(s)}$ of the ideal
$$\mathfrak n_{\nabla}:=\textup{Ker}(\mathfrak l_{\nabla}\rightarrow \mathfrak d)$$
coincides with the $\mathbb Z$-submodule of $\mathfrak e$ generated by those differences
$$\nabla_{\phi_{i_1}^{(s_1)}}\nabla_{\phi_{i_2}^{(s_2)}}\cdots
\nabla_{\phi_{i_m}^{(s_m)}}-\nabla_{\phi_{j_1}^{(r_1)}}\nabla_{\phi_{j_2}^{(r_2)}}\cdots
\nabla_{\phi_{j_n}^{(r_n)}}
$$
with the properties that 
$$s_1+...+s_m=r_1+...+r_n=s\ \ \textup{and}\ \ \phi_{i_1}^{(s_1)}\cdots \phi_{i_m}^{(s_m)}=\phi_{j_1}^{(r_1)}\cdots \phi_{j_n}^{(r_n)}.$$
Note that $\mathfrak n^{(c)}=0$.
For $e\in \mathbb N$ set
$$\mathfrak a^{(\geq ec)}:=\bigoplus_{m\geq e}\mathfrak a^{(mc)}.$$
Then the images of $\overline{\mathcal S}$ and $\overline{\mathcal R}$ are  contained in $\mathfrak a^{(\geq 2c)}$. Hence we can view  
$\overline{\mathcal R}$ as a Lie $2$-cocycle of $\mathfrak d$ with values in the Lie $\mathfrak d$-module $\mathfrak a^{(\geq 2c)}$. Similarly we can view $\overline{\mathcal S}$ as a Hochschild $2$-cocycle of $\mathfrak d$ with values in the $\mathfrak d$-$\mathfrak d$-bimodule $\mathfrak a^{(\geq 2c)}$.
\end{remark}

\begin{definition}\label{cae}
The {\bf conormal $\mathfrak d$-module} of $\nabla$ is the $\mathfrak d$-module
$$\mathfrak n_{\nabla}/\mathfrak n^2_{\nabla}.$$
The {\bf reduced 
conormal $\mathfrak d$-module} of $\nabla$ is the  $\mathfrak d$-module
$$\overline{\mathfrak n}_{\nabla}:=\textup{Im}(\mathfrak n_{\nabla}\rightarrow \mathfrak a).$$
The {\bf Lie characteristic class} of $\nabla$ is the Lie characteristic class of $\mathfrak l_{\nabla}$ and is denoted by
$\kappa_{\nabla}$; so 
$$\kappa_{\nabla}:=\kappa(\mathfrak l_{\nabla})\in H^2_{\textup{Lie}}(\mathfrak d,\mathfrak n_{\nabla}/\mathfrak n^2_{\nabla}).$$ 
The {\bf Hochschild characteristic class} of $\nabla$ is the Hochschild characteristic class of $\mathfrak l_{\nabla}$ and is denoted by
$h_{\nabla}$; so 
$$h_{\nabla}:=h(\mathfrak l_{\nabla})\in H^2_{\textup{Hoch}}(\mathfrak d,\mathfrak n_{\nabla}/\mathfrak n^2_{\nabla}).$$
The {\bf reduced Lie characteristic class} of $\nabla$ is the image 
$$\overline{\kappa}_{\nabla}\in H^2_{\textup{Lie}}(\mathfrak d,\mathfrak a^{(\geq 2c)})$$
of the  Lie $2$-cocycle $\overline{\mathcal R}.$

\noindent The {\bf reduced Hochschild characteristic class} of $\nabla$ is the image 
$$\overline{h}_{\nabla}\in H^2_{\textup{Hoch}}(\mathfrak d,\mathfrak a^{(\geq 2c)})$$
of the Hochschild $2$-cocycle $\overline{\mathcal S}$
 \end{definition}

  \begin{remark}\label{maisuss}\ Consider the following observations.

  \

1) The classes $\kappa_{\nabla}$ and $h_{\nabla}$ only depend on some 
   ``purely discrete/combinatorial" data (namely on monoid homomorphisms) so they can be viewed as analogous to  ``topological" invariants of $\nabla$. 

   \

   Moreover we claim that these classes are induced by the  purely group theoretic cohomology class $h_{\mathfrak D}$ that only depends on the Weil monoid $\mathfrak D$ (and hence is independent of $\nabla$); the emerging picture is somewhat analogous to the Chern-Weil theory. Indeed in the context of a graded $\pi$-connection $\nabla$ we have a surjective monoid homomorphism $\mathbb M_{n,c}^+\rightarrow \mathfrak D$ which lifts canonically to a surjective graded monoid homomorphism $\mathbb M_{n,c}^+\rightarrow \mathfrak L_{\nabla}$ that sends every $i\in \{1,\ldots,n\}$ into $\nabla_{\phi_i^{(c)}}$.
   We get a surjective graded associative $\mathbb Z$-algebra homomorphism $\mathfrak l_{n,c}\rightarrow \mathfrak l_{\nabla}$.
   We obtain a commutative diagram with exact rows in the category of graded associative $\mathbb Z$-algebras 
       \begin{equation}\label{focstar}
   \begin{array}{ccccccccc}
0 & \rightarrow & \mathfrak n_{\mathfrak D}/\mathfrak n_{\mathfrak D}^2 & \rightarrow & \mathfrak l_{n,c}/\mathfrak n_{\mathfrak D}^2& \stackrel{r_{\mathfrak D}^{(2)}}{\rightarrow} & \mathfrak d & \rightarrow & 0\\
\ & \ & \downarrow& \ & \downarrow  & \ & \downarrow  & \ & \ \\
0 & \rightarrow & \mathfrak n_{\nabla}/\mathfrak n_{\nabla}^2 & \rightarrow & \mathfrak l_{\nabla}/\mathfrak n_{\nabla}^2& \stackrel{r^{(2)}}{\rightarrow} & \mathfrak d & \rightarrow & 0\end{array}
    \end{equation} 
   The left vertical arrow is a $\mathfrak d$-module homomorphism, the middle and left vertical arrows are surjective,  and the classes $\kappa_{\nabla}$ and $h_{\nabla}$ are the images 
   of the classes $\kappa_{\mathfrak D}$ and $h_{\mathfrak D}$, respectively,
   via the group homomorphisms
   $$H^2_{\textup{Lie}}(\mathfrak d,\mathfrak n_{\mathfrak D}/\mathfrak n_{\mathfrak D}^2)\rightarrow H^2_{\textup{Lie}}(\mathfrak d,\mathfrak n_{\nabla}/\mathfrak n_{\nabla}^2),\ \ \ H^2_{\textup{Hoch}}(\mathfrak d,\mathfrak n_{\mathfrak D}/\mathfrak n_{\mathfrak D}^2)\rightarrow H^2_{\textup{Hoch}}(\mathfrak d,\mathfrak n_{\nabla}/\mathfrak n_{\nabla}^2).$$
   
   \
    
  2) We have canonical isomorphisms of Lie $\mathfrak c$-modules
  and  of $\mathfrak c$-$\mathfrak c$-bimodules
  $$\mathfrak a^{(\geq 2c)}\simeq \mathfrak a\simeq \mathfrak c^{n^2}$$
  where the first isomorphism is given by ``right division by $\textup{Fr}^c$"
  and the second isomorphism was noted previously, see Equation \ref{erer}.
  So we can canonically identify $\overline{\kappa}_{\nabla}$ and $\overline{h}_{\nabla}$ with  elements (which could still denoted by) 
  \begin{equation}
  \overline{\kappa}_{\nabla}\in H^2_{\textup{Lie}}(\mathfrak d,\mathfrak c)^{n^2},\ \ \ 
  \overline{h}_{\nabla}\in H^2_{\textup{Hoch}}(\mathfrak d,\mathfrak c)^{n^2},\end{equation}
  respectively.
  The groups $H^2_{\textup{Lie}}(\mathfrak d,\mathfrak c)$ and 
  $H^2_{\textup{Hoch}}(\mathfrak d,\mathfrak c)$
  only depend on the isomorphism class of the graded monoid $\mathfrak D$ and on the field $k$ (in particular they do not depend on $\nabla$). It would be interesting to find an ``explicit" description of these groups. 
    
    \

  3) We have the implications
  $$h_{\nabla}=0\ \ \Rightarrow \ \ \kappa_{\nabla}=0\ \ \Rightarrow\ \ \overline{\kappa}_{\nabla}=0,$$
   $$h_{\nabla}=0\ \ \Rightarrow \ \ \overline{h}_{\nabla}=0\ \ \Rightarrow\ \ \overline{\kappa}_{\nabla}=0,$$
   $$\mathfrak n/\mathfrak n^2\ \ \textup{is a trivial $\mathfrak d$-module}\ \ \Rightarrow\ \ 
   \overline{\mathfrak n}\ \ \textup{is a trivial $\mathfrak d$-module}.$$

\

 4) We will later see that,  in remarkable cases, such as that of the Levi-Civita connection,   the  reduced conormal module $\overline{\mathfrak n}$ and the class $\overline{\kappa}_{\nabla}$  are ``generally" non-trivial (and hence $\mathfrak n_{\nabla}/\mathfrak n_{\nabla}^2$ and  the  classes $\kappa_{\nabla},h_{\nabla},\overline{h}_{\nabla}$ will ``generally" be non-trivial); cf. Remark \ref{sidecenu} and Corollary \ref{miarevenit} below.  One step in checking this is the following proposition. \end{remark}

\begin{proposition} \label{cincisute}
Assume  there exist $X,Y\in \mathfrak D^{(c)}$ such that $[X,Y]=0$ and $\overline{\mathcal R}(X,Y)\neq 0$. 
Then the reduced Lie characteristic class  $\overline{\kappa}_{\nabla}$ is non-trivial. Hence (by Remark \ref{maisuss}, assertion 3) the classes
$\kappa_{\nabla},h_{\nabla},\overline{h}_{\nabla}$ are also non-trivial.
\end{proposition}

{\it Proof}. 
Assume  $\overline{\kappa}_{\nabla}=0$. Then there exists a $\mathbb Z$-module map
$g:\mathfrak d\rightarrow \mathfrak a^{(\geq 2c)}$ such that
$$\overline{\mathcal R}(X,Y)=X g(Y)-Yg(X)-g([X,Y])=X g(Y)-Yg(X)$$
because $[X,Y]=0$. Since $X,Y\in \mathfrak d^{(c)}$ and $g(X),g(Y)\in \mathfrak a^{(\geq 2c)}$ it follows that $X g(Y)-Yg(X)\in \mathfrak a^{(\geq 3c)}$. 
On the other hand $\overline{\mathcal R}(X,Y)\in \mathfrak a^{(2c)}$.
This forces $\overline{\mathcal R}(X,Y)=0$, a contradiction.
\qed

\begin{remark}
In contrast with Corollary \ref{cincisute}
the map  $\overline{\mathcal R}:\mathfrak d \times \mathfrak d \rightarrow \mathfrak a$
is easily seen to be a Lie $2$-coboundary for the $\mathfrak d$-module $\mathfrak a$ so its class in the group $H^2_{\textup{Lie}}(\mathfrak d,\mathfrak a)$  vanishes. This follows from the fact that the bottom exact sequence in Equation (\ref{westw}) is split in the category of  Lie $\mathbb Z$-algebras. \end{remark}

\subsection{Formulae for connection mod $\pi$}
\label{wq1}
Let $\Gamma^{k(s)}_{ij}\in \mathcal A$ be the components of the Christoffel
symbols of the second kind attached to a 
 graded $\pi$-connection  and let $\overline{\Gamma}_{ij}^{k(s)}\in k=\mathcal A/\mathcal M$
be the images of $\Gamma^{k(s)}_{ij}$, where we recall that $\mathcal M:=(\pi,x-1)$. 

\

The map (\ref{brechtta1})
is entirely determined by the elements  $\overline{\Gamma}_{ij}^{k(s)}\in k$  and hence (given  a family of metrics $(q^{(s)})_{s\in c\mathbb N}$)  also by the images
$\overline{\Gamma}_{ijk}^{(s)}\in k$ of the components $\Gamma_{ijk}^{(s)}\in \mathcal A$  of the Christoffel
symbol of the first kind. Explicitly, if $\textup{Fr}:k\rightarrow k$ is the $p$-power Frobenius on $k$
and $\overline{\Gamma}_i^{(s)}$ is the row vector of length $n^2$ defined by
$$\overline{\Gamma}_i^{(s)}:=(\overline{\Gamma}^{1(s)}_{i1}, \ \overline{\Gamma}^{1(s)}_{i2}, \  \ldots,  \  \overline{\Gamma}^{n(s)}_{i,n-1},  \ \overline{\Gamma}^{n(s)}_{in})\in \textup{Mat}_{1\times n^2}(k)$$
 then, by Equation (\ref{cucucu}), and under the identification (\ref{oprottu}), we have an equality of  $(n^2+1)\times (n^2+1)$-matrices
\begin{equation}
\label{iiiu}
\overline{\nabla}_{\phi_i^{(s)}}=\left(
\begin{array}{ccc}
\delta_i^{(s)}\pi\cdot \textup{Fr}^s &\ &  \overline{\Gamma}_i^{(s)}\cdot \textup{Fr}^s\\
\ & \ & \ \\
0 & \ & 0 
\end{array}
\right)
\end{equation}
Here (and in similar formulae later) we view our $k$-linear spaces as  $R_{\pi}$-modules;
in particular left multiplication of $\textup{Fr}^s$ by elements of $R_{\pi}$ 
such as $\delta_i^{(s)}\pi$
is well defined.

\subsection{Formulae for curvature mod $\pi$} \label{wq2} 
   The  component of lowest degree $c$ the map (\ref{babanicaa}) is  determined as follows. 
 For each $i,j\in \{1,\ldots,n\}$  consider the row vector of length $n^2$ defined by
 \begin{equation}
\label{pushka}
\overline{R}_{ij}:=
(\overline{R}^1_{ij1},\  \overline{R}^1_{ij2},\  \ldots, \  \overline{R}^n_{ij,n-1},\  \overline{R}^n_{ijn})\in \textup{Mat}_{1\times n^2}(k),
\end{equation}
where
\begin{equation}\label{misssing}
\overline{R}^k_{ijl}:= \delta^{(c)}_i \pi \cdot (\overline{\Gamma}^{k(c)}_{jl})^{p^c}-\delta^{(c)}_j \pi \cdot (\overline{\Gamma}^{k(c)}_{il})^{p^c}
-\overline{\Gamma}_{i\star j,l}^{k(2c)}+\overline{\Gamma}_{j\star i,l}^{k(2c)}\in k.
\end{equation}

\

Using \cite[Lem. 2.16]{BM22} and Equation (\ref{staru}) we have
\begin{equation}\label{astro1}
\delta^{(c)}_i \pi\cdot (\delta^{(c)}_j \pi)^{p^c}-\delta^{(c)}_j \pi\cdot (\delta^{(c)}_i \pi)^{p^c}-\delta_{i\star j}^{(2c)}\pi +\delta_{j\star i}^{(2c)}\pi \equiv 0 \ \ \textup{mod}\ \ \pi.\end{equation}

\

Using the equalities  (\ref{misssing}), (\ref{astro1}), (\ref{ppuu}), and (\ref{iiiu}) we get that, under the identification (\ref{oprottu}), we have an equality of
the $(n^2+1)\times (n^2+1)$-matrices with $k$-coefficients:
\begin{equation}
\label{cattt}\overline{\mathcal R}(\phi^{(c)}_i,\phi^{(c)}_j)=\left(
\begin{array}{ccc}
0 & \ & \overline{R}_{ij}\cdot \textup{Fr}^{2c}\\
\ & \ & \ \\
0 & \ & 0 \end{array}
\right).\end{equation}

\

In particular, recalling the arithmetic Riemann curvature tensor in Definition \ref{arct}
 we have that
 $\overline{R}_{ijl}^k\in k=\mathcal A/\mathcal M$ are the images of $R_{ijl}^k\in \mathcal A$.
 Set $q:=q^{(c)}$ and consider the elements $\overline{q}_{ij}\in k$ be defined as
$$\overline{q}_{ij}:=q_{ij}\ \ \textup{mod}\ \ \pi.$$
Then we may consider the images 
$\overline{R}_{ijkl}\in k=\mathcal A/\mathcal M$ of the elements $R_{ijkl}\in \mathcal A$
hence
\begin{equation}
\label{blondiee}\overline{R}_{ijkl}:=\sum_m \overline{R}_{ijk}^m\overline{q}_{ml}^{p^{2c}}\ \ \in k,\end{equation}
for $i,j,k,l\in \{1,\ldots,n\}$. We refer to the collection of elements $\overline{R}_{ijkl}$ as the {\bf reduced Riemann curvature tensor}.

\

Note that we have the following antisymmetry in the first $2$ indices.
\begin{equation}
\label{sy1}
\overline{R}_{ijkl}+\overline{R}_{jikl}=0.\end{equation}
Note finally that for all $i,j,k\in \{1,\ldots, n\}$ we have
\begin{equation}
[\overline{\nabla}_{\phi_k^{(c)}},\overline{\mathcal R}(\phi^{(c)}_i,\phi^{(c)}_j)]=\left(
\begin{array}{ccc}
0 & \ & \delta_i^{(c)}\pi \cdot \overline{R}_{ij}^{p^c}\cdot \textup{Fr}^{3c}\\
\ & \ & \ \\
0 & \ & 0 \end{array}
\right).
\end{equation}
Hence for all $i,j,k$ we have
$$\overline{\mathcal R}(\phi^{(c)}_i,\phi^{(c)}_j)]\neq 0\ \ \Rightarrow\ \ [\overline{\nabla}_{\phi_k^{(c)}},\overline{\mathcal R}(\phi^{(c)}_i,\phi^{(c)}_j)]\neq 0.
$$
In particular, we get the following corollary.

\begin{corollary}\label{coralles}
 Assume there exist $X,Y\in \mathfrak D^{(c)}$ such that $\overline{\mathcal R}(X,Y)\neq 0$. Then  the reduced conormal $\mathfrak d$-module $\overline{\mathfrak n}_{\nabla}$ is non-trivial. So (by Remark \ref{maisuss}, assertion 3) the conormal  $\mathfrak d$-module
 $\mathfrak n_{\nabla}/\mathfrak n^2_{\nabla}$ is also non-trivial.\end{corollary}
 
 {\it Proof}. Indeed, for all $Z\in \mathfrak D^{(c)}$ we have $Z\cdot \overline{\mathcal R}(X,Y)\neq 0$. \qed

\begin{remark}
For $\mathfrak D^{(c)}$ abelian  we have that $i\star j=j\star i$ hence we get 
\begin{equation}\label{misssing2}
\overline{R}^k_{ijl}= \delta^{(c)}_i \pi \cdot (\overline{\Gamma}^{k(c)}_{jl})^{p^c}-\delta^{(c)}_j \pi \cdot (\overline{\Gamma}^{k(c)}_{il})^{p^c}\in k.
\end{equation}
\end{remark}

\

\subsection{Curvature of arithmetic Levi-Civita connection}\label{calcc}
We  place ourselves in the general setting of Section \ref{gencur} and we assume  $N=n$. 
In particular, we have at our disposal the Lie symbols $\ell^{(s)}$ for $s\in c\mathbb N$ attached to the associative symbols $\alpha^{(s)}$.
We  assume in addition that we are given, for each $s\in c\mathbb N$, a torsion symbol $L^{(s)}=(L^{1(s)},\ldots,L^{n(s)})$ of the second kind. For now we do not assume any compatibility between 
the torsion symbol $L^{(s)}$ and the  Lie symbol $\ell^{(s)}$. 

\

Assume in what follows that for every $s\in c\mathbb N$ we are given a metric $q^{(s)}\in \textup{GL}_n(R_{\pi})^{\textup{sym}}$; we can refer to  $q=q^{(c)}$ as the {\bf primary metric} and to $q^{(s)}$ with $s\geq 2c$ as  {\bf  secondary metrics}. Then, for each $s\in c\mathbb N$, one can consider  the arithmetic Levi-Civita connection
$$\Delta^{(s)\textup{LC}}=((\delta_1^{(s)})^{\textup{LC}},\ldots,(\delta_n^{(s)})^{\textup{LC}})$$ 
of degree $s$
attached to  $q^{(s)}\in \textup{GL}_n(R_{\pi})^{\textup{sym}}$  and to the torsion symbol
$L^{(s)}$. We obtain a  graded $\pi$-connection and hence
 an associated curvature $\mathcal R^{\textup{LC}}$ satisfying
$$\mathcal R^{\textup{LC}}(\phi^{(s)}_i,\phi^{(r)}_j)=(\phi_i^{(s)})^{\textup{LC}}(\phi_j^{(r)})^{\textup{LC}} -(\phi_j^{(r)})^{\textup{LC}}(\phi_i^{(s)})^{\textup{LC}} -
\sum_k \ell_{ij}^{k(s,r)}(\phi_k^{(s+r)})^{\textup{LC}}.$$ In particular,
$$\mathcal R^{\textup{LC}}(\phi^{(c)}_i,\phi^{(c)}_j)=(\phi_i^{(c)})^{\textup{LC}}(\phi_j^{(c)})^{\textup{LC}} -(\phi_j^{(c)})^{\textup{LC}}(\phi_i^{(c)})^{\textup{LC}} -
(\phi_{i\star j}^{(2c)})^{\textup{LC}}+(\phi_{i\star j}^{(2c)})^{\textup{LC}}.
$$

\

Later we will show how to choose the torsion symbols and the secondary metrics canonically depending only on the primary metric and the labeling symbols; cf.  Section \ref{ctsm}. 

\medskip

We consider, in what follows, the case $\mathfrak D^{(c)}$ is abelian, equivalently
 $$\ell^{(c)}:=(\ell^{1(c)},\ldots,\ell^{n(c)})=0.$$
   In this case the quantities $\overline{R}_{ijk}^l$
 and hence the components of the arithmetic Riemannian curvature $\overline{R}_{ijkl}$
 only depend on the primary metric $q=q^{(c)}$ and on $L^{(c)}(1)$ (but not
 on the secondary metrics or the torsion symbols $L^{(s)}$ with $s>c$).

\begin{theorem}\label{tery}
Assume $\mathfrak D^{(c)}$ is abelian (equivalently $\ell^{(c)}=0$) and assume in addition that $L^{(c)}(1)=0$. Let $q\in \textup{GL}_n(R_{\pi})^{\textup{sym}}$ and let $\overline{R}_{ijkl}$ be the reduced Riemann curvature tensor attached to $q=q^{(c)}$.
The following equalities hold in $k$:
$$
\overline{R}_{ijkl} =  -\frac{1}{2} \cdot (\delta_{jl}^{(2c)}q_{ik}-
\delta_{jk}^{(2c)}q_{il}-\delta_{il}^{(2c)}q_{jk}+\delta_{ik}^{(2c)}q_{jl}
) \ \ \  \textup{mod}\ \ \pi.$$
\end{theorem}

{\it Proof}. We have:
$$\begin{array}{rcll}
\overline{R}_{ijkl}   & = & \sum_m \overline{R}_{ijk}^m q_{ml}^{p^{2c}}\ \ \  \textup{mod}\ \ \pi & \textup{cf. Def. \ref{arct}}\\
\ & \ &\ &\ \\
\ & = & \sum_m(\delta_i^{(c)} \pi \cdot (\overline{\Gamma}_{jk}^{m(c)})^{p^c}-\delta^{(c)}_j \pi \cdot (\overline{\Gamma}_{ik}^{m(c)})^{p^c})q_{ml}^{p^{2c}}\ \ \  \textup{mod}\ \ \pi & \textup{cf. Eq. \ref{misssing2}}\\
\ & \ & \ & \ \\
\ & = & \delta^{(c)}_i \pi \cdot (\overline{\Gamma}^{(c)}_{jkl})^{p^c}- \delta^{(c)}_j \pi \cdot (\overline{\Gamma}^{(c)}_{ikl})^{p^c}\ \ \  \textup{mod}\ \ \pi & \textup{cf. Eq. \ref{loweringg}}\\
\ & \ & \ & \ \\
\ & = & -\frac{1}{2} \delta^{(c)}_i \pi \cdot (\delta^{(c)}_j q_{kl} +\delta^{(c)}_k q_{jl}-\delta^{(c)}_lq_{jk})^{p^c} & \\
\ & \ & \ & \ \\
\ & \  & + \frac{1}{2} \delta^{(c)}_j \pi \cdot (\delta^{(c)}_i q_{kl} +\delta^{(c)}_k q_{il}-\delta^{(c)}_lq_{ik})^{p^c}  \ \ \  \textup{mod}\ \ \pi & \textup{cf. Thm. \ref{LCC}}\\
\ & \ & \ & \ \\
\ & = & 
 -\frac{1}{2} \cdot (\delta_{jl}^{(2c)}q_{ik}-
\delta_{jk}^{(2c)}q_{il}-\delta_{il}^{(2c)}q_{jk}+\delta_{ik}^{(2c)}q_{jl}
) \ \ \  \textup{mod}\ \ \pi & \ 
\end{array}
$$
where for the last equality we used \cite[Lem. 2.16]{BM22}.
\qed

\begin{remark}
We take the opportunity here to correct a typo in \cite{Bu17}, p. 317: the right hand sides of the last two lines in Equation 7.2 must have a $-$ sign in front. \end{remark}

\begin{remark}
The formula for $\overline{R}_{ijkl}$ in Theorem \ref{tery} is similar to the ``order $2$ part" of the  classical formula for the Riemann curvature tensor; cf. Equation (\ref{zzoo}) in the  Section \ref{diskus}.
\end{remark}

\begin{corollary} \label{antisymm}
 Assume $\ell^{(c)}=L^{(c)}(1)=0$. 
Then for all $i,j,k,l$ we have:
$$
  \begin{array}{rcll}
  \overline{R}_{ijkl} + \overline{R}_{jikl} & = & 0  \ & \ \\
  \ & \ & \ & \ \\
  \overline{R}_{ijkl} + \overline{R}_{ijlk} & = & 0 &  \ \\
  \ & \ & \ & \ \\
  \overline{R}_{ijkl}+\overline{R}_{iklj}+\overline{R}_{iljk} & = & 0 & \ \\
  \ & \ & \ & \  \\
  \overline{R}_{ijkl} - \overline{R}_{klij} & = & 0. & \ 
  \end{array}
 $$
\end{corollary}

{\it Proof}.
The equalities follow  directly  from Theorem \ref{tery}.
By the way the last equality  is also a formal consequence of the first $3$ equalities; cf. the classical tensor manipulations; cf. \cite{CCL00}, p. 142.\qed

\begin{remark}
The symmetries  in Corollary \ref{antisymm}  are, of course, analogues of the classical symmetries of the Riemann tensor in differential geometry.  \end{remark}

\begin{corollary}\label{miarevenit}
Assume $\ell^{(c)}=L^{(c)}(1)=0$ and let $i,j\in \{1,\ldots,n\}$. Assume there exist $k,l\in \{1,\ldots,n\}$   such that
\begin{equation}
\label{madoarecapp}\delta_{jl}^{(2c)}q_{ik}-
\delta_{jk}^{(2c)}q_{il}-\delta_{il}^{(2c)}q_{jk}+\delta_{ik}^{(2c)}q_{jl}
\not\equiv 0\ \ \textup{mod}\ \ \pi.
\end{equation}
Then the following hold:

\

1) $\overline{\mathcal R}(\phi_i^{(c)},\phi_j^{(c)})\neq 0$.

\

2) The classes $\overline{\kappa}_{\nabla},\kappa_{\nabla},h_{\nabla},\overline{h}_{\nabla}$ are non-trivial.

\

3) The $\mathfrak d$-modules $\overline{\mathfrak n}_{\nabla}$ and $\mathfrak n_{\nabla}/ \mathfrak n_{\nabla}^2$ are non-trivial.
\end{corollary}

{\it Proof.} To check assertion 1  assume
 $\overline{\mathcal R}(\phi_i^{(c)},\phi_j^{(c)})= 0$  and seek a contradiction. By Equations (\ref{cattt}) and (\ref{pushka}), we get that $\overline{R}^u_{ijv}=0$ for all indices $u,v$.
Then by Equation (\ref{blondiee}) we get that $\overline{R}_{ijuv}=0$ for all indices $u,v$. We get a contradiction by Theorem \ref{tery}. Assertions 2 and 3  follow from Corollaries \ref{cincisute} and \ref{coralles}. \qed

\begin{remark}\label{sidecenu}
Fix $\mathfrak D^{(c)}$, $L^{(c)}$ with $n\geq 2$, $\ell^{(c)}=L^{(c)}(1)=0$ and allow the primary metric $q=q^{(c)}$ to vary in $\textup{GL}_n(R_{\pi})^{\textup{sym}}$. Then
the condition  that there exist $i,j,k,l$ such that (\ref{madoarecapp}) holds is satisfied
for $q$  ``sufficiently general" (and hence the  non-triviality assertions in assertions 1,2,3 of Corollary \ref{miarevenit} hold for $q$ ``sufficiently general"). 
To see this (and at the same time 
explain the meaning of ``sufficiently general") we proceed as follows.

\

 One can write 
$$q_{ij}\equiv \zeta_{ij0}+\zeta_{ij}\pi \bmod \pi^2$$
with $\zeta_{ij0},\zeta_{ij}\in R$ unique elements each of which is either $0$ or a root of unity in $R$. Define $c_{ij0},c_{ij}\in k$ by
$$c_{ij0}:=\zeta_{ij0} \bmod  \pi,$$
$$c_{ij}:=\zeta_{ij} \bmod \pi.$$
Recall that reduction modulo $p$ gives a bijection between the set consisting of $0$ and the roots of unity  in $R$ and the elements of the field $k$.
So the values of $q$ mod $\pi^2$ are in a natural bijection with the set\begin{equation}
\label{zzu}
\textup{GL}_n(k)^{\textup{sym}}\times \textup{Mat}_n(k)^{\textup{sym}}\end{equation}
via the map
$$(q \bmod \pi^2) \mapsto\ \ ((c_{ij0}),(c_{ij})).
$$

\

\noindent For all $i,j$ define $\epsilon_{ij}\in k^{\times}$ by the formula
$$\epsilon_{ij}:=(\delta_i^{(c)} \pi \cdot (\delta_j^{(c)}\pi)^{p^c} \bmod \pi)^{p^{-2c}}$$
which makes sense because $k$ is perfect. (Note that if $\mathfrak D^{(c)}$ is abelian then $\epsilon_{ij}=\epsilon_{ji}$ for all $i,j$.)
By \cite[Lem. 2.18]{BM22} the condition  (\ref{madoarecapp}) reads
\begin{equation}
\label{nomiz}
\epsilon_{jl}\cdot c_{ik}-
\epsilon_{jk}\cdot c_{il}-
\epsilon_{il}\cdot c_{jk}+
\epsilon_{ik}\cdot c_{jl}\neq 0.
\end{equation}
 So the set of classes mod $\pi^2$ of matrices $q$ that satisfy (\ref{madoarecapp}) for some $i,j,k,l$ is in a natural bijection with  a non-empty Zariski open subset of the set (\ref{zzu}); by the above discussion this subset is in fact equal to the product of $\textup{GL}_n(k)^{\textup{sym}}$ with the complement in the $k$-linear space $\textup{Mat}_n(k)^{\textup{sym}}$ of a proper linear subspace. By the way, the left hand side of Equation (\ref{nomiz}) is the {\bf Kulkarni-Nomizu product} of the symmetric matrices $(\epsilon_{ij})$ and $(c_{ij})$; cf. \cite[p.47]{Be87}.
\end{remark}

\subsection{Curvature of the arithmetic Chern connection}\label{cacc}

Consider again the general setting of Section \ref{gencur}
with $N$ and $n$, this time,   not necessarily equal. Let 
$(q^{(s)})_{s\in c\mathbb N}$ be a sequence of metrics.
One can consider for every $s\in c\mathbb N$ the  Chern connection
$$\Delta^{(s)\textup{Ch}}=((\delta_1^{(s)})^{\textup{Ch}},\ldots,(\delta_n^{(s)})^{\textup{Ch}})$$ 
attached to   $q^{(s)}$. We obtain a  graded $\pi$-connection and hence
 an associated curvature $\mathcal R^{\textup{Ch}}$ satisfying
$$\mathcal R^{\textup{Ch}}(\phi^{(s)}_i,\phi^{(r)}_j)=(\phi_i^{(s)})^{\textup{Ch}}(\phi_j^{(r)})^{\textup{Ch}} -(\phi_j^{(r)})^{\textup{Ch}}(\phi_i^{(s)})^{\textup{Ch}} -
\sum_k \ell_{ij}^{k(s,r)}(\phi_k^{(s+r)})^{\textup{Ch}},$$
in particular,
$$\mathcal R^{\textup{Ch}}(\phi^{(c)}_i,\phi^{(c)}_j)=(\phi_i^{(c)})^{\textup{Ch}}(\phi_j^{(c)})^{\textup{Ch}} -(\phi_j^{(c)})^{\textup{Ch}}(\phi_i^{(c)})^{\textup{Ch}} -
(\phi_{i\star j}^{(2c)})^{\textup{Ch}}+(\phi_{j\star i}^{(2c)})^{\textup{Ch}}.$$

\

Note that for each $i$ the map
 $(\phi_i^{(s)})^{\textup{Ch}}$ only depends on $\phi_i^{(s)}$ and $q^{(s)}$ (and not on the other elements of $\mathfrak D^{(s)}$ or on the labeling symbol $\omega^{(s)}$); also for each $i$ and $j$ the map  $(\phi_{i\star j}^{(2c)})^{\textup{Ch}}$ only depends on $q^{(2c)}$ and on $\phi_{ij}^{(2c)}=\phi_i^{(c)}\phi_j^{(c)}$ (but not on the other elements of $\mathfrak D^{(2c)}$ or on the labeling symbol $\omega^{(2c)}$), cf. Theorem \ref{ch}. 
 
 \

We would like to compute the quantities $\overline{\Gamma}_{jk}^{l(s)}$ and $\overline{R}^l_{ijk}$ in Sections \ref{wq1} and \ref{wq2}, respectively. By \cite[Cor. 3.41]{BM22}, we have the following formula for $s\in c\mathbb N$:
\begin{equation}
\label{bbaa}
\overline{\Gamma}_{jk}^{l(s)}= -\frac{1}{2} \sum_m \delta_j^{(s)}q^{(s)}_{km}(((q^{(s)})^{-1})_{ml})^{p^s}\ \ \textup{mod}\ \ \pi.\end{equation}

\begin{definition}
The metrics $q$ and $q'$ are {\bf conformally equivalent} if there exists $\lambda\in R_{\pi}$ such that $q'=\lambda\cdot q$.
\end{definition}

Conformally equivalent metrics will appear naturally in Subsection \ref{mass}. For conformally equivalent metrics $q:=q^{(c)}$ and $q^{(2c)}$ the Chern curvature map has a particularly simple form:

\begin{theorem} \label{wqq}
If   $q^{(2c)}\equiv \lambda \cdot q^{(c)}$ mod $\pi^2$ for some $\lambda\in R_{\pi}$ then we have the following equalities:
$$\overline{R}^k_{ijl}=\frac{1}{2} \cdot \lambda^{-p^{2c}}(\delta_{ij}^{(2c)}\lambda-\delta_{ji}^{(2c)}\lambda)\cdot \delta_{kl}\ \ \textup{mod}\ \ \pi,$$
where $\delta_{kl}$ is the Kronecker symbol. 
In particular if $q^{(2c)}\equiv  q^{(c)}$ mod $\pi^2$ or if $\mathfrak D^{(c)}$ is abelian
 the arithmetic Chern curvature map $\overline{\mathcal R}$ in Equation (\ref{babanicaa})
vanishes on $\mathfrak D^{(c)}\times \mathfrak D^{(c)}$.\end{theorem}

{\it Proof}. As usual set $q=q^{(c)}$.
By Equations (\ref{misssing}), (\ref{bbaa}), and (\ref{staru}) we have:
\begin{equation}
\label{maue1}
\begin{array}{rcl}
\overline{R}^k_{ijl} & = & -\delta^{(c)}_i \pi (\frac{1}{2}\sum_m \delta^{(c)}_j q_{lm}(q^{-1})_{mk}^{p^c})^{p^c}
+\delta^{(c)}_j \pi (\frac{1}{2}\sum_m \delta^{(c)}_i q_{lm}(q^{-1})_{mk}^{p^c})^{p^c}\\
\ & \ & \ \\
\ & \ & +\frac{1}{2}\sum_m \delta_{ij}^{(2c)} q_{lm}^{(2c)}((q^{(2c)})^{-1})_{mk}^{p^{2c}}-
\frac{1}{2}\sum_m \delta_{ji}^{(2c)} q_{lm}^{(2c)}((q^{(2c)})^{-1})_{mk}^{p^{2c}}\ \ \textup{mod}\ \pi.
\end{array}
\end{equation}
On the other hand we have
\begin{equation}
\label{maue2}
\delta_{ij}^{(2c)}q^{(2c)}_{lm}\equiv \delta_{ij}^{(2c)}(\lambda q_{lm})\equiv \lambda^{p^{2c}}\delta_{ij}^{(2c)}q_{lm}+q_{lm}^{p^{2c}}\delta_{ij}^{(2c)}\lambda\ \ \textup{mod}\ \ \pi
\end{equation}
and \begin{equation}
\label{maue3}
((q^{(2c)})^{-1})_{mk}=\lambda^{-1}(q^{-1})_{mk}.
\end{equation}
Finally recall from \cite[Lem. 2.16]{BM22} that
\begin{equation}
\label{maue4}
\delta^{(c)}_i \pi \cdot (\delta^{(c)}_j q_{lm})^{p^c}\equiv \delta_{ij}^{(2c)} q_{lm}\ \ \textup{mod}\ \ \pi.
\end{equation}
Plugging in  Equations (\ref{maue2}), (\ref{maue3}), (\ref{maue4}) into Equation (\ref{maue1}) we get the equation in the statement of the theorem.
\qed

\section{Canonical  torsion symbols and canonical secondary metrics}\label{ctsm}

In this section we explain how to canonically attach torsion symbols and secondary metrics to associative symbols. Throughout this Section  we assume $N=n$. 

\subsection{Torsion symbols attached to labeling symbols}\label{cts}

\begin{definition}\label{procol}
A tuple $\beta=(\beta^1,\ldots,\beta^n)$ of matrices $\beta^k\in \textup{Mat}_n(R_{\pi})$ is called a {\bf matrix symbol}.
The {\bf Lie symbol} attached to a matrix symbol $\beta$ is the tuple
$\ell:=(\ell^1,\ldots,\ell^n)$ defined by
$$\ell^k_{ij}:=\beta^k_{ji}-\beta^k_{ij}\in R_{\pi}.$$
One can view  $
\ell^k\in \textup{Mat}_n(\widehat{R_{\pi}[y]})$ so one can view $\ell$ as a torsion symbol,
in which case we call $\ell$ the {\bf additive torsion symbol} attached to $\beta$.

\

The {\bf multiplicative torsion symbol} attached to $\beta$ is the torsion symbol
$\ell^*:=(\ell^{1*},\ldots,\ell^{n*})$ with $\ell^{k*}=(\ell^{k*}_{ij})\in \textup{Mat}_n(\widehat{R_{\pi}[y]})$, defined by
$$\ell^{k*}_{ij}:=\beta_{ji}^k-\beta_{ij}^k+\sum_m \beta^k_{jm}((y_m)_{ki}-\delta_{ki})-
\sum_m \beta^k_{im}((y_m)_{kj}-\delta_{kj}).$$
\end{definition}

\begin{remark} In what follows a role will be played by  stochastic matrices.
We recall that a square 
matrix with entries in a field is {\bf  stochastic} if every row  has sum $1$. 
A matrix is called {\bf doubly stochastic} if both the matrix and its transpose are stochastic. 
We  recall, by the way,  the classical Birkchoff-von Neumann theorem \cite{VN} according to which every doubly stochastic matrix with real non-negative coefficients is a convex real linear combination of permutation matrices; this result is, of course,  not directly relevant to our $p$-adic setting. \end{remark}

\begin{remark}\label{kominsky} In the notation above we have:

\

1) $\ell^{k*}(1)=\ell^k$.

\

2) If the matrices $\beta^1,\ldots,\beta^n$ are stochastic (in particular if they are permutation matrices)  then 
$$\ell^{k*}_{ij}:=\beta_{ji}^k-\beta_{ij}^k+\delta_{kj}-\delta_{ki}+\sum_m \beta^k_{jm}(y_m)_{ki}-
\sum_m \beta^k_{im}(y_m)_{kj}.$$
\end{remark}

\

Recall from \cite[Eq. 4.4]{BM22} that for $S$ a $p$-adically complete ring containing $\mathbb Z[\pi]$ we defined a group structure on the set $\textup{Mat}_N(S)$ by
$$a+_{\pi} b:=a+b+\pi ab,\ \ \ a,b\in \textup{Mat}_N(S).$$
As in Part 1, loc.cit. we denote by $\mathfrak g(S)$ the set $\textup{Mat}_N(S)$ equipped with the group law $+_{\pi}$.
This group
 naturally arose from considerations related to the $\pi$-jet spaces of $\textup{GL}_N$.
  
  \

 \begin{definition}
 For any $p$-adically complete $R_{\pi}$-algebra $S$ and any finite multiplicative submonoid $H$ of the ring $\textup{Mat}_n(S)$ (in particular for $H$ a finite subgroup of $\textup{GL}_n(S)$) consider the 
$\mathbb Z_p[\pi]$-linear span $\textup{Span}_{\mathbb Z_p[\pi]}(H)$ of $H$ in $\textup{Mat}_g(S)$. So if  $\mathbb Z_p[\pi][H]$ denotes the $\mathbb Z_p[\pi]$-monoid algebra of $H$ (where $H$ is viewed as acting trivially on $\mathbb Z_p[\pi]$)
then 
$\textup{Span}_{\mathbb Z_p[\pi]}(H)$ is the image of the natural associative $\mathbb Z$-algebra map
 $\mathbb Z_p[\pi][H]\rightarrow \textup{Mat}_n(S)$.

 \

Upon identifying  $\textup{Mat}_g(S)$ and
$\mathfrak g(S)$ as sets, we have that 
$\textup{Span}_{\mathbb Z_p[\pi]}(H)$ is a subgroup 
of $(\mathfrak g(S),+_{\pi})$. 
If $H$ is the group of permutation matrices (so $H$ is the ``Weyl subgroup" of $\textup{GL}_n(S)$) we write  $\mathfrak w(S):=\textup{Span}_{\mathbb Z_p[\pi]}(H)$ and we refer to the latter as the {\bf Weyl subgroup} of $\mathfrak g(S)$. 
Note that every element of $\mathfrak w(S)$ is a $\mathbb Z_p[\pi]$-multiple of a doubly stochastic matrix with entries in $\mathbb Q_p(\pi)$.
\end{definition}

The interpretation of Definition \ref{procol} in terms of symmetry of $\pi$-connections   is given by the following proposition.

\begin{proposition}\label{burta}
Let $\Delta^{(s)G}=((\delta_1^{(s)})^G,\ldots,(\delta_n^{(s)})^G)$ be a $\pi$-connection  of degree $s$  on $G$ with Christoffel symbol of the second kind  $\Gamma^{(s)}=(\Gamma_1^{(s)},\ldots,\Gamma_n^{(s)})$ and associated matrices $\Lambda^{(s)}=(\Lambda_1^{(s)},\ldots, \Lambda_n^{(s)})$. Let $\beta=(\beta^1,\ldots,\beta^n)$ be any matrix symbol. Let $\ell^{(s)}:=(\ell^{1(s)},\ldots,\ell^{n(s)})$ be the  additive torsion symbol attached to $\beta$ and let $\ell^{(s)*}:=(\ell^{1(s)*},\ldots,\ell^{n(s)*})$ be the multiplicative torsion symbol attached to $\beta$  Then the following conditions are equivalent:

\

$1)$ The  $\pi$-connection $\Delta^{(s)G}$ 
is symmetric with respect to $\ell^{(s)}$, i.e.,
$$\Gamma^{k(s)}_{ij}-\Gamma^{k(s)}_{ji}=\ell^{k(s)}_{ij}.$$

\

$2)$ The matrices $\beta^k+\Gamma^{k(s)}$ are symmetric, i.e.,
$$(\beta^k+\Gamma^{k(s)})^t=\beta^k+\Gamma^{k(s)}.$$

\

\noindent Moreover the following are equivalent:

\

$1^*)$ The  $\pi$-connection  $\Delta^{(s)G}$ 
is symmetric with respect to $\ell^{*(s)}$, i.e.,
$$\Gamma^{k(s)}_{ij}-\Gamma^{k(s)}_{ji}= \ell^{k(s)*}_{ij}(\Lambda^{(s)}).$$

\

$2^*)$ The matrices $\beta^k+_{\pi}\Gamma^{k(s)}$ are symmetric, i.e.,
$$(\beta^k+_{\pi}\Gamma^{k(s)})^t=\beta^k+_{\pi}\Gamma^{k(s)}.$$

\end{proposition}

{\it Proof}. The equivalence of $1)$ and $2)$ is trivial. The equivalence of $1^*)$ and $2^*)$
follows by a direct computation using Definition \ref{procol}.
\qed

\begin{remark}
Recall that our  Lie symbol $\ell^{(s)}$ in Definition \ref{ccaann} was  attached to
the associative symbol $\alpha^{(s)}=(\alpha^{1(s)},\ldots,\alpha^{n(s)})$  in Definition \ref{ccaann} which is a matrix symbol. In their turn the associative  symbols $\alpha^{(s)}$ were attached to the labeling symbols $\omega^{(s)}$ and $\omega^{(2s)}$. By analogy with classical differential geometry and in view of the equivalence between $1)$ and $2)$ in Proposition \ref{burta} it is natural to view the symmetry of the matrices
$$\alpha^k+\Gamma^{k(s)}\in \textup{Mat}_n(\mathcal A)$$
as a condition defining the symmetry of our $\pi$-connection of degree $s$; in other words it is natural to take
the torsion symbol $L^{(s)}$ in Subsection \ref{calcc} to be equal to the 
additive torsion symbol $\ell^{(s)}$ attached  to 
the matrix symbol $\alpha^{(s)}$.

\

On the other hand  recall 
that the matrices $\alpha^{k(s)}$ are permutation matrices and hence they  are naturally elements of the Weyl subgroup $\mathfrak w(R_{\pi})$ of $\mathfrak g(R_{\pi})\subset \mathfrak g(\mathcal A)$. The discussion in \cite[Subsect. 4.3]{BM22}
shows that the matrices $\Gamma_i^{(s)}$ and $\Gamma^{k(s)}$ can also be naturally viewed as elements
of the group $\mathfrak g(\mathcal A)$. So we may consider the symmetry of the matrices 
$$\alpha^k+_{\pi}\Gamma^{k(s)}\in \mathfrak g(\mathcal A)$$
 as another natural condition defining the concept of symmetry for the corresponding $\pi$-connection of degree $s$.
In other words, it is also natural to take the torsion symbol $L^{(s)}$ in Subsection \ref{calcc} to be the multiplicative torsion symbol $\ell^{(s)*}$ attached to 
the matrix symbol $\alpha^{(s)}$.\end{remark}

\begin{definition}\label{canonicaltorsion}
 The additive torsion symbol $\ell^{(s)}$ attached to the associative symbol $\alpha^{(s)}$ is called the {\bf additive canonical torsion symbol} (of degree $s$). The multiplicative torsion symbol $\ell^{(s)*}$ attached to the associative symbol $\alpha^{(s)}$ is called the {\bf multiplicative canonical torsion symbol} (of degree $s$). Two torsion symbols will be called {\bf proportional} if one of them is a constant in $R_{\pi}\setminus \{0\}$ times the other. \end{definition}
 
 \begin{remark}\label{agencyy}
 If the Lie symbol $\ell^{(s)}=0$ and $L^{(s)}$ is proportional to either the additive or the multiplicative canonical torsion symbol then $L^{(s)}(1)=0$.
 \end{remark}

\begin{remark}\label{goldmine}
 Let $\mathfrak D^{(c)}$ be an involutive subset of $\mathfrak F^{(c)}$, let $s\in c\mathbb N$,  
 assume $\mathfrak D^{(s)}$ is non-abelian and let 
 $L^{(s)}$ be a torsion symbol  proportional to either the additive or the multiplicative canonical torsion symbol of degree $s$. Then we claim that 
the  trivial  $\pi$-connection  of degree $s$  is not symmetric with respect to $L^{(s)}$. Indeed, for the  trivial $\pi$-connection  of degree $s$  we have $\Gamma_{ij}^{k(s)}=0$ and $\Lambda_i^{(s)}=1$ for all $i,j,k$ so either of the  equations $1)$ or  $1^*)$ in Proposition \ref{burta} 
 reads 
 $$\ell_{ij}^{k(s)}=0$$
 for all $i,j,k$ which is equivalent to $\mathfrak D^{(s)}$ being abelian; this proves our claim.
 \end{remark}

 \subsection{Secondary metrics attached to  labeling symbols}\label{mass}
 For $s\in c\mathbb N$ it is natural to choose the  secondary metrics $q^{(s)}$ with $s\in c\mathbb N\setminus \{c\}$ compatible, in some way, with the primary metric $q:=q^{(c)}$. One natural choice is to take $q^{(s)}=q$ for all $s$. But there are other, more  natural, choices one of which, based on  labeling symbols, we now describe.

\

Let $q=q^{(c)}$ be a primary metric, fix an integer $s\in c\mathbb N\setminus \{c\}$
  and consider words of length $l:=s/c$,
  $$\mu_k^{(s)}:=m_{k1}m_{k2}\ldots m_{kl},\ \ \ k\in \{1,\ldots,n\},\ \ m_{k1},\ldots, m_{kl}\in \{1,\ldots,n\}.$$ 
  Now we take a clue from classical Riemannian geometry where every metric on a bundle $E$ induces metrics on $E^{\otimes s}:=E\otimes\ldots \otimes E$ for $s\geq 2$. 
   This suggests to attach to the $n\times n$ symmetric matrix $q$  the $n\times n$ symmetric matrix $q^{(s)}$ defined by the formula
\begin{equation}
\label{albu}
q_{ij}^{(s)}:=q_{m_{i1} m_{j1}}q_{m_{i2} m_{j2}}
\ldots q_{m_{il} m_{jl}}.\end{equation}

\

In particular, consider
   an involutive subset $\mathfrak D^{(c)}\subset \mathfrak F^{(c)}$
   and consider a coherent labeling symbol on $\mathfrak D$ (equivalently
    a labeling symbol 
   $\omega=\omega^{(c)}:\{1,\ldots,n\}\rightarrow \mathfrak D^{(c)}$ plus a  function
   $\gamma:c \mathbb N\rightarrow \{1,\ldots,n\}$). Then we may consider the words
   $$\mu_k^{(s)}:=\gamma(s-c)\gamma(s-2c)\ldots\gamma(2c)\gamma(c)k\in \mathbb M_n^{(s/c)}$$
   in which case
   $$q_{ij}^{(s)}=q_{\gamma(s-c)\gamma(s-c)}q_{\gamma(s-2c)\gamma(s-2c)}\ldots q_{\gamma(2c)\gamma(2c)}
   q_{\gamma(c)\gamma(c)}q_{ij}.$$
   Note that if 
   \begin{equation}
   \label{qii}
   q_{ii}\in R_{\pi}^{\times}\ \ \textup{for all}\ \ i\in \{1,\ldots,n\}\end{equation}
    then the matrix
   $q^{(s)}$ is invertible hence it is a metric which is conformally equivalent to $q$.
   
   \begin{remark}
    The condition (\ref{qii}) is reminiscent of the positive definite property $q>0$ of metrics $q$ in classical Riemannian geometry.
   \end{remark}
   
   So to every primary metric $q$ satisfying (\ref{qii}), every labeling symbol $\omega=\omega^{(c)}:\{1,\ldots,n\}\rightarrow \mathfrak D^{(c)}$ on $\mathfrak D^{(c)}$,  and every  function $\gamma:c\mathbb N\rightarrow \{1,\ldots,n\}$,  we have attached a family $(q^{(s)})_{s\in c\mathbb N\setminus \{c\}}$ of secondary metrics. One can further specialize to the situation when $\gamma$ is one of the constant functions $\gamma_h$ with image $h\in \{1,\ldots,n\}$ in which case we denote our $q^{(s)}$ by $q^{(s)}_h$, so
   $$(q^{(s)}_h)_{ij}=q_{hh}^{s/c-1} q_{ij}.$$
   
   \begin{definition}
\label{canmet}
For $q$ a primary metric satisfying (\ref{qii}) and $s\in c\mathbb N\setminus\{c\}$
the metrics $q^{(s)}_h$ are called the  {\bf  canonical secondary metrics} attached to $q$, $\omega$, and $h$.   \end{definition}

\section{Invariant polynomials}\label{invpol}
In this section we address the problem of how our curvature tensors change under a ``change of coordinates." 
 The  of ``group of coordinate changes" in our theory will be the permutation group $\Sigma_n$ of the set $\{1,\ldots,n\}$. 
We place ourselves in the setting of Sections \ref{cannn} and \ref{gencur}.

\subsection{Levi-Civita case}
We first consider the situation in Subsection \ref{calcc} and we concentrate on
 the reduced Riemann  curvature tensor $\overline{R}_{ijkl}$  attached to the (additive or multiplicative) canonical 
 torsion symbol (cf. Definition \ref{canonicaltorsion}) or more generally with respect to some $\gamma\in R_{\pi}$ times the (additive or multiplicative) 
 canonical torsion symbol.  
The additive canonical torsion symbol (i.e., the Lie symbol) and the multiplicative canonical torsion symbol are
 completely determined by $q^{(c)},q^{(2c)}$, by $\gamma$, and by  the labeling symbols
 $\omega^{(c)},\omega^{(2c)},\omega^{(4c)}$.
If $\mathfrak D^{(c)}$  is abelian then $\overline{R}_{ijkl}$
only depends on $q:=q^{(c)}$ and on the labeling symbol $\omega:=\omega^{(c)}$ but not
 on $q^{(2c)}$ or on $\gamma$ or on the labeling symbols $\omega^{(2c)}$ or $\omega^{(4c)}$; cf. Theorem \ref{tery} and Remark \ref{kominsky}, assertion 1.
To keep track of the dependence on the matrix $q$ and the labeling symbol $\omega$ we write, in the abelian case: 
$$\overline{R}_{ijkl}=\overline{R}_{ijkl}(q,\omega).$$
For every $\epsilon\in \Sigma_n$ we let $P_{\epsilon}$ be the permutation matrix corresponding to $\epsilon$, hence $(P_{\epsilon})_{ij}=\delta_{\epsilon(i)j}$ and
$(P_{\epsilon}qP_{\epsilon}^t)_{ij}=q_{\epsilon(i)\epsilon(j)}$. We 
may consider the labeling symbol 
$\omega\circ \epsilon$ and we write
$$\epsilon\cdot q:=P_{\epsilon}qP_{\epsilon}^t.$$
We immediately conclude from  Theorem \ref{tery} the following consequence. 

\begin{corollary}\label{sphericall}
 Assume $\mathfrak D^{(c)}$ is abelian. Then  for all $\epsilon\in \Sigma_n$ we have
$$\overline{R}_{ijkl}(\epsilon\cdot q,\omega \circ \epsilon)= \overline{R}_{\epsilon(i)\epsilon(j)\epsilon(k)\epsilon(l)}(q,\omega).$$
\end{corollary}

\

In view of this corollary it is natural to introduce a  ring of invariants as follows.
Consider  $n^4$ indeterminates $X_{ijkl}$ where $i,j,k,l\in \{1,\ldots,n\}$
and a symmetric $n\times n$ matrix of indeterminates $Q=(Q_{ij})$. Then consider the ring
$$\mathbb F_p[X,Q,\det(Q)^{-1}]:=\mathbb F_p[X_{ijkl}, Q_{ij}\ |\ i,j,k,l\in \{1,\ldots,n\}][\det(Q)^{-1}]$$
 and the ideal $\mathfrak a_{\textup{Riem}}$ generated by  all ``symmetries of the Riemann curvature" i.e., by the linear polynomials
 \begin{equation}
  \begin{array}{l}
  X_{ijkl} + X_{jikl},\\
  \ \\
  X_{ijkl} + X_{ijlk},\\
   \ \\
  X_{ijkl}+X_{iklj}+X_{iljk}.
 \end{array}
  \end{equation}

  \

  We recall the classical fact (cf.  \cite{CCL00}, p. 142) that the ideal $\mathfrak a_{\textup{Riem}}$ also contains the polynomials
  $$X_{ijkl}-X_{klij}.$$
  On the ring  $\mathbb F_p[X,Q,\det(Q)^{-1}]$ consider the $\Sigma_n$-action defined by considering for every $\epsilon\in \Sigma_n$ the $\mathbb F_p$-algebra homomorphism
  $$X_{ijkl}\mapsto X_{\epsilon(i)\epsilon(j)\epsilon(k)\epsilon(l)},\ \ Q_{ij}\mapsto Q_{\epsilon(i)\epsilon(j)}.$$
  The ideal $\mathfrak a_{\textup{Riem}}$ is clearly invariant under this action so we have a $\Sigma_n$-action on the  factor ring $\mathbb F_p[X,Q,\det(Q)^{-1}]/\mathfrak a_{\textup{Riem}}$ and hence we may consider the ring of invariants
  $$\mathfrak I_{\textup{Riem}}:=(\mathbb F_p[X,Q,\det(Q)^{-1}]/\mathfrak a_{\textup{Riem}})^{\Sigma_n}.$$
 For every class $[f]\in \mathbb F_p[X,Q,\det(Q)^{-1}]/\mathfrak a_{\textup{Riem}}$ of a polynomial
$$f\in \mathbb F_p[X,Q,\det(Q)^{-1}]$$ 
 it makes sense to evaluate $[f]$ at
 the tuple $(\ldots,q_{ij},\ldots,\overline{R}_{ijkl},\ldots)$ where $\overline{R}_{ijkl}$ are the components of the reduced Riemannian curvature attached to $q=(q_{ij})$ to get an element
 $$[f](\ldots, q_{ij}, \ldots,\overline{R}_{ijkl},\ldots)  \in \ k:$$
 indeed the value $f(\ldots, q_{ij},\ldots,\overline{R}_{ijkl},\ldots)$ only depends on $[f]$ because the elements of the ideal $\mathfrak a_{\textup{Riem}}$ vanish at $(\overline{R}_{ijkl})$. Note that if $[f]\in \mathfrak I_{\textup{Riem}}$ then one has $f-\epsilon f\in \mathfrak a_{\textup{Riem}}$.
  The ring of invariants $\mathfrak I_{\textup{Riem}}$ is a finitely generated $\mathbb F_p$-algebra of the same Krull dimension as $\mathbb F_p[X,Q,\det(Q)^{-1}]/\mathfrak a_{\textup{Riem}}$.

 \begin{definition}
 For $I\in \mathfrak I_{\textup{Riem}}$.
 the element 
 \begin{equation}
 \label{scalar}
I(q,\overline{R}):= I(\ldots, q_{ij},\ldots, \overline{R}_{ijkl}(q,\omega),\ldots)\in k,\ \ \ s\in \{1,\ldots,m\},\end{equation}
  is called the {\bf $I$-scalar curvature} of $(q,\omega)$. \end{definition}

 In view of Corollary \ref{sphericall}, we get the following invariance:
 
 \begin{corollary} Assume $\mathfrak D^{(c)}$ is abelian.
 For all $q$, all $\omega$, 
 all $\epsilon\in \Sigma_n$, and all $I\in \mathfrak I_{\textup{Riem}}$ we have
$$I(\epsilon\cdot q,\overline{R}(\epsilon\cdot q,\omega\circ \epsilon)=
I(q,\overline{R}(q,\omega)).$$
 \end{corollary}
 
   \begin{example}
  Note that the  coefficients of the characteristic polynomial of $Q$ are elements
  of $\mathbb F_p[Q,\det(Q)^{-1}]^{\Sigma_n}$ and hence of
   $\mathfrak I_{\textup{Riem}}$. But, of course, the ring $\mathbb F_p[Q,\det(Q)^{-1}]^{\Sigma_n}$
  is  larger than the algebra generated by the the  coefficients of the characteristic polynomial of $Q$. In particular, one can trivially show that the $\mathbb F_p$-linear space
  of $\Sigma_n$-invariant elements in the vector space $\sum_{ij} \mathbb F_p Q_{ij}$ has a basis consisting of the following two polynomials:
  $$\sum_i Q_{ii}\ \ \textup{and}\ \  \sum_{i,j} Q_{ij}.$$

  \

  Similarly we have the following $\Sigma_n$-invariant polynomials
  $$\sum_{i, j} X_{ijij}\ \ \textup{and}\ \  \sum_{i,j,k} X_{ijik};$$
  hence their images
   in $\mathbb F_p[X,Q,\det(Q)^{-1}]/\mathfrak a_{\textup{Riem}}$ belong to $\mathfrak I_{\textup{Riem}}$.
  
  \

  More generally, for an integer $d\geq 1$,  let us say that a subset $D\subset \{1,\ldots,n\}^d$ is a {\bf principal diagonal} if there exist two distinct integers $s,t\in \{1,\ldots,d\}$ such that
  $$D=\{(i_1,\ldots,i_d)\in \{1,\ldots,n\}^d\ |\ i_s=i_t\}.$$
  A subset $D \subset \{1,\ldots,n\}^d$ will be called  a {\bf  diagonal} if it is either equal to 
  $\{1,\ldots,n\}^d$ or it is an intersection of principal diagonals. 
  Finally define $Q^{ij}$ to be the components of the inverse of $Q=(Q_{ij})$.
  Here are some basic examples of $\Sigma_n$-invariant polynomials linear in the $X$ indeterminates:
  
  \

  1) For every diagonal $D\subset \{1,\ldots,n\}^6$  the polynomials
  $$\sum_{(i,j,k,l,m,u)\in D} Q_{ij}X_{klmu},\ \ \sum_{(i,j,k,l,m,u)\in D} Q^{ij}X_{klmu}$$
  are $\Sigma_n$-invariant hence their images
   in $\mathbb F_p[X,Q,\det(Q)^{-1}]/\mathfrak a_{\textup{Riem}}$ belong to $\mathfrak I_{\textup{Riem}}$. 
  
  \

   2) Similarly for every diagonal $D\subset \{1,\ldots,n\}^8$  the polynomials
  $$\sum_{(i,j,k,l,m,u,v,w)\in D} Q_{ij}Q_{kl}X_{muvw},\ \ \sum_{(i,j,k,l,m,u,v,w)\in D} Q^{ij}Q_{kl}X_{muvw},$$
  $$ \sum_{(i,j,k,l,m,u,v,w)\in D} Q^{ij}Q^{kl}X_{muvw}$$
  are $\Sigma_n$-invariant hence their images
   in $\mathbb F_p[X,Q,\det(Q)^{-1}]/\mathfrak a_{\textup{Riem}}$ belong to $\mathfrak I_{\textup{Riem}}$.
    In particular, let $S\in \mathfrak I_{\textup{Riem}}$ be the image of the polynomial
    $$\sum_{i,j,k,l} Q^{ij}Q^{lk}X_{likj}.$$ 
    Then the $S$-scalar curvature is given by
  $$S(\ldots,q_{ij},\ldots,\overline{R}_{ijkl},\ldots)=\sum_{i,j,k,l}q^{ij}q^{lk} \overline{R}_{likj},$$
  where $(q^{ij})$ is the inverse of $(q_{ij})$. 
  This is a direct analogue of the classical scalar curvature; cf Equation (\ref{rws}) in the Section
  \ref{diskus}. 
   
   \

   3) The above examples can be, of course, generalized to the case when the degree in the $Q$ indeterminates is arbitrary. It would be interesting to explicitly determine a minimal set of generators for the ring     
   $\mathfrak I_{\textup{Riem}}$. \end{example}

  \begin{example}
  For $n=2$ we have
  $$\mathbb F_p[X,Q,\det(Q)^{-1}]/\mathfrak a_{\textup{Riem}}\simeq \mathbb F_p[X_{1212}.
  Q_{11},Q_{12},Q_{22}, \det(Q)^{-1}],$$
  The element $X_{1212}$ is in $\mathfrak I_{\textup{Riem}}$ 
and the $X_{1212}$-scalar curvature is given by
$$\overline{R}_{1212}=-\frac{1}{2}(\delta^{(2)}_{22}q_{11}-2\delta^{(2)}_{12}q_{12}+\delta_{11}^{(2)}q_{22})\ \ \textup{mod}\ \ \pi.$$
This quantity is of course non-zero for $q$ ``sufficiently general."
  \end{example}

 \subsection{Chern case}
 We next consider the curvature of the Chern connection; cf. Subsection \ref{cacc}.
 To simplify our discussion we consider the case $N=n$, $c=1$,  a  primary metric 
 $$q=q^{(1)}\in \textup{GL}_n(R_{\pi})^{\textup{sym}}$$
 all of whose diagonal entries are invertible
 and a  labeling symbol 
 $$\omega=\omega^{(1)}:\{1,\ldots,n\}\rightarrow \mathfrak D^{(1)}.$$
 Fix $h\in \{1,\ldots,n\}$.
 We may consider 
 the canonical labeling symbol
  $$\omega^{(2)}=\omega_h^{(2)}:\{1,\ldots,n\}\rightarrow \mathfrak D^{(2)}$$
  attached to $\omega$ and $h$; 
  cf. Definition \ref{cfff}.
  Hence
  $$\omega^{(2)}(i)=\phi^{(2)}_i=\phi_h\phi_i.$$
   Furthermore we can consider the canonical secondary metric
   $q^{(2)}=q^{(2)}_h$ attached to $q,\omega,\omega^{(2)}$; cf. Definition \ref{canmet}. Hence
  $$q^{(2)}=q_{hh}\cdot q.$$
  For each $h$ the elements 
  $$\overline{R}^l_{ijk}=\overline{R}^l_{ijk}(q,\omega,h)\in k$$
  in Theorem \ref{wqq} are given by:
  $$\overline{R}^l_{ijk}(q,\omega,h)=\frac{1}{2} \cdot q_{hh}^{-p^2}(\delta_{ij}^{(2)}q_{hh}-\delta_{ji}^{(2)}q_{hh})\cdot \delta_{kl}\ \ \textup{mod}\ \ \pi.$$
 
\begin{definition}
The {\bf arithmetic Chern curvature tensor} is the collection of elements $\overline{F}_{hij}\in k$,
$$\overline{F}_{hij}:= \overline{F}_{hij}(q,\omega):=\frac{1}{2} \cdot q_{hh}^{-p^2}(\delta_{ij}^{(2)}q_{hh}-\delta_{ji}^{(2)}q_{hh})\ \ \textup{mod}\ \ \pi,$$
with $h,i,j\in \{1,\ldots,n\}$.
\end{definition}

\

As in the Levi-Civita case for every $\epsilon\in \Sigma_n$ we let $P_{\epsilon}$ be the permutation matrix corresponding to $\epsilon$ so
$\epsilon\cdot q:=P_{\epsilon}qP_{\epsilon}^t$; this matrix has again the property that all its diagonal entries are invertible.
 Also for $\epsilon \in \Sigma_n$ we
may consider the labeling symbol 
$\omega\circ \epsilon$. Then 
 for all $\epsilon \in \Sigma_n$  we have
\begin{equation}
\label{soluu}
\overline{F}_{hij}(\epsilon\cdot q,\omega \circ \epsilon)= \overline{F}_{\epsilon(h)\epsilon(i)\epsilon(j)}(q,\omega).\end{equation}

\

In view of this latter equation it is natural to introduce a ring of invariants as follows.
Consider first the ring of polynomials in $n^3$ indeterminates
$$\mathbb F_p[Y,Q,\det(Q)^{-1}]:=\mathbb F_p[Y_{hij},Q_{ij},\det(Q)^{-1}\ |\ h\in \{1,\ldots,n\}, \ i,j\in \{1,\ldots,n\}]$$ and the ideal $\mathfrak a_{\textup{Ch}}$ generated by  the linear polynomials
$$Y_{hij}+Y_{hji}.$$
  On the ring $\mathbb F_p[Y,Q,\det(Q)^{-1}]$ consider the $\Sigma_n$-action defined by considering for every $\epsilon\in\Sigma_n$ the $\mathbb F_p$-algebra homomorphism given by
  $$Y_{hij}\mapsto Y_{\epsilon(h)\epsilon(i)\epsilon(j)},\ 
  \ Q_{ij}\mapsto Q_{\epsilon(i)\epsilon(j)}.$$

  \
  
  The ideal $\mathfrak a_{\textup{Ch}}$ is clearly invariant under this action so we have a $\Sigma_n$-action on the  factor ring and hence we may consider the ring of invariants
  $$\mathfrak I_{\textup{Ch}}:=(\mathbb F_p[Y,Q,\det(Q)^{-1}]/\mathfrak a_{\textup{Ch}})^{\Sigma_n}.$$
 For every class $[f]\in \mathbb F_p[Y,Q,\det(Q)^{-1}]/\mathfrak a_{\textup{Ch}}$ of a polynomial
$f\in \mathbb F_p[Y]$ 
 it makes sense to evaluate $[f]$ at
 the components of the arithmetic Chern curvature tensor $\overline{F}_{hij}$ to get an element
 $$[f](\ldots,q_{ij},\ldots,\overline{F}_{hij},\ldots)  \in \ k.$$
  Note that if $[f]\in \mathfrak I_{\textup{Ch}}$ then one has $f-\epsilon f\in \mathfrak a_{\textup{Ch}}$.
  The ring of invariants $\mathfrak I_{\textup{Ch}}$ is a finitely generated $\mathbb F_p$-algebra of the same Krull dimension as $\mathbb F_p[Y,Q,\det(Q)^{-1}]/\mathfrak a_{\textup{Ch}}$.

 \begin{definition}
 For $I\in \mathfrak I_{\textup{Ch}}$.
 the element 
 \begin{equation}
 \label{scalar2}
I(q,\overline{F}(q,\omega)):= I(\ldots, q_{ij},\ldots, \overline{F}_{hij}(q,\omega),\ldots)\in k,\ \ \ s\in \{1,\ldots,m\},\end{equation}
  is called the {\bf $I$-scalar curvature} of $(q,\omega)$. \end{definition}

 In view of by Equation (\ref{soluu}), we get the following invariance:
 
  \begin{corollary} For all metrics $q$ with invertible diagonal entries, all labeling symbols $\omega$, 
 all $\epsilon\in \Sigma_n$,
 and all $I\in \mathfrak I_{\textup{Ch}}$ we have
$$I(\epsilon\cdot q,\overline{F}(\epsilon\cdot q,\omega\circ \epsilon))=
I(q,\overline{F}(q,\omega)).$$
 \end{corollary}
 
 \begin{example}
  We have the following $\Sigma_n$-invariant polynomial
  $$\sum_{i, j} Y_{iij};$$
  hence its image $I$
   in $\mathbb F_p[Y,Q,\det(Q)^{-1}]/\mathfrak a_{\textup{Ch}}$ belongs to $\mathfrak I_{\textup{Ch}}$ and satisfies
   $$I(q,\overline{F}(q,\omega))=\sum_{i,j}\overline{F}_{iij}.$$
   The invariant $I$ can be viewed as an analogue of the trace of the curvature of the Chern connection in classical differential geometry.
 As in the Levi-Civita case one can consider, more generally, $\Sigma_n$-invariant elements in  $\mathbb F_p[Y,Q,\det(Q)^{-1}]/\mathfrak a_{\textup{Ch}}$ involving the variables $Q_{ij}$
 and it would be interesting to find a minimal set of generators for $\mathfrak I_{\textup{Ch}}$. \end{example}

 \section{Multiplicative curvature}
 
 The curvature of a  graded $\pi$-connection  introduced and studied in the previous sections
 should be viewed as an ``additive" object and has a ``multiplicative" analogue which we explain in what follows. 
  
 \subsection{Abstract setting} \label{abbs}
 We begin by describing an abstract setting which will later be specialized to our $\pi$-connection  context. 
 
 \begin{definition}\label{simonn1}
 Let $\mathfrak G$ be a group with operation denoted by $\circ$, inverse denoted by $(\  )^{\circ(-1)}$, identity denoted by $e$,
 and elements denoted by $X,Y,\ldots$. By a {\bf monoid over} $\mathfrak G$ we understand a monoid $\mathfrak E$ with operation denoted by $\circ$ and with identity denoted by $e$, equipped with a homomorphism of monoids with identity
 \begin{equation}
 \label{ramonn}
 \mathfrak E\rightarrow \mathfrak G,\ \ \ f\mapsto f_*.
 \end{equation}
 Denote by $\mathfrak R$ the preimage in $\mathfrak E$ of $e\in \mathfrak G$ under (\ref{ramonn}). We say that $\mathfrak E$ is a {\bf split monoid} if in addition to the above data  
 one is given a 
   homomorphism of monoids with identity
 \begin{equation}
 \label{ramonna}
 \mathfrak G\rightarrow \mathfrak E,\ \ \ X\mapsto X^*
 \end{equation}
 which is left inverse to (\ref{ramonn}) such that the map
 \begin{equation}
 \label{ramonnel}
 \mathfrak G\times \mathfrak R\rightarrow \mathfrak E,\ \ \ (X,g)\mapsto g\circ X^*
 \end{equation}
 is a bijection. We denote the inverse of (\ref{ramonnel}) by
 \begin{equation}
 \label{ramonnela}
 \mathfrak E\rightarrow \mathfrak G\times \mathfrak R,\ \ \ f\mapsto (f_*,\widetilde{f}).
 \end{equation}
 \end{definition}
 
 \begin{remark}
 Given a split monoid $\mathfrak E$ over a group $\mathfrak G$ there is a natural left action of $\mathfrak G$ on $\mathfrak R$ defined by the map
 \begin{equation}
 \mathfrak G\times \mathfrak R\rightarrow \mathfrak R,\ \ \ (X,g)\mapsto
 g^{(X)}:=X^*\circ g \circ (X^{\circ(-1)})^*. \end{equation}
 The monoid operation $\circ$  on $\mathfrak G\times \mathfrak R$
 induced from the operation $\circ $ on  $\mathfrak E$ via the bijection (\ref{ramonnel}) is then given by the {\bf semidirect product} rule
 $$(X_1,g_1) \circ (X_2,g_2)=(X_1\circ X_2,g_1\circ g_2^{(X_1)}).$$
 So for all $f_1,f_2\in \mathfrak E$ we have
 $$\widetilde{f_1\circ f_2}=\widetilde{f_1}\circ (\widetilde{f_2})^{(f_{1*})}.$$
 \end{remark}
 
 \begin{definition}\label{simonn2}
 Assume $\mathfrak E$ is a split monoid over a group $\mathfrak G$. By a {\bf compatible group structure} on $\mathfrak R$ we understand a group structure on $\mathfrak R$ with operation denoted by $\cdot$, inverse denoted by $(\ )^{-1}$, and identity denoted by $1$ such that the following properties hold:

 \

 1) ({\bf left distributivity}) For all $g_1,g_2,g_3\in \mathfrak R$ we have
 $$g_1\circ (g_2\cdot g_3)=(g_1\circ g_2)\cdot (g_1\circ g_3).$$

\

 2) For all $X\in \mathfrak G$ and all $g_1,g_2\in \mathfrak E$ we have
 $$(g_1\cdot g_2)^{(X)}=g_1^{(X)}\cdot g_2^{(X)}.$$\end{definition}

\

Denote by $\mathfrak E\times_{\mathfrak G}\mathfrak E$  the fiber product i.e., the set of all pairs $(f_1,f_2)\in \mathfrak E\times \mathfrak E$ such that $f_{1*}=f_{2*}$.
 Given a compatible group structure as above there are unique operations 
 $$\mathfrak E\times_{\mathfrak G}\mathfrak E\rightarrow \mathfrak E,\ \ (f_1,f_2)\mapsto f_1\cdot f_2$$
$$\mathfrak E\rightarrow \mathfrak E,\ \ \ f\mapsto f^{-1}$$
satisfying the conditions
$$\widetilde{f_1\cdot f_2}=\widetilde{f_1}\cdot \widetilde{f_2},\ \ \ \widetilde{f^{-1}}=(\widetilde{f})^{-1}.$$
One immediately checks that we again have the left distributivity property:
$$f_1\circ (f_2\cdot f_3)=(f_1\circ f_2)\cdot (f_1\circ f_3)\ \ \textup{for}\ \  (f_1,(f_2,f_3))\in \mathfrak E\times(\mathfrak E\times_{\mathfrak G}\mathfrak E).$$
Also note the formula:
$$f_1\circ (f_2^{-1})=(f_1\circ f_2)^{-1},\ \ f_1,f_2\in \mathfrak E.$$

\begin{remark}
The set $\mathfrak R$ equipped with the operations $\circ$ and $\cdot$ is a near-ring in the sense of \cite{C92}. However the set $\mathfrak E$ with the operations $\circ$ and $\cdot$ is not a near-ring, as the operation $\cdot$ on $\mathfrak E$ is only defined on $\mathfrak E\times_{\mathfrak G}\mathfrak E$ and not on $\mathfrak E\times\mathfrak E$.
\end{remark}

  Consider now a submonoid $\mathfrak G'$ of $\mathfrak G$ 
 and let $\mathfrak E'$ be the preimage of $\mathfrak G'$ under the projection $\mathfrak E\rightarrow \mathfrak G$. 
 Note that $\mathfrak E'$ is a submonoid of $(\mathfrak E,\circ)$ and it has an induced
 operation $\cdot :\mathfrak E'\times_{\mathfrak G'}\mathfrak E'\rightarrow \mathfrak E'$ (with respect to which $\circ$ is, again, left distributive) and an induced ``inverse" operation $\mathfrak E'\rightarrow \mathfrak E'$.
 Consider 
  a section $\nabla:\mathfrak G'\rightarrow \mathfrak E'$, $X\mapsto \nabla_X$, of the projection $\mathfrak E'\rightarrow \mathfrak G'$. So for $X,Y\in \mathfrak G'$ we have the formula
  $$\widetilde{\nabla_X \circ \nabla_Y}=\widetilde{\nabla_X}\circ (\widetilde{\nabla_Y})^{(X)}.$$

\

\begin{definition}\label{ohnoy}
For every section $\nabla:\mathfrak G'\rightarrow \mathfrak E'$, $X\mapsto \nabla_X$, of the projection $\mathfrak E'\rightarrow \mathfrak G'$ 
 we consider the maps
  $$
  \begin{array}{ll}
  \mathcal S^*:\mathfrak G'\times \mathfrak G'\rightarrow
 \mathfrak E', & \mathcal S^*(X,Y):=(\nabla_X\circ \nabla_Y)
 \cdot (\nabla_{X\circ Y})^{-1},\\
 \ & \ \\
 \mathcal R^*:\mathfrak G'\times \mathfrak G'\rightarrow
 \mathfrak R, &  \mathcal R^*(X,Y)  :=  \widetilde{\mathcal S^*(X,Y)}\cdot (\widetilde{\mathcal S^*(Y,X)})^{-1}.\end{array}$$
\end{definition}

Note that $\mathcal S^*$ is well defined because 
$(\nabla_X \circ \nabla_Y,\nabla_{X\circ Y})\in \mathfrak E'\times_{\mathfrak G'}\mathfrak E'$. 

\

\begin{remark} The following formulae hold:

 $$
 \begin{array}{rcl}
 \widetilde{\mathcal S^*(X,Y)}  & := &  (\widetilde{\nabla_X \circ \nabla_Y})\cdot (\widetilde{\nabla_{X\circ Y}})^{-1}=(\widetilde{\nabla_X}\circ (\widetilde{\nabla_Y})^{(X)})\cdot 
 (\widetilde{\nabla_{X\circ Y}})^{-1}\\
 \ & \ & \ \\
 \mathcal R^*(X,Y)&=&  (\widetilde{\nabla_X\circ \nabla_Y}) \cdot 
  (\widetilde{\nabla_{X\circ Y}})^{-1}
  \cdot (\widetilde{\nabla_{Y\circ X}})
  \cdot
 (\widetilde{\nabla_Y\circ \nabla_X})^{-1}.
 \end{array}$$
\end{remark}

\begin{remark} 
Assume $X,Y\in \mathfrak G'$. Then the following are equivalent:

\

i) $\widetilde{\mathcal S^*(X,Y)}=1$.

\

ii) $\nabla_{X\circ Y}=\nabla_X\circ \nabla_Y$.

\end{remark}

\begin{remark}\label{conditto}
The following conditions are equivalent:

\

i) $\widetilde{\mathcal S^*(X,Y)}=1$ for all $X,Y\in \mathfrak G'$.

\

ii) $\nabla$ is a monoid homomorphism from $(\mathfrak G',\circ)$ to $(\mathfrak E',\circ)$.

\

iii) $X\mapsto \widetilde{\nabla_X}$ is a  $1$-cocycle from $(\mathfrak G',\circ)$ to $(\mathfrak R,\circ)$.
\end{remark}

\begin{remark}
Assume $X,Y\in \mathfrak G'$ are such that $X\circ Y=Y\circ X$. Then the following formula holds:
$$ \mathcal R^*(X,Y) =
  (\widetilde{\nabla_X\circ \nabla_Y})
  \cdot
 (\widetilde{\nabla_Y\circ \nabla_X})^{-1}.$$
Hence if $X\circ Y=Y\circ X$ the following are equivalent:

\smallskip

i) $\mathcal R^*(X,Y)=1$.

\smallskip

ii) $\nabla_X \circ \nabla_Y=\nabla_Y \circ \nabla_X$.

\end{remark}

In our applications the map $\nabla$ will arise from a  graded $\pi$-connection while $\mathcal R^*$ will correspond to the ``multiplicative curvature" of the  graded $\pi$-connection (see Definition \ref{eei} below).

 \subsection{Relative endomorphisms}
 As an intermediate step between the abstract setting described above and our application
 to the context of  $\pi$-connections we describe a formalism of ``relative endomorphisms."
 
 \

 Assume we are given a  $\mathbb Z_p$-algebra $A_{\mathbb Z_p}$ and let $A:=A_{\mathbb Z_p}\otimes_{\mathbb Z_p} R_{\pi}$. 
 The sets 
 $$\textup{End}_{\mathbb Z_p-\textup{alg}}(\widehat{A})\ \ \ \textup{and}\ \ \ 
 \mathfrak G:= \textup{Aut}_{\mathbb Z_p-\textup{alg}}(R_{\pi})$$
  have structures of monoids (where the second is, of course, a group) under  composition. In this subsection composition will be denoted 
  by ``$\circ$" rather than  by juxtaposition in order to avoid confusion with matrix multiplication which will appear later and which will be denoted by ``$\cdot$"; we denote by $e$ the identity maps in these monoids; also we will denote by $\phi,\phi_1,\phi_2,\ldots$ the elements of $\mathfrak G$ and the inverse operation in $\mathfrak G$ will be denoted by $\phi\mapsto \phi^{\circ (-1)}$ in order to later avoid confusion with the inverse operation on matrices.

  \

    We let $\mathfrak E$ be the set of all 
  $f\in  \textup{End}_{\mathbb Z_p-\textup{alg}}(\widehat{A})$ that send $R_{\pi}$ onto itself.
    For every $f\in \mathfrak E$ we denote by 
 $f_*\in \mathfrak G$ the restriction of $f$ to $R_{\pi}$;
 we get a  monoid homomorphism
 \begin{equation}
 \label{mona1}
 \mathfrak E\rightarrow \mathfrak G,\ \ f\mapsto f_*;\end{equation}
 hence the preimage of the identity $e\in \mathfrak G$ in $\mathfrak E$ equals
  $$\mathfrak R:= 
   \textup{End}_{R_{\pi}-\textup{alg}}(\widehat{A})\subset \mathfrak E.$$

   \

 \noindent On the other hand for every $\mathbb Z_p$-algebra automorphism $\phi\in \mathfrak G$, $\phi:R_{\pi}\rightarrow R_{\pi}$, we have an induced endomorphism $\phi^*:\widehat{A}\rightarrow \widehat{A}$,
 $\phi^*:=\widehat{1\otimes \phi}\in\mathfrak E$.
The map
 \begin{equation}
 \label{mona2}
 \mathfrak G\rightarrow \mathfrak E,\ \ \phi\mapsto \phi^*
 \end{equation}
 is a monoid homomorphism which is a right inverse for (\ref{mona1}).
 Clearly the map
 $$\mathfrak R\times \mathfrak G\rightarrow \mathfrak E,\ \ \ (\phi,g)\mapsto g\circ \phi^*$$
 is a bijection whose inverse we denote by 
 $$f\rightarrow (f_*,\widetilde{f}).$$
  So $\mathfrak E$, equipped with the data above, is a split monoid over $\mathfrak G$
  in the sense of Definition \ref{simonn1}.  For every $f\in \mathfrak E$ the homomorphism $\widetilde{f}\in \mathfrak R$
 can be viewed as a ``relative" version of $f$; the homomorphisms $f$ and $\widetilde{f}$ coincide on $A_{\mathbb Z_p}$ but $f$ is only a $\mathbb Z_p$-algebra map whereas $\widetilde{f}$ is an $R_{\pi}$-algebra map.
 
 \
 
  {\it  From now on we let $A_{\mathbb Z_p}:=\mathbb Z_p[x,\det(x)^{-1}]$.} 
  
  \

  \noindent So we have  
  $$A=R_{\pi}[x,\det(x)^{-1}]=\mathcal O(G),\ \ \ \widehat{A}=\mathcal A=\widehat{\mathcal O(G)}.$$
 For every $F\in \mathcal A$ and every $\phi\in \mathfrak G$  we write $F^{(\phi)}:=\phi^* (F)$. 
  We adopt a similar notation for matrices with entries in $\mathcal A$. 
  If for $i\in \{1,2\}$ we consider $f_i\in \mathfrak E$ and we write
  \begin{equation}
  \label{lkk}
   f_i(x)=F_i:=F_i(x)\in \textup{GL}_n(\mathcal A),\ \ \ \phi_i:=f_{i*}\end{equation}
   then
 \begin{equation}
 \label{binee}
 (\widetilde{f_1\circ f_2})(x)=(f_1\circ f_2)(x)=
 f_1(F_2(x))=F_2^{(\phi_1)}(f_1(x))=F_2^{(\phi_1)}(F_1(x)).\end{equation}
On the other hand
\begin{equation}
\label{bineee}
(\widetilde{f_1}\circ \widetilde{f_2})(x)=\widetilde{f_1}(\widetilde{f_2}(x))=\widetilde{f_1}(f_2(x))=
\widetilde{f_1}(F_2(x))=F_2(\widetilde{f_1}(x))=F_2(F_1(x)).\end{equation}
Moreover, for $g\in \mathfrak R$, $F:=g(x)$, $\phi\in \mathfrak G$, we have
\begin{equation}\label{brrh}
g^{(\phi)}(x)=\phi^*(F(x))=F^{(\phi)}.
\end{equation}

\

In addition to carrying a monoid structure defined by composition, the set
of endomorphisms $\mathfrak R$ carries, in our case here, a group structure coming from that of $G$. Indeed we have a bijection
  \begin{equation}
  \label{biji}
  \mathfrak R\simeq G(\mathcal A),\ \ g\mapsto g(x)\end{equation}
  where $G(\mathcal A)=\textup{GL}_n(\mathcal A)$ is the set of
   $\mathcal A$-points of the $R_{\pi}$-group scheme $G$;
   now $G(\mathcal A)$ has a group structure defined by that of $G$ so  by transport of structure we get a group structure on 
   $\mathfrak R$ 
   whose multiplication operation will be denoted by ``$\cdot$", whose identity element will be denoted by $1$,  and whose inverse operation will be denoted by 
   $(\ )^{-1}$. Note that, with the notation in Equation (\ref{lkk}) we have
\begin{equation}
(\widetilde{f_1}\cdot \widetilde{f_2})(x)=
\widetilde{f_1}(x)\cdot \widetilde{f_2}(x)=F_1(x)\cdot F_2(x).
\end{equation}

\

Moreover the two operations on $\mathfrak R$ satisfy the following compatibility relation: for all $g_1,g_2,g_3\in \mathfrak R$ we have
\begin{equation}
\label{douaop}
g_1\circ (g_2\cdot g_3)=(g_1\circ g_2)\cdot (g_1\circ g_3).\end{equation}
   Finally note that the action of $\mathfrak G$ on the set $\mathfrak R$ also respects the group operation ``$\cdot$" i.e., for all $g_1,g_2\in \mathfrak R$ and all $\phi\in \mathfrak G$ we have
   \begin{equation}\label{sechea2}
   (g_1\cdot g_2)^{(\phi)}=g_1^{(\phi)}\cdot g_2^{(\phi)},
   \end{equation}
   as one readily checks using Equation (\ref{brrh}).
So the data above defines a compatible group structure on $\mathfrak R$ in the sense of Definition \ref{simonn2}.  Consequently we may consider various submonoids $\mathfrak G'$ of $\mathfrak G$, various sections $\nabla:\mathfrak G'\rightarrow \mathfrak E'$ of the projection $\mathfrak E'\rightarrow \mathfrak G'$, and the corresponding maps $\mathcal S^*, \mathcal R^*$ in Definition \ref{ohnoy}. We will specialize 
this construction to the case of $\pi$-connections   in the next subsection. 

\

 \subsection{Definition  of multiplicative curvature and computation mod $\pi$}
 Assume we are given an involutive subset $\mathfrak D^{(c)}\subset \mathfrak F^{(c)}$ and labeling symbols
 $\omega^{(s)}$ for $s\in c\mathbb N$ as in Subsection \ref{cannn}. We have at our disposal  the higher $\pi$-Frobenius lifts $\phi_i^{(s)}$  of degree $s$ on $R_{\pi}$, where $i\in \{1,\ldots,n\}$. We also have at our disposal the operations $(i,j)\mapsto (i \star j)_{s,r}$ on $\{1,\ldots,n\}$; cf. Equation (\ref{istarjsr}).

\

Take, in Definition \ref{ohnoy}, the monoid $\mathfrak G'$ to be  the submonoid $\mathfrak D=\cup_{s\in c\mathbb N} \mathfrak D^{(s)}$ of $\mathfrak G$ and recall the monoid $\mathfrak E'$ in loc. cit., which was defined as the inverse image of $\mathfrak G'$ in $\mathfrak E$.
Note that the monoid 
$\mathfrak E'$ is naturally graded and 
$\mathfrak E$ in Equation (\ref{Esbt}) is a graded submonoid of $\mathfrak E'$.
Every  graded $\pi$-connection may be viewed as a map of sets
$$\nabla:\mathfrak D\rightarrow \mathfrak E\subset \mathfrak E',\ \ \ X\mapsto \nabla_X,\ \ \nabla_{\phi^{(s)}_i}=
(\phi^{(s)}_i)^G$$  
which is a section of  $\mathfrak E'\rightarrow \mathfrak D$ and 
is graded in the sense that it sends $\mathfrak D^{(s)}$ into $\mathfrak E^{(s)}$ for all $s\in c\mathbb N$.

  \begin{definition}\label{eei}
  Given a  graded $\pi$-connection $\nabla:\mathfrak D\rightarrow \mathfrak E'$ consider the maps
  $$
  \begin{array}{ll}
  \mathcal S^*:\mathfrak D\times \mathfrak D\rightarrow
 \mathfrak E', & \mathcal S^*(X,Y):=(\nabla_X\circ \nabla_Y)
 \cdot (\nabla_{X\circ Y})^{-1},\\
 \ & \ \\
 \mathcal R^*:\mathfrak D\times \mathfrak D\rightarrow
 \mathfrak R, &  \mathcal R^*(X,Y)  :=  \widetilde{\mathcal S^*(X,Y)}\cdot (\widetilde{\mathcal S^*(Y,X)})^{-1}.\end{array}$$
 The map $\mathcal R^*$ is called the {\bf multiplicative curvature} of $\nabla$.
\end{definition}

\begin{remark} Definition \ref{eei} is  a specialization of Definition \ref{ohnoy}. So all the remarks following that definition apply to our context here. Additionally, if one assumes $\mathfrak D^{(c)}$ abelian and sets $X,Y\in \mathfrak D$,  then we have that
$\mathcal R^*(X,Y)=1$ if and only if $\mathcal R(X,Y)=0$.
\end{remark}

\
  
Explicitly, for all $i,j\in \{1,\ldots,n\}$ and all $s,r\in c\mathbb N$, writing
\begin{equation}
  \label{carru}
  \widetilde{\phi^{(s)}_i},\ \widetilde{\phi^{(s)}_i\circ \phi^{(r)}_j},\ldots\in \mathfrak R\end{equation}
in place of
  $$\widetilde{(\phi^{(s)}_i)^G},\ \widetilde{(\phi^{(s)}_i)^G\circ (\phi^{(r)}_j)^G},\ldots\in \mathfrak R$$
we have
$$\mathcal R^*(\phi^{(s)}_i,\phi^{(r)}_j):=
(\widetilde{\phi^{(s)}_i \circ \phi^{(r)}_j})
\cdot(\widetilde{\phi^{(s+r)}_{(i\star j)_{s,r}}})^{-1}
\cdot 
(\widetilde{\phi^{(s+r)}_{(j\star i)_{r,s}}})
\cdot (\widetilde{\phi^{(r)}_j \circ \phi^{(s)}_i})^{-1}.
$$

\

\noindent In particular, if $\mathfrak D^{(c)}$ is abelian we have
$$\mathcal R^*(\phi^{(c)}_i,\phi^{(c)}_j)=
(\widetilde{\phi^{(c)}_i \circ \phi^{(c)}_j})\cdot (\widetilde{\phi^{(c)}_j \circ \phi^{(c)}_i})^{-1}.
$$

\

We next ``specialize" the multiplicative curvature to points in $G(R_{\pi})$.
Indeed, for every point $a\in G(R_{\pi})$ and every $g\in \mathfrak R$, we have a well defined 
point $g(a)\in G(R_{\pi})$: if $g(x)=F\in \textup{GL}_n(\mathcal A)$
then $g(a):=F(a)$.
 Clearly, for every point $a\in G(R_{\pi})$
the point $\mathcal R^*(\phi^{(s)}_i,\phi^{(r)}_j)(a)$
 reduces mod $\pi$ to the identity $1\in G(k)$;
so $\mathcal R^*(\phi^{(s)}_i,\phi^{(r)}_j)(a)$ belongs to the group $G^1(R_{\pi})$ which we can identify with the group $\mathfrak g(R_{\pi})$ via the isomorphism in \cite[Eq. 4.7]{BM22}.
Under this identification  we can write
$$\mathcal R^*(\phi^{(s)}_i,\phi^{(r)}_j)(a)\in \mathfrak g(R_{\pi}).
$$
 Writing
$$\phi_i^{(s)}(x)=x^{(p^s)}\cdot \Lambda_i^{(s)},\ \ 
\phi_j^{(r)}(x)=x^{(p^r)}\cdot \Lambda_j^{(s)}$$
and using  formula (\ref{binee}) we have that
$$
\phi_i^{(s)}\phi^{(r)}_j (x)= (x^{(p^s)}\cdot \Lambda_i^{(s)})^{(p^r)}\cdot (\Lambda_j^{(r)})^{(\phi_i^{(s)})}(x^{(p^s)}\cdot \Lambda_i^{(s)}).$$

\

\noindent In particular, we have
$$\widetilde{\phi_i^{(s)}}(1)=\Lambda_i^{(s)}(1)=1+\pi \Gamma^{(s)t}_i(1)$$
and
$$\begin{array}{rcl}
\widetilde{\phi_i^{(s)}\phi^{(r)}_j}(1) & = & (\Lambda_i^{(s)}(1))^{(p^r)}\cdot (\Lambda_j^{(r)})^{(\phi_i^{(s)})} (\Lambda_i^{(s)}(1))\\
\ & \ & \ \\
\ & \equiv & (\Lambda_j^{(r)})^{(\phi_i^{(s)})}(\Lambda_i^{(s)}(1))\ \ \textup{mod}\ \ \pi^2\\
\ & \ & \ \\
\ & \equiv & 1+\phi_i^{(s)}\pi \cdot (\Gamma^{(r)t}_j)^{(\phi_i^{(s)})}(\Lambda_i^{(s)}(1))\ \ \textup{mod}\ \ \pi^2\\
\ & \ & \ \\
\ & \equiv & 1+\pi \cdot \delta_i^{(s)}\pi \cdot (\Gamma^{(r)t}_j)^{(\phi_i^{(s)})}(1)\ \ \textup{mod}\ \ \pi^2\\
\ & \ & \ \\
\ & \equiv & 1+\pi \cdot \delta_i^{(s)}\pi \cdot (\Gamma^{(r)t}_j)^{(\phi_i^{(s)})}\ \ \textup{mod}\ \ \mathcal M^2\\
\ & \ & \ \\
\ & \equiv & 1+\pi \cdot \delta_i^{(s)}\pi \cdot (\Gamma^{(r)t}_j(1))^{(\phi_i^{(s)})}\ \ \textup{mod}\ \ \mathcal M^2\\
\ & \ & \ \\
\ & \equiv & 1+\pi \cdot \delta_i^{(s)}\pi \cdot (\Gamma^{(r)t}_j(1))^{(p^s)}\ \ \textup{mod}\ \ \mathcal M^2.\end{array}$$
Hence
$$(\widetilde{\phi_i^{(s)}\circ \phi^{(r)}_j})(1)\equiv 1+\pi \cdot \delta_i^{(s)}\pi \cdot (\Gamma^{(r)t}_j(1))^{(p^s)}\ \ \textup{mod}\ \ \pi^2
$$
because the homomorphism $R_{\pi}/\pi^2R_{\pi}\rightarrow \widehat{A}/\mathcal M^2$ is injective. We get the following congruence in $\textup{Mat}_n(R_{\pi})$:
\begin{equation}
\label{tyu}
\begin{array}{rcl}
\mathcal R^*(\phi^{(s)}_i,\phi^{(r)}_j)(1) & \equiv &
1+\pi(\delta_i^{(s)}\pi \cdot (\Gamma^{(r)t}_j(1))^{(p^s)}-
\delta_j^{(r)}\pi \cdot (\Gamma^{(s)t}_i(1))^{(p^r)})\\
\ & \ & \ \\
\ & \  & +\pi(\Gamma_{j\star i}^{(s+r)t}(1)-\Gamma_{i\star j}^{(s+r)t}(1))\ \ \ 
\textup{mod}\ \ \ \pi^2.\end{array}
\end{equation}

\

\noindent Let $\overline{R}^{k*}_{ijl}\in \textup{Mat}_n(k)$ be the image of the $kl$-entry $(\mathcal R^*(\phi_i^{(c)},\phi_j^{(c)})(1))_{kl}$ of the matrix 
$$\mathcal R^*(\phi_i^{(c)},\phi_j^{(c)})(1)\in G^1(R_{\pi})\simeq \mathfrak g(R_{\pi})$$ via the reduction map $$\mathfrak g(R_{\pi})\rightarrow \mathfrak g(k)=\textup{Mat}_n(k).$$
 Also recall the elements $\overline{R}^k_{ijl}$ in Equation (\ref{misssing}).
 
 \medskip
 
 The following theorem  shows that the multiplicative curvature and the ``additive" curvature ``coincide mod $\pi$ at the origin." 

\begin{theorem}\label{ursuh}
For all $i,j,k,l$ we have equalities
$$\overline{R}^{k*}_{ijl}=\overline{R}^k_{ijl}.$$
\end{theorem}

{\it Proof}.
The statement follows directly from Equations (\ref{tyu}) and (\ref{misssing}) using the fact that  the image
of $\Gamma^{k(s)}_{ij}(1)\in \textup{Mat}_n(R_{\pi})$  in $\textup{Mat}_n(k)$ equals $\overline{\Gamma}^{k(s)}_{ij}$.
\qed

\bigskip

It is interesting to see what happens if one considers the following  alternative  version of the multiplicative curvature. 

\begin{definition}
Define the {\bf naive multiplicative curvature} to be the map
 $$\mathcal R^{**}:\mathfrak D\times \mathfrak D\rightarrow
 \mathfrak R$$
 given by the formula
$$\mathcal R^{**}(X,Y):=
(\widetilde{\nabla_X} \circ \widetilde{\nabla_Y})\cdot 
(\widetilde{\nabla_Y} \circ \widetilde{\nabla_X})^{-1}.
$$
\end{definition}

In other words,
$$\mathcal R^{**}(\phi^{(s)}_i,\phi^{(r)}_j):=
(\widetilde{\phi^{(s)}_i} \circ \widetilde{\phi^{(r)}_j})\cdot 
(\widetilde{\phi^{(r)}_j} \circ \widetilde{\phi^{(s)}_i})^{-1}.
$$

Let 
$$\overline{R}^{k**}_{ijl}\in \textup{Mat}_n(k)$$ be the image of the $kl$-entry $(\mathcal R^{**}(\phi_i^{(c)},\phi_j^{(c)})(1))_{kl}$ of the matrix 
$$\mathcal R^{**}(\phi_i^{(c)},\phi_j^{(c)})(1)\in G^1(R_{\pi})\simeq \mathfrak g(R_{\pi})$$ via the reduction map $$\mathfrak g(R_{\pi})\rightarrow \mathfrak g(k)=\textup{Mat}_n(k).$$

\begin{theorem}
For all $i,j,k,l$ we have equalities
$$\overline{R}^{k**}_{ijl}=\overline{\Gamma}^{k(c)}_{il}-\overline{\Gamma}^{k(c)}_{jl}.$$
In particular if $\mathfrak D^{(c)}$ is abelian and $\nabla$ is symmetric with respect to the (additive or multiplicative) 
canonical torsion symbol  then
$$\overline{R}^{k**}_{ijl}=0.$$
\end{theorem}

{\it Proof}.
The first statement follows by direct computation similar to the one leading to Theorem \ref{ursuh}. The second statement follows from formula (\ref{uf2}).
\qed

 \section{Gauge group}\label{torsors} Recall as before that we view $\textup{GL}_N$ as a trivial torsor under itself and indeed there are no non-trivial torsors for the algebraic group $\textup{GL}_N$ over a field by Hilbert's Theorem 90. All the more, the  ``gauge group" underlying our theory to be not $\textup{GL}_N(R_{\pi})$. In this section, we will replace the ``gauge group" by a smaller group as follows.

\subsection{Ad-invariant metrics}

Taking a clue from \cite{Bu17}  we will take the  ``gauge group" underlying our theory to be a group $\mathfrak W$ which will be reviewed in what follows. This group can be viewed as the ``general linear group over the algebraic closure of the field with one element."

 \begin{notation}\label{defofW}
 Let $W\subset \textup{GL}_N(R_{\pi})$ be the subgroup consisting of all permutation matrices
 (the ``classical Weyl group" of $\textup{GL}_N$) and let $T^{\delta}\subset \textup{GL}_N(R)$ be the subgroup of diagonal matrices with diagonal entries roots of unity in $R$. Let $\mathfrak W=WT^{\delta}=T^{\delta}W$ the subgroup generated by 
  $T^{\delta}$ and $W$.
 (In \cite{Bu17} this group was denoted by $N^{\delta}$ but, in order to avoid confusion with our other use of the letter $N$ we are adopting here a different notation.)

 \

 \noindent Note that for all $a\in \mathfrak W$ and $b\in \textup{Mat}_N(A)$ with $A$ an $R$-algebra and for all $s\in \mathbb N$ we have
 \begin{equation}
 \label{neumflat}
 (ab)^{(p^s)}=a^{(p^s)}b^{(p^s)},\ \ \ (ba)^{(p^s)}=b^{(p^s)}a^{(p^s)}.
 \end{equation}
 Moreover if $a\in W$ then $a^t=a^{-1}$.
 Fix a tuple $\Phi^{(s)}=(\phi_1^{(s)},\ldots,\phi^{(s)}_n)$ of higher $\pi$-Frobenius
 lifts  of degree $s$ on $R_{\pi}$. Then for all $a\in \mathfrak W$ we have
 \begin{equation}
 \label{rayy200}
 \phi_i^{(s)}(a)=a^{(p^s)}.\end{equation}

 \

 For every $w\in \mathfrak W$ consider the left multiplication $R_{\pi}$-algebra automorphism 
 \begin{equation}
 \label{twG}
 ^Gw:\mathcal A\rightarrow \mathcal A,\ \ \ ^Gw(x)=wx.\end{equation}
 The above formula defines a right action 
 $$\mathfrak W\rightarrow \textup{Aut}_{R_{\pi}-\textup{alg}}(\mathcal A),\ \ 
 \ w \mapsto \ ^Gw.$$\end{notation}
 
 We have the following compatibility between these maps and the Levi-Civita (respectively Chern) connections:
 
 \begin{proposition}\label{gaugecomp1}
 Consider a metric $q^{(s)}\in \textup{GL}_n(R_{\pi})^{\textup{sym}}$, an element $w\in \mathfrak W$, and the metric
 $\tilde{q}^{(s)}:=w^tq^{(s)}w$. 
 Let $L^{(s)}$ be a torsion symbol and let $\Delta^{\textup{(s)LC}}$ be the Levi-Civita connection attached to $q^{(s)}$ and $L^{(s)}$.  Similarly,  $\tilde{\Delta}^{\textup{(s)LC}}$ be the Levi-Civita connection attached to $\tilde{q}^{(s)}$ and $L^{(s)}$.
 Consider 
 $(\phi_1^{(s)\textup{LC}},\ldots,\phi_n^{(s)\textup{LC}})$ and  $(\tilde{\phi}_1^{(s)\textup{LC}},\ldots,\tilde{\phi}_n^{(s)\textup{LC}})$ the  tuples of $\pi$-Frobenius lifts of degree $s$ attached to  $\Delta^{\textup{(s)LC}}$ and $\tilde{\Delta}^{\textup{(s)LC}}$, respectively. 

 \

 For all $i\in \{1,\ldots,n\}$ we have the following equality of maps $\mathcal A\rightarrow \mathcal A$:
 \begin{equation}
 \label{pptt}
 \tilde{\phi}_i^{(s)\textup{LC}}\circ \ ^Gw=\ ^Gw\circ \phi^{(s)\textup{LC}}_i.\end{equation}
 \end{proposition}
 
 \
 
 \begin{proposition}\label{gaugecomp2}
 Consider a metric $q^{(s)}\in \textup{GL}_N(R_{\pi})^{\textup{sym}}$, an element $w\in \mathfrak W$, and the metric
 $\tilde{q}^{(s)}:=w^tq^{(s)}w$. 
  Let  $\Delta^{\textup{(s)Ch}}$ be the Chern connection attached to $q^{(s)}$ and let $\tilde{\Delta}^{\textup{(s)Ch}}$ be the Chern connection attached to $\tilde{q}^{(s)}$.
 Let 
 $(\phi_1^{(s)\textup{Ch}},\ldots,\phi_n^{(s)\textup{Ch}})$ and  $(\tilde{\phi}_1^{(s)\textup{Ch}},\ldots,\tilde{\phi}_n^{(s)\textup{Ch}})$ be the  tuples of $\pi$-Frobenius lifts of degree $s$ attached to  $\Delta^{\textup{(s)Ch}}$ and $\tilde{\Delta}^{\textup{(s)Ch}}$, respectively. 

 \

 For all $i\in \{1,\ldots,n\}$ we have the following equality of  maps  $\mathcal A\rightarrow \mathcal A$:
 $$\tilde{\phi}_i^{(s)\textup{Ch}}\circ \ ^Gw=\ ^Gw\circ \phi^{(s)\textup{Ch}}_i.$$
 \end{proposition}

 \noindent We proceed directly with the proofs of these propositions. 

 \ 
 
 \noindent {\it Proof of Propositions \ref{gaugecomp1} and \ref{gaugecomp2}}.
 We prove Proposition \ref{gaugecomp1}; Proposition \ref{gaugecomp2} is proved similarly.
 Consider the matrices $A_i^{(s)},B^{(s)},\Lambda_i^{(s)}$ in Equation (\ref{1955}) corresponding to $\Delta^{\textup{(s)LC}}$ and similarly the matrices $\tilde{A}_i^{(s)},\tilde{B}^{(s)},\tilde{\Lambda}_i^{(s)}$  corresponding to $\tilde{\Delta}^{\textup{(s)LC}}$. 
 Set $$\Lambda_{i,w}^{(s)}:=\ ^Gw(\Lambda_i^{(s)}),\ \ 
 \Lambda_{w}^{(s)}:=\ ^Gw(\Lambda^{(s)}).$$
 Using Equations (\ref{neumflat}) and (\ref{rayy200}) and formulae (\ref{1955}) 
 one trivially gets
 $$^Gw(A_i^{(s)})=\tilde{A}_i^{(s)},\ \ ^Gw(B^{(s)})=\tilde{B}^{(s)}.$$
Applying $^Gw$ to  Equation (\ref{LALB}) we  get
 $$\Lambda_{i,w}^{(s)t}\tilde{A}_i^{(s)}\Lambda_{i,w}^{(s)}=\tilde{B}^{(s)}.$$
 Similarly, applying $^Gw$ to  Equation (\ref{uf22}) we  get 
 $$(\Lambda_{i,w}^{(s)}-1)_{kj}-(\Lambda_{j,w}^{(s)}-1)_{ki}=\pi\cdot L_{ij}^{k(s)}(\Lambda_w^{(s)}).
$$
By the uniqueness of the Levi-Civita connection attached to $\tilde{q}^{(s)}$ and $L^{(s)}$ we get
$$\tilde{\Lambda}_i^{(s)}=\Lambda_{i,w}^{(s)}.$$
Using the last equation in (\ref{1955}) one immediately gets the equality (\ref{pptt}).
 \qed

\bigskip

As an application of the above discussion we explore, 
in what follows,  an arithmetic analogue of bi-invariant metrics on Lie groups.
Assume we are in the setting of Subsection \ref{calcc}. In particular  we assume  $N=n$ and we assume we have an involutive set $\mathfrak D^{(c)}$ of cardinality $n$, with corresponding group $\mathfrak S$, and  canonical Lie algebra 
$\mathfrak d$. Also for all $s\in c\mathbb N$ we have a labeling symbol $\omega^{(s)}$, a torsion symbol $L^{(s)}$ of the second kind, a metric $q^{(s)}$,  and the attached arithmetic Levi-Civita connection
$\Delta^{(s)\textup{LC}}$ with $\pi$-Frobenius lifts $\phi_i^{\textup{(s)LC}}$. Also we have 
 the associated curvature $\mathcal R^{\textup{LC}}$.

 \

 \noindent For each $s$ we consider the bilinear symmetric map 
 $$\langle\ \ ,\ \ \rangle:\mathfrak d\times \mathfrak d\rightarrow R_{\pi}$$
 defined by 
 $$\langle \phi_i^{(s)},\phi_j^{(s)}\rangle:=q_{ij}^{(s)}\ \ \textup{and}\ \ \langle \phi_i^{(s)},\phi_j^{(r)}\rangle:=0\ \ \textup{for}\ \ r\neq s.$$
 Note that for all  $\sigma\in \mathfrak S$ and all $X\in \mathfrak d^{(s)}$ (respectively $X\in \mathfrak D^{(s)}$) we have
 $\sigma X \sigma^{-1}\in \mathfrak d^{(s)}$ (respectively $\sigma X \sigma^{-1}\in \mathfrak D^{(s)}$). Write
 $$\sigma \phi_i^{(s)}\sigma^{-1}=\phi_{\sigma^{(s)}(i)}.$$
 We obtain a group homomorphism
 $$\textup{Ad}^{(s)}:\mathfrak S\rightarrow \Sigma_n,\ \
  \textup{Ad}^{(s)}(\sigma)(i):=\sigma^{(s)}(i),\ \ \textup{for}\ \ 
 \sigma\in \mathfrak S\ \ \textup{and}\ \  i\in \{1,\ldots,n\},$$
 where $\Sigma_n$ is the symmetric group.  Recall that for any permutation $\epsilon\in \Sigma_n$ we let $P_{\epsilon}$ be the corresponding permutation matrix defined by
  $(P_{\epsilon})_{ij}=\delta_{\epsilon(i)j}$; hence 
$(P_{\epsilon}q^{(s)}P_{\epsilon}^t)_{ij}=q^{(s)}_{\epsilon(i)\epsilon(j)}$.

 \begin{definition}\label{adinv}
 We say $q^{(s)}$ is {\bf \textup{Ad}-invariant} if for all $\sigma\in \mathfrak S$ and all $X,Y\in \mathfrak d^{(s)}$ we have
 $$\langle \sigma X \sigma^{-1},\sigma Y \sigma^{-1}\rangle=\langle X,Y\rangle.$$
 Equivalently, for all $i,j$
 $$q^{(s)}_{\sigma^{(s)}(i)\sigma^{(s)}(j)}=q^{(s)}_{ij}.
 $$
 Equivalently,
 $$P_{\textup{Ad}^{(s)}(\sigma)}q^{(s)}P_{\textup{Ad}^{(s)}(\sigma)}^t=q^{(s)}.$$
 \end{definition}
 
 Recall from Equation \ref{twG} the definition of the $R_{\pi}$-automorphisms $^Gw$ of $\mathcal A$ for $w\in \mathfrak W$.
  Propositions \ref{gaugecomp1} and \ref{gaugecomp2}  imply the following statements (where the notation is from loc.cit):
 
 \begin{corollary}
 Assume that $q^{(s)}$ is \textup{Ad}-invariant for some $s\in c\mathbb N$.
 Let $(\phi_i^{(s)})^G$ be either $\phi_i^{\textup{(s)LC}}$ or $\phi_i^{\textup{(s)Ch}}$
 and let $\mathcal R$ be either $\mathcal R^{\textup{LC}}$ or $\mathcal R^{\textup{Ch}}$.
 Then for all $i\in \{1,\ldots,n\}$ and all $\sigma\in \mathfrak S$ we have that the ring endomorphisms 
  $$(\phi_i^{(s)})^G:\mathcal A\rightarrow \mathcal A\ \ \ \textup{and}\ \ \ 
  ^G(P_{\textup{Ad}^{(s)}(\sigma)}):\mathcal A\rightarrow \mathcal A$$ commute. In particular, if $q^{(s)}$ is \textup{Ad}-invariant for all $s\in c\mathbb N$ then
  for all $i,j\in \{1,\ldots,n\}$, all $r,s\in c\mathbb N$,  and all $\sigma\in \mathfrak S$
   the $\mathbb Z$-module endomorphisms
  $$\mathcal R(\phi^{(s)}_i,\phi_j^{(r)}):\mathcal A\rightarrow \mathcal A\ \ \ \textup{and}\ \ \ 
  ^G(P_{\textup{Ad}^{(s)}(\sigma)}):\mathcal A\rightarrow \mathcal A$$ commute. 
   \end{corollary}

\begin{example}\label{399}
Assume we are under the assumptions of  Example \ref{ogarcenusiu}.

\

1) If  $n=2$  then $\textup{Ad}^{(c)}$ is the trivial homomorphism hence every metric $q^{(c)}$ is $\textup{Ad}^{(c)}$-invariant.

\

2) If  $n=3$   
and the action of $\phi$ on $\mathfrak S$ by conjugation is trivial then, again, $\textup{Ad}^{(c)}$ is the trivial homomorphism and hence every metric $q^{(c)}$ is $\textup{Ad}^{(c)}$-invariant.

\

 3) If $n=3$ and  the action of $\phi$ on $\mathfrak S$ by conjugation
   is non-trivial then one can easily see that
   $$P_{\textup{Ad}^{(c)}}(\sigma_1)=1,\ \ \ 
   P_{\textup{Ad}^{(c)}}(\sigma_2)=\left(
   \begin{array}{ccc} 0 & 1 & 0 \\ 0 & 0 & 1\\ 1 & 0 & 0\end{array}\right),\ \ \ 
   P_{\textup{Ad}^{(c)}}(\sigma_3)=\left(
   \begin{array}{ccc} 0 & 0 & 1 \\ 1 & 0 & 0\\ 0 & 1 & 0\end{array}\right).
   $$
   Hence a metric $q^{(c)}$ is $\textup{Ad}^{(c)}$-invariant if and only it has the form
   $$q^{(c)}=\left(
   \begin{array}{ccc} a & b & b \\ b & a & b\\ b & b & a\end{array}\right),\ \ \ a,b \in R_{\pi},\ \ \ a^3 +2b^3-3ab^2\in R^{\times}.$$
So we have a ``$2$-parameter family" of $\textup{Ad}^{(c)}$-invariant metrics, rather than a ``$6$-parameter family" as in situation 2).
\end{example}

\subsection{Torsors}
 
 Throughout this section  we freely use the terminology of Subsections \ref{cannn} and \ref{gencur}.
 The symbol $H^1_{\textup{gr}}$ is used to denote group cohomology. 
 We fix in what follows  a tuple $\Phi^{(s)}=(\phi_1^{(s)},\ldots,\phi^{(s)}_n)$ of higher $\pi$-Frobenius
 lifts  of degree $s$ on $R_{\pi}$ and we let the integer $N$ be
 not necessarily equal to $n$.
 
 \

Since $H^1_{\textup{gr}}(\mathfrak G(L/K), \textup{GL}_N(L))=1$
 for every finite Galois extension $L/K$ (Hilbert's Theorem $90$) there seems to be no ``naive" analogue, in our context, for non-trivial torsors under $\textup{GL}_N$. 
  However, as we shall see in this Subsection, one can develop an analogue of torsors (and their relation with our arithmetic Riemannian concepts) for the ``gauge group" 
 $\mathfrak W$ (cf. Notation \ref{defofW})  that has non-trivial Galois cohomology. We recall that $\mathfrak W=WT^{\delta}$ where $W$ is the group of permutation matrices and $T^{\delta}$ is the group of diagonal matrices with entries roots of unity in $R$. 
 We also recall the formulae
 (\ref{neumflat}) and (\ref{rayy200}).

 \

 \begin{notation}
 Let  $\mathfrak T\subset \mathfrak G(K_{\pi}/K)$ be a subgroup of order $m$ normalized by all $\phi^{(s)}_i$'s. (E.g.,
  one can take the situation in Example \ref{boooo} with $\mathfrak T$ any invariant subgroup of $\mathfrak S$; but our $\mathfrak T$ need not be related to $\mathfrak S$ in general.)
 Write $\mathfrak T=\{\tau_1,\ldots,\tau_m\}$. We view $\mathfrak W$  as a trivial $\mathfrak T$-module. 
  We denote by 
 $$Z^1_{\textup{gr}}(\mathfrak T,\mathfrak W)=\textup{Hom}_{\textup{gr}}(\mathfrak T,\mathfrak W),$$
 the set of (group) $1$-cocycles and consider the cohomology pointed set
 $$
 H^1_{\textup{gr}}(\mathfrak T,\mathfrak W)=\textup{Hom}_{\textup{gr}}(\mathfrak T,\mathfrak W)/\mathfrak W,$$
 where $\mathfrak W$ acts on $\textup{Hom}_{\textup{gr}}(\mathfrak T,\mathfrak W)$ by conjugation; so two cocycles are  cohomologous if and only if they are $\mathfrak W$-conjugate.  The homomorphism $\mathfrak W\rightarrow \mathfrak W/T^{\delta}\simeq W$ has a section in the category of groups with trivial $\mathfrak T$-action (a section is given by the inclusion) so the map
 $$H^1_{\textup{gr}}(\mathfrak T,\mathfrak W)\rightarrow H^1_{\textup{gr}}(\mathfrak T,W)=\textup{Hom}_{\textup{gr}}(\mathfrak T,W)/W$$
 has a section in the category of sets, in particular it is surjective. Also since $T^{\delta}$ is abelian one gets that the map
 $$\textup{Hom}_{\textup{gr}}(\mathfrak T,T^{\delta})/W\rightarrow H^1_{\textup{gr}}(\mathfrak T,\mathfrak W) $$
 is injective.

 \

 The pointed set $ H^1_{\textup{gr}}(\mathfrak T,\mathfrak W)$ is generally non-trivial which stands in stark contrast with Hilbert's Theorem 90 and opens up the possibility of considering non-trivial torsors in our context. 
  For all
 $i\in \{1,\ldots,n\}$ and $j\in \{1,\ldots,m\}$ we write
 \begin{equation}
 \label{rayy2}
 \tau_{j,i}:=\phi_i^{(s)}\tau_j(\phi_i^{(s)})^{-1}.\end{equation}
 Denote by $[\Phi^{(s)},\mathfrak T]$ the subgroup of $\textup{Aut}_{\mathbb Z_p-\textup{alg}}(R_{\pi})$ generated by the set of commutators
 $$\{\phi_i^{(s)}\tau_j(\phi_i^{(s)})^{-1}\tau_j^{-1}\ |\  i\in \{1,\ldots,n\},\ \ j\in \{1,\ldots,m\}\}.$$
 Note that $[\Phi^{(s)},\mathfrak T]$ is a subgroup of $\mathfrak T$.\end{notation}

 \begin{definition}
 A cocycle $u\in Z^1_{\textup{gr}}(\mathfrak T,\mathfrak W)$, $\tau_i\mapsto u_{\tau_i}$, is said to be {\bf $\Phi^{(s)}$-invariant} if  for all $i\in \{1,\ldots,n\}$ and $j\in \{1,\ldots,m\}$ we have
 \begin{equation}
 \label{rayy3}
 u_{\tau_{j,i}}=u_{\tau_j};\end{equation}
 equivalently, if
 \begin{equation}
 \label{rayy381}
 [\Phi^{(s)},\mathfrak T]\subset \textup{Ker}(u).
 \end{equation}
 \end{definition}
 
 \

 \begin{remark}\label{rayy129} If $u$ and $v$ are two cohomologous cocycles in $Z^1_{\textup{gr}}(\mathfrak T,\mathfrak W)$ and $u$ 
 is $\Phi^{(s)}$-invariant then $v$ is also  $\Phi^{(s)}$-invariant.  \end{remark}
 
  \begin{definition}\label{remtorsor}
 For every cocycle $u\in Z^1_{\textup{gr}}(\mathfrak T,\mathfrak W)$ consider the ring of invariants 
 $$R_{\pi}^u:=\{a\in R_{\pi}\ |\ \tau_j(a)=a,\ j\in \{1,\ldots,m\}\}.$$
 (Of course $R_{\pi}^u$ only depends on $\mathfrak T$ and not on the choice of the cocycle $u$.)
 If  $u$ is given by $\tau_i\mapsto u_{\tau_i}$  define a  left $\mathfrak T$-action on $\mathcal A=R_{\pi}[x,\det(x)^{-1}]^{\widehat{\ }}$, 
  $$\mathfrak T\rightarrow \textup{Aut}_{R_{\pi}^u-\textup{alg}}(\mathcal A),\ \ \ 
  \tau_j\mapsto \tau_j^{G,u},$$
 by letting the restriction of $\tau_j^{G,u}$  to $R_{\pi}$  be $\tau_j$ and by letting
 \begin{equation}
 \label{rayy4}
 \tau_j^{G,u}(x):=u_{\tau_j}^{-1}x.\end{equation}
 If $u$ is understood from context we write $\tau_j^G=\tau_j^{G,u}$.
 We denote by $\mathcal A^u$ the ring of invariants 
 \begin{equation}
 \label{eqtor}
 \mathcal A^u=\{f\in \mathcal A\ |\ \tau_j^G(f)=f,\ j\in \{1,\ldots,m\}\}.\end{equation}
 It is a $p$-adically complete  algebra over $R_{\pi}^u$. 
 Furthermore we consider the left action
 $$G(R_{\pi}) \rightarrow \textup{Aut}_{R_{\pi}-\textup{alg}}(\mathcal A),\ \ \ 
  g\mapsto g^G,$$
 defined by 
 \begin{equation}
 \label{rayy4.5}
 g^G(x):=xg.\end{equation}\end{definition}
 
 \begin{remark}
 Let  $u,v\in Z^1_{\textup{gr}}(\mathfrak T,\mathfrak W)$ be cohomologous with $v=a^{-1}ua$ for some  $a\in \mathfrak W$. Then
 the $R_{\pi}$-algebra automorphism $\alpha$ of $\mathcal A$ defined by $x\mapsto \alpha(x):=ax$ trivially satisfies the following conditions:
 
 \
 
 1)  $\alpha$ is $\textup{GL}_N(R_{\pi})$-equivariant;
 
 \
 
 2) $\alpha$ commutes with the $\mathfrak T$-actions on $\mathcal A$ 
 defined by $u$ and $v$ in the sense that
 $$\alpha \circ \tau_j^{G,u}=\tau_j^{G,v}\circ \alpha\ \ \textup{for all} \ \ j\in \{1,\ldots,m\};$$
 
 \
 
 3)  $\alpha$ commutes with the trivial  $\pi$-connection  of degree $s$ in the sense that
 $$\alpha \circ (\phi_{i,0}^{(s)})^G=(\phi_{i,0}^{(s)})^G\circ \alpha\ \ \textup{for all} \ \ i\in \{1,\ldots,n\}.
 $$
 \end{remark}
 
Conversely we have the following proposition.

 \begin{proposition}
 Assume $u,v\in Z^1_{\textup{gr}}(\mathfrak T,\mathfrak W)$ and $\alpha$ is an $R_{\pi}$-algebra automorphism  of $\mathcal A$ satisfying the conditions 1), 2), 3) above.
  Then there exists $a\in \mathfrak W$ such that $\alpha(x)=ax$ and $v=a^{-1}ua$; in particular $u$ and $v$ are cohomologous.\end{proposition}
 
 {\it Proof}.  Indeed, condition 1) easily implies that there exists $a\in \textup{GL}_N(R_{\pi})$ such that $\alpha(x)=ax$. Condition 2) easily implies that
 $v_{\tau_j}=a^{-1}u_{\tau_j} a^{\tau_j}$ for all $j$. Condition 3) easily implies 
 that $\phi_i(a)x^{(p^s)}=(ax)^{(p^s)}$ for all $i$. Setting $x=1$ we get 
 $\phi_i^{(s)}(a)=a^{(p^s)}$ for all $i$. Hence $a^{(p^s)}x^{(p^s)}=(ax)^{(p^s)}$ which easily implies that $a$ is in the group generated by $W$ and the group of diagonal matrices with entries in $R_{\pi}$. On the other hand by  \cite[Lem. 2.17]{BM22}  applied to the entries of $a$, we get that every entry of $a$ is either $0$ or a root of unity in $R$. Hence $a\in \mathfrak W$ and   $a^{\tau_j}=a$ which ends the proof.\qed

 \begin{remark}
 The induced   $G(R_{\pi}^u)$-action on $\mathcal A$
  commutes with the $\mathfrak T$-action on $\mathcal A$
  in the sense that
 $$g^G\circ \tau_j^G=\tau_j^G\circ g^G$$
 for all $g\in G(R_{\pi}^u)$ and all $j\in \{1,\ldots,m\}$.
 Consequently we have an induced action
 $$G(R_{\pi}^u)\rightarrow \textup{Aut}_{R_{\pi}^u-\textup{alg}}(\mathcal A^u).
 $$
 We intuitively view the $R_{\pi}^u$-algebra
 $\mathcal A^u$ with its $G(R_{\pi}^u)$-action 
 as ``the torsor" attached to $u$.

 \

 If $u,v\in Z^1_{\textup{gr}}(\mathfrak T,\mathfrak W)$ are cohomologous with $v=a^{-1}ua$ for some  $a\in \mathfrak W$ then we have a natural $R_{\pi}^u$-algebra isomorphism 
 \begin{equation}
 \label{rayy440}
 \mathcal A^u \rightarrow \mathcal A^v\end{equation}
 induced by the $R_{\pi}$-algebra isomorphism $\mathcal A\rightarrow \mathcal A$, $x\mapsto ax$. The isomorphism (\ref{rayy440}) is equivariant with respect to the corresponding 
 $G(R_{\pi}^u)$-actions on $\mathcal A^u$ and $\mathcal A^v$, respectively.
 (Recall $R_{\pi}^u=R_{\pi}^v$.)
 \end{remark}
 
\begin{notation}
For $F=F(x)\in \mathcal A$ we write $F^{(\tau_i)}$ for the element of $\mathcal A$ obtained from $F$ (viewed as a $p$-adic limit of rational functions in the entries of $x$) by applying $\tau_i$
 to the coefficients of $F$ (and leaving $x$ unchanged). Then  we have the formula
 \begin{equation}
 \label{rayy5}
 \tau_i^G(F(x))=F^{(\tau_i)}(u_{\tau_i}^{-1}x).\end{equation}\end{notation}
 
 \

 Assume now we have a $\pi$-connection  $\Delta^{(s)G}:=((\delta_1^{(s)})^G,\ldots,(\delta_n^{(s)})^G)$ of degree $s$ on $G$ with attached higher $\pi$-Frobenius lifts 
 $\Phi^{(s)G}=(\phi_1^{(s)G},\ldots,\phi_n^{(s)G})$
  of degree $s$ extending $\Phi^{(s)}$.   
  
 \begin{definition}
The 
 cocycle $u\in Z^1_{\textup{gr}}(\mathfrak T,\mathfrak W)$ is {\bf compatible with the
  $\pi$-connection} $\Delta^{(s)G}$ if for all $i\in \{1,\ldots,n\}$ and $j\in \{1,\ldots,m\}$  the following equality holds in $\textup{End}_{\mathbb Z_p-\textup{alg}}(\mathcal A)$:
 \begin{equation}
 \label{umflatura}
 (\phi_i^{(s)})^G\circ \tau_j^G=(\tau_{j,i})^G\circ (\phi_i^{(s)})^G.
 \end{equation}
 \end{definition}
 
 \begin{remark}
 Assume the 
 cocycle $u\in Z^1_{\textup{gr}}(\mathfrak T,\mathfrak W)$ is  compatible with the
  $\pi$-connection $\Delta^{(s)G}$. 
 Then for all $i\in \{1,\ldots,n\}$ we have 
 $$(\phi_i^{(s)})^G(\mathcal A^u)\subset \mathcal A^u$$
 and the  homomorphism
 $$\phi_i^{(s)u}:\mathcal A^u\rightarrow \mathcal A^u$$
 induced by $(\phi_i^{(s)})^G$ 
 is trivially seen to be a higher $\pi$-Frobenius lift of degree $s$.
 Hence $\mathcal A^u$ has a natural structure of  partial $\delta$-ring of degree $s$ compatible with the corresponding induced structure of partial $\delta$-ring of degree $s$ on $R_{\pi}^u$.
 \end{remark}
 
 In what follows we write, as usual, cf.  Equation (\ref{1955}):
 \begin{equation}
 \label{rayy6}
 (\phi_i^{(s)})^G(x)=x^{(p^s)}\Lambda_i^{(s)}(x).\end{equation}

 \begin{lemma}\label{aeroplann}
 Assume $u\in Z^1_{\textup{gr}}(\mathfrak T,\mathfrak W)$ is  $\Phi^{(s)}$-invariant and assume
 \begin{equation}
 \label{bubba}
 \tau_j^G \Lambda_i^{(s)}=\Lambda_i^{(s)}\end{equation}
  for all $i\in \{1,\ldots,n\}$ and $j\in \{1,\ldots,m\}$. Then the cocycle $u$ is compatible with the $\pi$-connection  $\Delta^{(s)G}=((\delta_1^{(s)})^G,\ldots,(\delta_n^{(s)})^G)$. \end{lemma}
 
 {\it Proof}.
 Clearly the equality (\ref{umflatura}) holds on $R_{\pi}$. It remains to check equality (\ref{umflatura}) on the entries of $x$. 
 Using Equation (\ref{neumflat}) we have the following computation:
 $$\begin{array}{rcll}
 \tau_{j,i}^G((\phi_i^{(s)})^G(x)) & = & \tau_{j,i}^G(x^{(p^s)}\Lambda_i^{(s)}) & \textup{cf}. \ (\ref{rayy6})\\
 \ & \ & \ &\  \\
 \ & = & (u_{\tau_{j,i}}^{-1}x)^{(p^s)} \tau_{j,i}^G(\Lambda_i^{(s)}) & \textup{cf}. \ (\ref{rayy4})\\
  \ & \ &\ & \ \\
 \ & = & (u_{\tau_{j,i}}^{(p^s)})^{-1} x^{(p^s)} \Lambda_i^{(s)} & \textup{cf}. \ (\ref{neumflat})\\
  \ & \ & \ & \ \\
  \ & = & (u_{\tau_j}^{(p^s)})^{-1} x^{(p^s)} \Lambda_i^{(s)} & \textup{cf}. \ (\ref{rayy3})\\
  \ & \ & \ & \ \\
 \ & = & (\phi_i^{(s)}(u_{\tau_j}))^{-1} x^{(p^s)} \Lambda_i^{(s)} & \textup{cf}. \  (\ref{rayy200})\\
 \ & \ & \ & \ \\
 \ & = & (\phi_i^{(s)})^G(u_{\tau_j}^{-1}x) & \textup{cf}. \ (\ref{rayy6})\\
 \ & \ & \ & \ \\
 \ & = & (\phi_i^{(s)})^G(\tau_j^G(x)) & \textup{cf}. \ (\ref{rayy4}).\end{array}
 $$
 \qed
 
 \begin{remark}
 By Lemma \ref{aeroplann}
 if $u\in Z^1_{\textup{gr}}(\mathfrak T,\mathfrak W)$ is $\Phi^{(s)}$-invariant then $u$ is compatible with the trivial  $\pi$-connection  $((\delta_{1,0}^{(s)})^G,\ldots,(\delta_{n,0}^{(s)})^G)$ of degree $s$. So $(\phi_{i,0}^{(s)})^G$ induce higher $\pi$-Frobenius lifts
 $$\phi_{i,0}^{(s)u}:\mathcal A^u\rightarrow \mathcal A^u$$ 
 and hence a structure of partial $\delta$-ring of degree $s$ on $\mathcal A^u$ which we call the {\bf trivial} structure.

 \

 Moreover for cohomologous cocycles $u,v\in  Z^1_{\textup{gr}}(\mathfrak T,\mathfrak W)$, $v=a^{-1}ua$, $a\in \mathfrak W$, the isomorphism (\ref{rayy440}) is easily seen to be compatible, in the obvious sense, with the 
 corresponding higher $\pi$-Frobenius lifts $\phi_{i,0}^{(s)u}$ and $\phi_{i,0}^{(s)v}$ 
  on $\mathcal A^u$ and $\mathcal A^v$ respectively; in other words
  (\ref{rayy440}) is an isomorphism of partial $\delta$-rings of degree $s$
  (with respect to the trivial structures). This latter fact fails if instead of $a\in \mathfrak W$ with $v=a^{-1}ua$
  we take an arbitrary $a\in \textup{GL}_N(R_{\pi})$ satisfying the same condition $v=a^{-1}ua$; and this is one of the  motivations for replacing $\textup{GL}_N$ with $\mathfrak W$
  as a ``gauge group."
   \end{remark}
 
 \begin{definition} 
 We say that a cocycle $u\in Z^1_{\textup{gr}}(\mathfrak T,\mathfrak W)$ is {\bf compatible with a metric}
 $q^{(s)}\in \textup{GL}_N(R_{\pi})$ if for all $j\in \{1,\ldots,m\}$ we have an equality:
 \begin{equation}\label{iarr}
 \tau_j q^{(s)}=u_{\tau_j}^t q^{(s)}  u_{\tau_j}.
 \end{equation}
 \end{definition}
 
 \begin{remark}\ 
 
 \

 1)
 If $u$ takes values in $W$ then $u$ is compatible with any scalar metric $q$ (i.e., any $q$ in the center of $\textup{GL}_N(R_{\pi})$).
 
 \

 2) If the metric $\tilde{q}$ is conformally equivalent to a metric $q$, $\tilde{q}=\lambda\cdot q$,
  $\lambda\in R^{\times}$ and if the cocycle $u$ is compatible with $q$ then $u$ is compatible with $\tilde{q}$.
 \end{remark}
 
 \begin{example} We discuss cocycles in the context of   Example \ref{ogarcenusiu}.  For convenience we recall the relevant objects. For simplicity we denote the action of Galois elements by superscripts.
 
 \

 1) Assume $\pi^2=p$ and  $\mathfrak T:=\mathfrak G(K_{\pi}/K)=\{\tau_1,\tau_2\}$ where $\tau_1=\textup{id}$ and 
 $\tau_2 \pi=-\pi$. Let $\Phi^{(s)}:=\{\phi_1,\phi_2\}$ with  $\phi_1 \pi=\pi$ and  $\phi_2=\phi_1 \tau_2=\tau_2 \phi_1$.  We have $\tau_{j,i}=\tau_j$ for all $i,j\in \{1,2\}$ so every cocycle $u\in Z^1_{\textup{gr}}(\mathfrak T,W)$ is $\Phi^{(s)}$-invariant.
  Consider the cocycle
 $$u_{\tau_1}=\left(\begin{array}{cc} 1 & 0\\ 0 & 1\end{array}\right),
  \ \ \ u_{\tau_2}=\left(\begin{array}{cc} 0 & 1\\ 1 & 0\end{array}\right).$$
  Let  $q\in \textup{GL}_2(R_{\pi})^{\textup{sym}}$ be a metric. Then  $u$ is compatible with $q$ if and only if $q$ is of the form
  $$q=\left(\begin{array}{ll} a & b\\ b & a^{\tau}\end{array}\right),\ \ a\in R_{\pi},\ \ b\in R.
  $$
 
 \

  2) Assume $\pi^3=p\neq 3$ and $\mathfrak T:=\mathfrak G(K_{\pi}/K)=\{\tau_1,\tau_2,\tau_3\}$
  where $\tau_1=1$, $\tau_2=\tau$, $\tau_3=\tau^2$, $\pi^{\tau}=\zeta_3\pi$, $\pi^{\tau^2}=\zeta_3^2\pi$, $\zeta_3$ a cubic root of $1$. Let $\Phi^{(s)}=\{\phi_1,\phi_2,\phi_3\}$
  with $\phi_1\pi=\pi$, $\phi_2=\phi_1\tau_2$, $\phi_3=\phi_1 \tau_3$.
  
  If $p\equiv 1$ mod $3$ we have $\tau_{j,i}=\tau_j$ for all $i,j\in \{1,2,3\}$
  so every cocycle $u\in Z^1_{\textup{gr}}(\mathfrak T,W)$ is  $\Phi^{(s)}$-invariant. Define the cocycle
  $$u_{\tau_1}=\left(\begin{array}{ccc} 1 & 0 & 0\\ 0 & 1 & 0\\ 0 & 0 & 1\end{array}\right),\ \  
  u_{\tau_2}=\left(\begin{array}{ccc} 0 & 1 & 0\\ 0 & 0 & 1\\ 1 & 0 & 0\end{array}\right),\ \ 
u_{\tau_3}=\left(\begin{array}{ccc} 0 & 0 & 1\\ 1 & 0 & 0\\ 0 & 1 & 0\end{array}\right).$$
 Let  $q\in \textup{GL}_3(R_{\pi})^{\textup{sym}}$ be a metric. Then  $u$ is compatible with $q$ if and only if $q$ is of the form
   $$q=\left(\begin{array}{lll} a & b^{\tau} & b\\ b^{\tau} & a^{\tau} & b^{\tau^2}\\
  b & b^{\tau^2} & a^{\tau^2}\end{array}\right),\ \ a,b\in R_{\pi}.
  $$

If $p\equiv 2$ mod $3$ we have $\tau_{2,i}=\tau_3$ so the only cocycle $u\in Z^1_{\textup{gr}}(\mathfrak T,W)$ that is $\Phi^{(s)}$-invariant is the trivial cocycle $u=1$.
 \end{example}

 Going back to the general case consider, in what follows,  the matrices (cf. Equation (\ref{1955})):
 \begin{equation}
 \label{rayy7}
 A_i^{(s)}:=x^{(p^s)t}\phi_i^{(s)}(q) x^{(p^s)},\ \ B^{(s)}:=(x^tq^{(s)}x)^{(p^s)}.\end{equation}

 \begin{lemma}\label{spanissh}
 If the cocycle $u\in Z^1_{\textup{gr}}(\mathfrak T,\mathfrak W)$ is  $\Phi^{(s)}$-invariant and is compatible with a metric
 $q^{(s)}\in \textup{GL}_N(R_{\pi})^{\textup{sym}}$ then the following equalities hold:
 $$\tau_j^G A_i^{(s)}=A_i^{(s)},\ \ \ \tau_j^G B^{(s)}=B^{(s)}.$$
 \end{lemma}

{\it Proof}. First we have:
$$
\begin{array}{rcll}
\tau_j^G B^{(s)} & = & \tau_j^G((x^tq^{(s)}x)^{(p^s)}) & \textup{cf}.\ (\ref{rayy7})\\
\ & \ & \ &\ \\
\ & = & ((u_{\tau_j}^{-1}x)^t \tau_j(q^{(s)}) u_{\tau_j}^{-1}x)^{(p^s)}& \textup{cf}.\ (\ref{rayy4})\\
\ & \ & \ &\ \\
\ & = & (x^t (u^t_{\tau_j})^{-1}\tau_j(q^{(s)}) u_{\tau_j}^{-1}x)^{(p^s)}& \ \\
\ & \ & \ &\ \\
\ & = & (x^t q^{(s)}x)^{(p^s)} & \textup{cf}\ (\ref{iarr})\\
\ & \ & \ & \ \\
\ & = & B^{(s)} & \textup{cf}.\ (\ref{rayy7}).
\end{array}
$$
Secondly we have:
$$
\begin{array}{rcll}
(\tau_{j,i})^G A_i^{(s)} & = & (\tau_{j,i})^G(x^{(p^s)t}\phi_i^{(s)}(q) x^{(p^s)}) & \textup{cf}.\ (\ref{rayy7}) \\
\ & \ & \ &\ \\
\ & = & (u_{\tau_{j,i}}^{-1}x)^{(p^s)t}\tau_{j,i}(\phi_i^{(s)}(q)) (u_{\tau_{j,i}}^{-1}x)^{(p^s)}& 
\textup{cf}.\ (\ref{rayy4})\\
\ & \ & \ &\ \\
\ & = & x^{(p^s)t}(u_{\tau_{j,i}}^{(p^s)t})^{-1}\tau_{j,i}(\phi_i^{(s)}(q))  
(u_{\tau_{j,i}}^{(p^s)})^{-1}x^{(p^s)}& \textup{cf}.\ (\ref{neumflat})\\
\ & \ & \ &\ \\
\ & = & x^{(p^s)t}(u_{\tau_j}^{(p^s)t})^{-1}\tau_{j,i}(\phi_i^{(s)}(q))  
(u_{\tau_j}^{(p^s)})^{-1}x^{(p^s)}& \textup{cf}.\ (\ref{rayy3})\\
\ & \ & \ &\ \\
\ & = & x^{(p^s)t}((\phi_i^{(s)}(u_{\tau_j})^t)^{-1}(\phi_i^{(s)}(\tau_j(q))  
(\phi_i^{(s)}(u_{\tau_j}))^{-1}x^{(p^s)}& \textup{cf}.\ (\ref{rayy200})\\
\ & \ & \ &\ \\
\ & = & x^{(p^s)t} \phi_i^{(s)}(((u_{\tau_j})^t)^{-1} \tau_j(q)  
u_{\tau_j}^{-1})x^{(p^s)}& \ \\
\ & \ & \ &\ \\
\ & = & x^{(p^s)t} \phi_i^{(s)}(q^{(s)})x^{(p^s)} &\textup{cf}.\ (\ref{iarr}) \\
\ & \ & \ & \ \\
\ & = & A_i^{(s)}& \textup{cf}.\ (\ref{rayy7}) .
\end{array}
$$
Since, for fixed $i$,  the map $\tau_j\mapsto \tau_{j,i}$
is a bijection from $\mathfrak T$ to itself it follows that $(\tau_j)^G A_i^{(s)}=A_i^{(s)}$ for all $i,j$.
\qed

\begin{theorem}
Assume the cocycle $u\in Z^1_{\textup{gr}}(\mathfrak T,\mathfrak W)$ is  $\Phi^{(s)}$-invariant and is compatible a metric
 $q^{(s)}\in \textup{GL}_N(R_{\pi})^{\textup{sym}}$. Then the following hold:

 \
 
 1) The cocycle $u$ is compatible with the arithmetic Chern connection  attached to the metric $q^{(s)}$. 
 
 \

 2) Let $n=N$ and let $L^{(s)}$ be a torsion symbol 
 in $\textup{Mat}_n(\widehat{R[y]})^n$. Then 
 the cocyle $u$ is compatible with the 
 arithmetic Levi-Civita connection attached to the metric $q^{(s)}$ and to $\frac{p}{\pi}L^{(s)}$.
\end{theorem}

{\it Proof}.
To check 1) recall that the arithmetic Chern connection attached to $q^{(s)}$ is the {\it unique} 
$\pi$-connection  of degree $s$ for which the corresponding matrices $\Lambda_i^{(s)}$ satisfy the equalities
$$\begin{array}{rcl}
\Lambda_i^{(s)}A_i^{(s)}\Lambda_i^{(s)} & = & B^{(s)},\\
\ & \ & \ \\
A_i^{(s)}\Lambda_i^{(s)} & = & \Lambda_i^{(s)t}A_i^{(s)}.\end{array}$$
Applying $\tau_j^G$ to these equalities,
setting $\Lambda_{i,j}^{(s)}:=\tau_j^G(\Lambda_i^{(s)})$, 
 and using Lemma \ref{spanissh} we get
$$\begin{array}{rcl}
\Lambda_{i,j}^{(s)}A_i^{(s)}\Lambda_{i,j}^{(s)} & = & B^{(s)},\\
\ & \ & \ \\
A_i^{(s)}\Lambda_{i,j}^{(s)} & = & \Lambda_{i,j}^{(s)t}A_i^{(s)}.\end{array}$$
By uniqueness we get
$$\tau_j^G(\Lambda_i^{(s)})=\Lambda_i^{(s)}$$
for all $i,j$. Then assertion 1 follows from Lemma \ref{aeroplann}. 

\

To check assertion 2 recall that the arithmetic Levi-Civita connection attached to $q^{(s)}$ and 
$\frac{p}{\pi}L=(\frac{p}{\pi}L^{k(s)}_{ij})$ is the {\it unique} 
$\pi$-connection  of degree $s$  for which the corresponding matrices $\Lambda_i^{(s)}$ satisfy the equalities
$$\begin{array}{rcl}
\Lambda_i^{(s)}A_i^{(s)}\Lambda_i^{(s)} & = & B^{(s)},\\
\ & \ & \ \\
(\Lambda_i^{(s)}-1)_{kj}-(\Lambda_j^{(s)}-1)_{ki} & = & \pi\cdot  \frac{p}{\pi}\cdot L_{ij}^{k(s)}(\Lambda^{(s)})\\
\ & \ & \ \\
\ & = & 
p L_{ij}^{k(s)}(\Lambda_1^{(s)},\ldots,\Lambda_n^{(s)})\end{array}.$$
Applying $\tau_l^G$ to these equalities, setting $\Lambda_{i,l}^{(s)}:=\tau_l^G(\Lambda_i^{(s)})$, and using Lemma \ref{spanissh} we get
$$\begin{array}{rcl}
\Lambda_{i,l}^{(s)}A_i^{(s)}\Lambda_{i,l}^{(s)} & = & B^{(s)}\\
\ & \ & \ \\
 (\Lambda_{i,l}^{(s)}-1)_{kj}-(\Lambda_{j,l}^{(s)}-1)_{ki} & = &
p L_{ij}^{k(s)}(\Lambda_{1,l}^{(s)},\ldots,\Lambda_{n,l}^{(s)}),\end{array}$$
because $L^{(s)}$ has coefficients in $R$ and $\tau_l$ acts as the identity on $R$.
By uniqueness we get
$$\tau_l^G(\Lambda_i^{(s)})=\Lambda_i^{(s)}$$
for all $i,l$. Then 2) follows, again, from Lemma \ref{aeroplann}.
\qed

\section{Legendre symbols}
As in \cite{Bu19} metric connections are related, in the simplest cases, to Legendre symbols and can be viewed, roughly,  as matrix analogues of Legendre symbols. We explain this in our PDE context here. We note, this suggests even deeper connections between the cohomology classes introduced in Section \ref{nucah} and class field theory. 

\

We introduce the following notation: for every finite extension $L/\mathbb Q_p$ we write
$$\mathcal A_L:=\mathcal O_{L}[x,\det(x)^{-1}]^{\widehat{\ }}.$$

\

Let $\pi$ be a prime element in a Galois extension 
 $E$ of $\mathbb Q_p$. Let  $F$ be the maximum unramified extension of $\mathbb Q_p$ contained in $E$ and let the cardinality of the residue field of $F$ be $p^f$. Recall 
 that $\pi$ is a  root of an Eisenstein polynomial with coefficients in $\mathcal O_F$ of some degree $n$ and $E=F(\pi)$; also $K_{\pi}\simeq E\otimes_F K$.
Write 
$$\mathfrak G(K_{\pi}/K)=\{\sigma_1,\ldots,\sigma_n\},$$
and fix  a Frobenius automorphism $\phi\in \mathfrak F^{(1)}(K_{\pi}/\mathbb Q_p)$ of degree $1$.
 Note that the restriction homomorphism
$$\mathfrak G(K_{\pi}/K)\rightarrow \mathfrak G(E/F)$$
is an isomorphism. We continue to denote by $\sigma_i$ the corresponding elements
in $G(E/F)$ and by $\phi$ the restriction of $\phi$ to $E$. 

\

 We may consider the higher $\pi$-Frobenius lifts $\phi_i^{(s)}:=\phi^s\sigma_i$ on $R_{\pi}$; every $\pi$-Frobenius lift on $R_{\pi}$ of degree $s$ has this form. For each $s\geq 1$ we  view
 $R_{\pi}$ with its induced structure of partial $\delta$-ring  of degree $s$ 
 defined by the tuple $(\phi_1^{(s)},\ldots,\phi^{(s)}_n)$.
 Note that 
$\phi^f$ is the identity on $F$. 

\

Now consider a symmetric matrix $q^{(s)}$ with entries in the ring of integers $\mathcal O_{E}$ of $E$ and consider an arbitrary  $\pi$-connection $\Delta^{(s)G}$ of degree $s$ on $G$, with attached $\pi$-Frobenius lifts
$$\Phi^{(s)G}=((\phi_1^{(s)})^G,\ldots,(\phi_n^{(s)})^G),$$
   that is metric with respect to $q^{(s)}$. Let $A^{(s)}_i,B^{(s)}$ be the associated matrices in Equations  (\ref{1955}). 
The matrices $A_i^{(s)},B^{(s)}$ have entries in the ring $\mathcal A_E$.
From now on we assume that the matrices $\Lambda_i^{(s)}$ also have entries in $\mathcal A_E$. 
This is the case (in view of the proofs of Theorems \ref{ch} and \ref{LCC}) in each of the following two cases:

\

{\bf (a)}  $\Delta^{(s)G}$ is  the arithmetic Chern connection $\Delta^{(s)\textup{Ch}}$ attached  to $q^{(s)}$.

\

{\bf (b)}   $\Delta^{(s)G}$   is  the arithmetic Levi-Civita connection $\Delta^{(s)\textup{LC}}$
attached to $q^{(s)}$ and to a torsion symbol $L^{(s)}$ that has components with coefficients in $\mathcal O_E$. 

\

By Equation
(\ref{LALB}) we have $\Lambda_i^{(s)t}A_i^{(s)}\Lambda_i^{(s)}=B^{(s)}$. Taking determinants we get
$$\det(\Lambda_i^{(s)})^2=\det(B^{(s)})\cdot \det(A_i^{(s)})^{-1}\equiv \det((q^{(s)})^{(p^s)})\cdot (\det(\phi_i^{(s)}(q^{(s)})))^{-1}\ \ \textup{mod}\ \ (x-1).$$
Assume now that $q^{(s)}$ is diagonal with entries $q^{(s)}_1,\ldots,q^{(s)}_n$ and set $D:=q^{(s)}_1\ldots q^{(s)}_n$. We get
\begin{equation}
\label{frantz}
\det(\Lambda_i^{(s)})^2\equiv  \frac{D^{p^s}}{\phi_i^{(s)}(D)}
\ \ \textup{mod}\ \ (x-1).
\end{equation}
Set
\begin{equation}
\lambda_i:=\det(\Lambda_i^{(s)})_{|x=1}\in \mathcal O_{E}.\end{equation}
Then we have
\begin{equation}
\label{franzjosef}
\lambda_i^2=  \frac{D^{p^s}}{\phi_i^{(s)}(D)}.
\end{equation}
For $i,j\in \{1,\ldots,n\}$ define $\phi(i),i\circ j\in \{1,\ldots,n\}$ by the equalities
\begin{equation}
\label{defcircc}
\phi^{-1} \sigma_i \phi=\sigma_{\phi(i)}\ \ \textup{and}\ \ \sigma_i\sigma_j=\sigma_{i\circ j}\end{equation}
in $\mathfrak G(E/F)$ (equivalently in $\mathfrak G(K_{\pi}/K)$). In particular $i\mapsto \phi(i)$,  $j\mapsto i\circ j$, and $j\mapsto i\circ j$ are bijections from $\{1,\ldots,n\}$ to itself; we denote by $i\mapsto \phi^s(i)$ the $s$-th iterate of $i\mapsto \phi(i)$.
Let
$$\textup{\bf N}:=\textup{\bf N}_{E/F}:E^{\times}\rightarrow F^{\times}$$
 be the norm map. 
 
 \

 For $f|s$ we have that $\phi^s$ is the identity on $F$ so in this case  we have
$$
\textup{\bf N}\left(\frac{D^{p^s}}{\phi_i^{(s)}(D)}\right) =   \frac{\textup{\bf N}(D)^{p^s}}{\textup{\bf N}(\phi_i^{(s)}(D))}= \frac{\textup{\bf N}(D)^{p^s}}{\prod_{k=1}^n \sigma_k
\phi_i^{(s)}(D)}=\frac{\textup{\bf N}(D)^{p^s}}{\prod_{k=1}^n \phi^s \sigma_{\phi^s(k)\circ i}(D)}
$$
\medskip
$$=\frac{\textup{\bf N}(D)^{p^s}}{\phi^s(\prod_{k=1}^n  \sigma_{\phi^s(k)\circ i}(D))}
=\frac{\textup{\bf N}(D)^{p^s}}{\phi^s(\textup{\bf N}(D))}=\frac{\textup{\bf N}(D)^{p^s}}{\textup{\bf N}(D)}=\textup{\bf N}(D)^{p^s-1}.
$$
Note that one has a natural norm map (induced by $\textup{\bf N}$, still denoted by $\textup{\bf N}$):
$$\textup{\bf N}:\mathcal A_E\rightarrow \mathcal A_F.$$

\

By Equation (\ref{franzjosef}) if $f|s$ we have
\begin{equation}
\label{oa}
\textup{\bf N}(\lambda_i)^2= \textup{\bf N}(D)^{p^s-1}.
\end{equation}
Note by the way that the right hand side of Equation (\ref{oa}) is independent of $i$.
Consider  the {\bf Legendre symbol}
$$\mathcal O_F^{\times}\rightarrow \{\pm 1\},\ \ a\mapsto \left(\frac{a}{p}\right),$$
which is $1$ if $a$ is a square in $\mathcal O_F/(\pi)$ (equivalently,   in $\mathcal O_F$) and is $-1$ otherwise. 
Clearly for all $a\in \mathcal O_F$ we have that 
$$a^*:=\left(\frac{a}{p}\right) \cdot a^{\frac{p^s-1}{2}}$$
is the unique element in $\mathcal O_F^{\times}$ such that $a^*\equiv 1$ mod $p$ and $(a^*)^2=a^{p^s-1}$.

\

We have the following result, where we also recall the hypotheses used in the previous computations.

\begin{theorem}\label{legendre1}
Assume $n=d$,  $f|s$, $q^{(s)}$ is diagonal with entries in $\mathcal O_{E}$, and $\Lambda_i^{(s)}$ have coefficients in $\mathcal O_{E}$. Then
for all $i$ we have:
\begin{equation}
\label{llp}
\textup{\bf N}(\lambda_i)=\left(\frac{\textup{\bf N}(D)}{p}\right) \cdot \textup{\bf N}(D)^{\frac{p^s-1}{2}}.\end{equation}
\end{theorem}

{\it Proof}.
In view of Equation (\ref{oa}) it is enough to check that
 $\textup{\bf N}(\lambda_i)\equiv 1$ mod $p$. Note that since $\det(\Lambda_i^{(s)})\equiv 1$ mod $\pi$ in $\mathcal O_{E}[x,\det(x)^{-1}]^{\widehat{\ }}$
 we have $\lambda_i\equiv 1$ mod $\pi$ in $\mathcal O_{E}$.
 Hence $\textup{\bf N}(\lambda_i)\equiv 1$ mod $p$ in $\mathcal O_F$.
  \qed.

\

\begin{remark}
We  take the opportunity here to correct an omission in
\cite[Equation (2.7)]{Bu19}: for that equation to hold one needs the additional hypothesis 
that $q$ is diagonal.
\end{remark}

In what follows we drop the hypothesis that $f|s$ but we add the hypothesis
$$F(\sqrt{D})\subset E.$$
For notational simplicity we assume $\sigma_1$ in $\mathfrak G(K_{\pi}/K)$ is the identity: $\sigma_1=1$.
Dividing the Equation (\ref{franzjosef}) by the same equation for $i=1$ we have:
\begin{equation}
\label{mariatheresa}
\frac{\lambda_1^2}{\lambda_i^2}=\frac{\phi_i^{(s)}(D)}{\phi_1^{(s)}(D)}=\phi^s\left(
\frac{\sigma_i D}{D}\right)=\left(\phi^s\left(
\frac{\sigma_i \sqrt{D}}{\sqrt{D}}\right)\right)^2.
\end{equation}

\

\

\noindent Together, we use this to prove the following theorem. 

\

\begin{theorem}\label{recc}
Assume $n=d$,   $q^{(s)}$ is diagonal with entries in $\mathcal O_{E}$, $\Lambda_i^{(s)}$ have coefficients in $\mathcal O_{E}$, and $F(\sqrt{D})\subset E$. Then for all $i$ we have an equality:
\begin{equation}
\label{josef2}
\frac{\lambda_1}{\lambda_i}=\phi^s\left(
\frac{\sigma_i \sqrt{D}}{\sqrt{D}}\right).
\end{equation}
\end{theorem}

{\it Proof}.
 Since the squares of the right and left hand sides of Equation (\ref{josef2})
are equal (cf. Equation \ref{mariatheresa}) it is enough to show that both sides are equivalent to $1 \bmod \pi$.

\

First, since $\sigma_i$ is the identity on $F$ it induces the identity on the residue field of $\mathcal O_F$ which is the same as the residue field of $\mathcal O_{E}$ hence $\sigma_i \sqrt{D}\equiv \sqrt{D}$ mod $\pi$. Hence the right hand side of Equation (\ref{josef2}) is $\equiv 1$ mod $\pi$. On the other hand, since $\det(\Lambda_i^{(s)})\equiv \det(\Lambda_1^{(s)}) \equiv 1$ mod $\pi$ in $\mathcal A_E$
 we have $\lambda_i\equiv \lambda_1\equiv 1$ mod $\pi$ so the left hand side of Equation (\ref{josef2}) is also $\equiv 1$ mod $\pi$.\qed
 
 \

 \begin{remark}
 If $D\in F$ and $\sqrt{D}\not\in F$ the right hand side of Equation (\ref{josef2}) is the
 {\bf Legendre symbol} $(\sigma_i,F(\sqrt{D})/F)$ defined as the
  image of $\sigma_i$ via the  projection
 $$\mathfrak G(E/F)\rightarrow \mathfrak G(F(\sqrt{D})/F)\simeq \{\pm 1\}.$$
 \end{remark}
 
 \
 
 \begin{corollary}
 Under the assumptions of Theorem \ref{recc} and assuming $D\in F$ and $\sqrt{D}\not\in F$ we have that for all $i$ and $j$,
 $$\frac{\lambda_i}{\lambda_j}=\frac{(\sigma_i,F(\sqrt{D})/F)}{(\sigma_j,F(\sqrt{D})/F)}.
 $$
 \end{corollary}

\appendix 

\section{Analogies with classical differential geometry}
\label{diskus}
We explain here some analogies between our arithmetic formalism and the formalism of classical differential geometry. 

\

\subsection{Classical differential geometry}The material on classical differential geometry is, of course, well known and was included only in order to introduce notation that can be used to emphasize the analogies with the arithmetic setting.

\

\subsubsection{Manifolds and vector fields}
Let $M$ be an  $n$-dimensional connected manifold and 
denote by $\mathfrak X(M)$ the $C^{\infty}(M)$-module (and $\mathbb R$-Lie algebra) of vector fields, i.e. of $\mathbb R$-derivations on $C^{\infty}(M)$. For $X\in \mathfrak X(M)$ and a point $P\in M$ we denote by $X_P$ the corresponding tangent vector at $P$.
Suppose $M$ is parallelizable i.e., the tangent bundle is trivial. E.g., $M$ could be an open set in the Euclidean space $\mathbb R^n$ or a Lie group. In our analogy with arithmetic differential geometry we may want to consider  a distinguished point  $P_0\in M$; e.g., in the case of a Lie group we could take $P_0$ to be the identity element. 

\

If $M=I\subset \mathbb R$ is an interval containing $0$  and $\delta_t=\frac{d}{dt}$ is the natural vector field (where $t$ is the coordinate on $\mathbb R$) then for $\lambda\in I$ and $f\in C^{\infty}(I)$ (satisfying an appropriate analyticity condition) we have the Taylor formula
\begin{equation}
\label{taylorule}
f(\lambda)=\left(\sum_{i=0}^{\infty} \frac{\lambda^i\delta_t^i f}{i!}\right)(0).
\end{equation}
So, at least formally, evaluation at $\lambda$ of $f$ is given by the evaluation at $0$ of an (infinite!) linear combination of the functions $\delta_t^if$.

\

Going back to a general parallelizable manifold $M$, consider a basis $X_1,\ldots,X_n$ of $\mathfrak X(M)$. Then the Lie bracket satisfies
\begin{equation}
\label{bba}
[X_i,X_j]=\sum_k f_{ij}^k X_k\end{equation}
for some $f_{ij}^k\in C^{\infty}(M)$. 
 
 \  

\subsubsection{Connections}
Assume one is given a connection on $M$ i.e., an $\mathbb R$-linear map
$$\nabla:\mathfrak X(M)\rightarrow \textup{End}_{\mathbb R-\textup{mod}}(\mathfrak X(M)),\ \ X\mapsto \nabla_X$$
satisfying
$$\begin{array}{rcl}
\nabla_{fX} Y & = & f\nabla_X Y\\
\ & \ & \ \\
\nabla_X(fY) & = & X(f)Y+f\nabla_XY\end{array}$$
for $f\in C^{\infty}(M)$ and $X,Y\in \mathfrak X(M)$.
One sets 
\begin{equation}
\label{bau}
\nabla_{X_i}X_j=\sum_k \Gamma_{ij}^k X_k\end{equation}
for  functions $$\Gamma_{ij}^k\in C^{\infty}(M)$$ which we here refer to  as the {\bf Christoffel symbols of the second kind} with respect to $\{X_i\}$. (This terminology is classically reserved to the case when $X_i$ are coordinate vector fields.)

\

A connection $\nabla$
is {\bf symmetric} if
\begin{equation}
\label{offf}
\nabla_{X}Y-\nabla_{Y}X=[X,Y]\end{equation}
for all vector fields $X,Y\in\mathfrak X(M)$.
The the symmetry condition (\ref{offf}) reads
\begin{equation}
\label{llaa}
\Gamma_{ij}^k-\Gamma_{ji}^k=f_{ij}^k.
\end{equation} 

\

\subsubsection{Parallel transport and geodesics}
Consider a  {\bf curve}
 $$c:I\rightarrow M$$ i.e., an immersion where $I$ is an open interval in $\mathbb R$ containing $0$.
A {\bf vector field along $c$} is a ``smooth"
map that attaches to every $t\in I$ a tangent vector at $c(t)$. Every vector field $Y$ on $M$ induces a vector field which we denote by $Y_c$ along $c$ defined by $t\mapsto Y_{c(t)}$ .
Another important example is the {\bf velocity} of $c$, defined by
$$t\mapsto T_{t}(\delta_t),$$ where $T_t:T_tI\rightarrow T_{c(t)}M$ denotes the tangent map and $\delta_t:=\frac{d}{dt}$, 
is a vector field along $c$. 

\

For every curve $c$
 choose (locally) a vector field $X({c})$ on $M$ with the property that 
 $$X(c)_{c}=T_{t}(\delta_t)$$
 for all
$t\in I$. Assume we are given a connection $\nabla$, 
a vector field $Y$ on $M$ and a curve $c$. The {\bf derivative} of $Y_{c}$ is defined to be the vector field along $c$, denoted by 
 $Y'_{c}$, given by
$$t\mapsto (Y'_{c})_{c(t)}:=(\nabla_{X({c})}Y)_{c(t)}.$$
The map $Y'_{c}$ only depends on $c$ and $Y_{c}$  but not on $X(c)$ or $Y$; in fact if $X_i$ are coordinate vector fields, $c$ is given in local coordinates around $P_0$ by an $n$-tuple of functions $(c_1(t),\ldots,c_n(t))$
and $Y=\sum_j w_j X_j$ then we have
\begin{equation}
\label{derofY}
(Y_{c}')_{c(t)}=\sum_k\left(\delta_t(w_k\circ c)(t)+\sum_{ij}\Gamma^k_{ij}(c(t))w_j(c(t))\delta_tc_i\right)X_k.\end{equation}

\

We say that $Y_{c}$ is {\bf parallel} (with respect to $\nabla$) if
\begin{equation}
\label{parallelll}
Y'_{c}=0;
\end{equation}
equivalently if
\begin{equation}
\label{paralleluli}
\delta_t(w_k\circ c)(t)+\sum_{ij}\Gamma^k_{ij}(c(t))(\delta_tc_i)w_j(c(t))=0.
\end{equation}
For every tangent vector $v_0$ at $P_0$ there is a unique vector field along $c$ that is parallel and
whose value at $P_0$ is $v_0$; the map sending $v_0$ into the above parallel vector field is referred to as {\bf parallel transport}.

\

Recall that $c$
is called a {\bf geodesic} through a distinguished point $P_0\in M$ if $c(0)=P_0$ and we have
\begin{equation}
\label{nablaXX}
(\nabla_{X({c})} X({c}))_{c(t)}=0\end{equation}
 for $t\in I$; in other words if $X(c)_{c}$ is parallel.
 So this definition is independent on the choice of $X({c})$ and in fact
 if $X_i$ are coordinate vector fields and  $c$ is given in coordinates around $P_0$ by an $n$-tuple of functions $(c_1(t),\ldots,c_n(t))$ then $c$ is a geodesic through $P_0$ if and only if $c_i(0)=0$ and 
\begin{equation}\label{claseqgeo}
\delta_t^2 c+\sum_{ij}\Gamma^k_{ij}(c(t))\delta_t c_i \delta_t c_j=0.
\end{equation}
For every tangent vector $v_0$ at $P_0$ there is a unique geodesic on $M$ whose tangent vector at $P_0$ is $v_0$.
 
 \

\subsubsection{Curvature}
The {\bf curvature} of a connection $\nabla$ is defined as
\begin{equation}
  \label{nuc}
  \mathcal R(X,Y):=\nabla_X\nabla_Y-\nabla_Y\nabla_X-\nabla_{[X,Y]}:\mathfrak X(M)\rightarrow \mathfrak X(M),\end{equation}
  for $X,Y\in \mathfrak X(M)$. 
  
  \medskip
  
  \subsubsection{Metrics}
    A {\bf metric} on $M$ is given by a symmetric positive definite 
  form
  $$(X,Y)\mapsto \langle X,Y\rangle$$
  where $X,Y$ are vector fields on $M$ hence by a symmetric positive definite 
  $n\times n$ matrix  
  $q=(q_{ij})$, $q>0$, with entries $$q_{ij}:=\langle X_i,X_j\rangle.$$ 
  In case $M=\mathcal F$ is a Lie group a metric on $\mathcal F$ is called {\bf bi-invariant} if it is invariant under both left and right translations. 
  Bi-invariant metrics on $\mathcal F$  are in bijection with the positive definite symmetric bilinear maps
  $$\beta:\mathfrak f\times \mathfrak f\rightarrow \mathbb R$$
  on the Lie algebra $\mathfrak f$ of $\mathcal F$
  satisfying the ``\textup{Ad}-invariance" condition
  $$\beta(\textup{Ad}(g)X,\textup{Ad}(g)Y)=\beta(Y,Z),\ \ \ g\in \mathcal F,  \ \ Y,Z\in \mathfrak f;$$
  equivalently
  $$\beta([X,Y],Z)+\beta(Y,[X,Z])=0,\ \ \ X,Y,Z\in \mathfrak f.$$
  
  \

  As an example, if $\mathcal F=\textup{GL}_m(\mathbb R)$, $\mathfrak f=\mathfrak g\mathfrak l_m(\mathbb R)$ (hence $n=m^2$)
  then the Non-degenerate \textup{Ad}-invariant symmetric bilinear   forms satisfying this condition  are parameterized by pairs of constants 
$(\lambda_1,\lambda_2)\in \mathbb R^{\times}\times \mathbb R^{\times}$ and are given by the formula:
 \begin{equation}
 \label{glgl}
 \beta(a,b)=\lambda_1\cdot \text{tr}(a_0b_0) +\lambda_2\cdot
 \text{tr}(a)\text{tr}(b),\end{equation}
 where 
 $$a_0:=a-\frac{1}{m}\text{tr}(a)\cdot 1_m$$
 denotes the traceless part of a matrix $a$ and
 $1_m$ is the identity matrix.
 Recall that to see this one considers
  the Lie algebra decomposition
 $${\mathfrak g}{\mathfrak l}_m(\mathbb R)\simeq {\mathfrak s}{\mathfrak l}_m(\mathbb R)\oplus {\mathfrak g}{\mathfrak l}_1(\mathbb R),\ \ \ 
 a\mapsto (a_0,a-a_0)$$
 where ${\mathfrak g}{\mathfrak l}_1(\mathbb R)$ is the center of ${\mathfrak g}{\mathfrak l}_m(\mathbb R)$;
 by the irreducibility of the adjoint action of $\textup{SL}_m$ on the first factor one gets that
  any \textup{Ad}-invariant  non-degenerate symmetric bilinear form $\beta$ on ${\mathfrak g}{\mathfrak l}_m(\mathbb R)$ must be an orthogonal sum of two such forms on the two factors;
  also by the irreducibility above plus Schur's Lemma one gets that
    all \textup{Ad}-invariant non-degenerate symmetric bilinear forms on the first factor are 
 multiples of each other. Cf. \cite[pg. 158]{sternberg}.

  \
  
  Going back to the general case  of a connection $\nabla$ on a general $M$
  we say $\nabla$ is {\bf metric}  with respect to a metric $\langle\ ,\ \rangle$ if 
  \begin{equation}
  \label{metcom}
  X\langle Y,Z\rangle=\langle \nabla_X Y,Z\rangle +\langle Y,\nabla_X Z\rangle
  \end{equation}
  for all vector fields $X,Y,Z\in \mathfrak X(M)$. 
  
  \

  Given a connection $\nabla$ and a metric the  {\bf Christoffel symbols
   of the first kind} of $\nabla$ are defined as
  \begin{equation}
  \label{dinno}
  \Gamma_{ijk}:=\sum_m \Gamma_{ij}^mq_{mk}.\end{equation}
  It is also convenient to ``lower the indices" of $f^k_{ij}$ by defining:
  \begin{equation}
  f_{ijk}:=\sum_m f_{ij}^m q_{mk}.
  \end{equation}

\medskip

\subsubsection{Levi-Civita connection}
  Given a metric there is a unique connection that is symmetric and metric and is called the {\bf Levi-Civita connection}; it is given by 
  \begin{equation}
  \label{lilu}
  \Gamma_{ijk}=\frac{1}{2}(X_iq_{jk}+X_jq_{ik}-X_kq_{ij})-\frac{1}{2}
  (f_{kij}+f_{ijk}-f_{jki}).
  \end{equation}
  In the special case when $\{X_i\}$ are coordinate vector fields formula (\ref{lilu}) becomes 
   \begin{equation}
  \label{lilu2}
  \Gamma_{ijk}=\frac{1}{2}(X_iq_{jk}+X_jq_{ik}-X_kq_{ij}).
  \end{equation}
  If $\mathcal R$ is the curvature of the Levi-Civita connection $\nabla$
  one sets
  \begin{equation}
  \mathcal R(X_i,X_j)X_k=\sum_m R^m_{kij}X_m
  \end{equation}
  with smooth functions $R^m_{kij}$ 
  and one defines the {\bf Riemann curvature tensor} by 
  \begin{equation}
  \label{nnuu}
  R_{lkij}:=\sum_m R^m_{kij}q_{lm}.
  \end{equation}
For $\{X_i\}$ coordinate vector fields we have
\begin{equation}
R^l_{kij} =  X_i \Gamma^l_{kj}-X_j \Gamma^l_{ki}+\sum_m \Gamma^l_{mi}\Gamma^m_{kj}-
\sum_m \Gamma^l_{mj}\Gamma^m_{ik}
\end{equation}
and the Riemann curvature tensor satisfies the usual symmetries 
  \begin{equation}
  \label{summery}
  \begin{array}{rcll}
  R_{ijkl} + R_{jikl} & = & 0  \ & \ \\
  \ & \ & \ & \ \\
  R_{ijkl} + R_{ijlk} & = & 0 &  \ \\
  \ & \ & \ & \ \\
  R_{ijkl}+R_{iklj}+R_{iljk} & = & 0 & \ \\
  \ & \ & \ & \  \\
  R_{ijkl} - R_{klij} & = & 0. & \ 
  \end{array}
  \end{equation}
  Moreover one defines in this case the {\bf Ricci} tensor
  \begin{equation}
  R_{ij}:=\sum_k R^k_{ikj}=\sum_{k,l}q^{lk}R_{likj}
  \end{equation}
  and the {\bf scalar curvature}
  \begin{equation}
  \label{rws}
  S:=\sum_{i,j} q^{ij}R_{ij}=\sum_{i,j,k,l} q^{ij}q^{lk}R_{likj},
  \end{equation}
  where $(q^{ij})$ is the inverse of $q=(q_{ij})$.
  Note that for $\{X_i\}$ coordinate vector fields one has:
   \begin{equation}
  R_{ijkl}=R^{(2)}_{ijkl}+R^{(1)}_{ijkl},
  \end{equation}
  where $R^{(2)}_{ijkl}$ is the {\bf order $2$ part} of $R_{ijkl}$ defined by
  \begin{equation}
  \label{zzoo}
  R_{ijkl}^{(2)}=\frac{1}{2}(X_iX_k q_{jl}+X_lX_j q_{ki}-X_iX_l q_{kj}-X_kX_j q_{il}),
  \end{equation}
 and $R^{(1)}_{ijkl}$ is the {\bf order $\leq 1$ part} of $R_{ijkl}$ 
 which is a universal polynomial in 
 the functions
 $$\{q_{st},\det(q)^{-1},X_mq_{st}\ |\  \ m,s,t\in \{1,\ldots,n\}\}.$$
 
 \medskip
 
 \subsubsection{Involutive distributions}
 Some of the above concepts can be generalized to the case when instead of considering a manifold $M$ with trivial tangent bundle $TM$ 
 with basis $X_1,\ldots,X_n$ one considers a manifold
 $M$ with an involutive distribution 
 \begin{equation}
 \label{tingg}
 \mathfrak D\subset \mathfrak X(M)\end{equation}
  with basis vector fields $X_1,\ldots,X_n$,
 $n\leq \dim(M)$. Given such an involutive distribution one may consider the smallest Lie algebra $\mathfrak l$  of 
 $\textup{End}_{\mathbb R-\textup{mod}}(\mathfrak X(M))$ containing all endomorphisms $\nabla_X$ for $X\in \mathfrak D$.  The commutator $\mathfrak h=[\mathfrak l,\mathfrak l]$  corresponds to  the classical  {\bf holonomy} Lie algebra of $\nabla$. 
 
 \medskip
 
 \subsubsection{Cohomology}
  Finally  one can ask  for  cohomology groups that are  natural recipients of cohomology classes
 attached to our geometric concepts.
 For instance, dropping the assumption that $M$ has a trivial tangent bundle and assuming $M$ is a compact algebraic complex manifold recall the {\bf Atiyah  sequence} \cite{A57} of complex algebraic vector bundles
 \begin{equation}
 \label{reberluiza}
 0\rightarrow \textup{End}(T)\rightarrow Q \rightarrow T\rightarrow 0\end{equation}
 where $T$ is the algebraic tangent bundle of $M$ and $Q$ is the bundle whose sheaf of  sections is the sheaf of  invariant
 algebraic vector fields on the principal $\textup{GL}_n$-bundle associated to $T$.  The obstruction to the splitting of the Atiyah extension ``gives rise" to the characteristic classes of $M$ living in the de Rham cohomology
 $$H^*_{DR}(X,\mathbb C).$$ 
 On the other hand the
  main examples of parallelizable manifolds we have in mind are  Lie groups. 
 If $M=\mathcal F$ is a compact Lie group we recall Cartan's theorem according to which 
 we have a canonical  isomorphism 
 \begin{equation}
 \label{ecartan}
 H^*_{DR}(\mathcal F)\simeq H^*(\mathfrak f,\mathbb R)\end{equation}
 between  the de Rham cohomology of $\mathcal F$ and  the Lie algebra cohomology 
  of the Lie algebra $\mathfrak f$ of $\mathcal F$.

\subsection{Comparison between classical and arithmetic differential geometry}
The analogy between classical differential geometry 
and our arithmetic differential geometry is as follows.

\

Our ring $R_{\pi}$ is the analogue of the 
 ring of smooth functions $C^{\infty}(M)$ on a parallelizable  manifold $M$ (possibly with a distinguished point $P_0\in M$). The higher
 $\pi$-Frobenius lifts $\phi_i^{(s)}$ (alternatively  the higher $\pi$-derivations $\delta_i^{(s)}$) on $R_{\pi}$
 are the analogues of the vector fields 
  $X_i$ forming a basis for the module of vector fields. 
 
\

 Since $R_{\pi}$ is a local ring one is tempted to replace, in this analogy, the manifold $M$ by the ``germ" of $M$ at  $P_0$. However, in some sense, $R_{\pi}$ and various objects attached to it 
 behave, in our theory, 
 as if they possessed some ``global" features. This is suggested, for instance, by the  fact that the $\pi$-Frobenius lifts on $R_{\pi}$ generally do not commute (and this cannot be remedied by a ``change of basis" since the only changes of bases in our theory correspond to permutations of the $\pi$-Frobenius lifts). Moreover, even in case the $\pi$-Frobenius lifts on $R_{\pi}$ do commute, our theory leads to
  non-trivial cohomology classes. 
  Contrast this with the fact that every point on a manifold has a neighborhood that has a ``trivial" topology and admits a commuting basis of vector fields.

  \

  So thinking of $R_{\pi}$ as the analogue of a ``germ of a manifold at a point" may be, in many ways, misleading. In a number of situations, the best analogy seems to be between $R_{\pi}$ and the ring of smooth functions on a Lie group. In particular, the ring $R$ can be viewed as as the analogue of the ring of smooth functions on $\mathbb R$. Our definition in 
  \cite[Def. 2.22]{BM22} can be viewed as an arithmetic analogue of the  Taylor formula (\ref{taylorule}) with 
  expressions $P(a,\delta^{(s)} a,\ldots)$ for
  polynomials $P\in R\{t\}$ defined in loc. cit.  being analogous to expressions $\sum_i \lambda^i\delta^i_t f/i!$.
  
  \

  Let $P$ be the  frame bundle on $M$.
 The vector fields $X_i$ define an identification  $$P\simeq M\times \textup{GL}_n.$$
  Then
 our ring $\mathcal A=R_{\pi}[x,\det(x)^{-1}]^{\widehat{\ }}$ is the analogue of 
  the ring $C^{\infty}(P)$ and the indeterminates $x_{ij}\in \mathcal A$ are analogues of the projections  $$x_{ij}:P\rightarrow \textup{GL}_n\rightarrow \mathbb R$$ defined by the coordinate functions on $\textup{GL}_n$.
  
  \

  If $X_i^{P,0}$ are the vector fields on $P$ lifting $X_i$ and killing the coordinate functions $x_{jk}$ then
  the higher $\pi$-derivations on $\mathcal A$ defining a $\pi$-connection are analogues 
  of the vector fields $X^P_i$ on $P$ lifting $X_i$ naturally induced by the ``dual" of the connection $\nabla$, i.e., 
  $$X^P_i:=X_i^{P,0}-\sum_{jkl}\Gamma_{il}^jx_{lk} \frac{\partial}{\partial x_{jk}}.$$
  
  \

  Our  Lie symbol $(\ell_{ij}^{k(s)})$ in Definition \ref{ccaann}, cf. Equation (\ref{nygov}), is the analogue of the family $(f_{ij}^k)$ in Equation (\ref{bba}). Note that our $\ell_{ij}^{k(s)}$'s
  in the context of Definition \ref{ccaann} belong to $\{-1,0,1\}$ so they are killed by the higher $\pi$-derivations and should be viewed as ``constants" in the theory. The analogous situation in classical differential geometry is that in which the functions $f_{ij}^k$ are constant; this can always be achieved locally (by taking coordinate vector fields in which case  $f_{ij}^k=0$) and sometimes globally (as with Lie groups in which case $f_{ij}^k$ are the structure constants of the Lie algebra of our Lie group). 
  
  \

  Our Christoffel symbols of the second kind $\Gamma_{ij}^{k(s)}$ 
  are the analogues of the classical
  Christoffel symbols
   $\Gamma_{ij}^k$ in Equation (\ref{bau}): one can view Equation (\ref{penu})
   as an analogue of Equation (\ref{bau}). Our Christoffel symbols of the first kind 
   are analogous to the classical Christoffel symbols of the first kind in Equation (\ref{dinno}).
   
   \

   Next note that a curve $c:I\rightarrow M$ in a manifold is given (in coordinates $\xi_1,\ldots,\xi_n$ on $M$) by an $n$-tuple $(c_1(t),\ldots,c_n(t))$ of functions in $C^{\infty}(I)$; since $R$ can be viewed as the analogue of $C^{\infty}(I)$ our definition of a curve  in the arithmetic case 
  as being a vector $c\in R^n$ (see \cite[Def. 3.7]{BM22})
   is analogous to the definition of a curve in a manifold.
   Recall that in Part 1, loc.cit. we fixed an $R$-module basis of $R_{\pi}$ and hence curves $c\in R^n$ are identified with
    the $R$-linear maps 
   $$c^*:R_{\pi}\rightarrow R;$$ 
   these maps should be viewed as
    analogous to the restriction maps
   $$c^*:C^{\infty}(M) \rightarrow C^{\infty}(I).$$
    Note however that $c^*$ in the classical case are ring homomorphisms whereas in our arithmetic context the maps $c^*$ are  $R$-module homomorphisms. (Of course, there are no ring homomorphisms $R_{\pi}\rightarrow R$ which is why we made the somewhat bold move of taking $R$-module homomorphisms $c^*:R_{\pi}\rightarrow R$ as analogues of restriction maps between rings of functions; this is consistent with the fact that, in order to adapt concepts to  geometries over the field with one element, one sometimes has to sacrifice   compatibility with one of the two operations: addition or multiplication.)
  In the same vein a vector field along a curve $c:I\rightarrow M$ in a manifold has the form
  $$\sum_{i=1}^n w_i(t)\left(\frac{\partial}{\partial \xi_i}\right)_{c(t)}.$$
  for an $n$-tuple $(w_1(t),\ldots,w_n(t))$ of functions in $C^{\infty}(I)$.
  This makes pairs $(c,w)\in R^n\times R^n$ analogous to vector fields along a curve in a manifold.
  
   \

   The definition (\ref{parallelll}) of classical parallel vector fields along a curve $c:I\rightarrow M$ is analogous to our definition \cite[Def. 3.8]{BM22} of parallel vectors in $R^n$ and similarly
   the definition (\ref{nablaXX}) of classical geodesics $c:I\rightarrow M$ is analogous to our definition \cite[Def. 3.10]{BM22} of geodesics $c\in R^n$ in the arithmetic case.
    Note that  
    the classical equation for parallelism (\ref{paralleluli}) is analogous to  \cite[Eq. 3.12]{BM22}; similarly
    the classical equation (\ref{claseqgeo}) for geodesics is analogous to \cite[Eq. 3.13]{BM22}. 
   The theorems \cite[Thm. 
     3.16 and 3.17]{BM22} are analogous to the classical differential geometric theorems claiming the 
    the existence of parallel transport and the
   existence and uniqueness of geodesics, respectively. Note however that in the arithmetic case the exponential map constructed via geodesics is not defined on the whole ``tangent spaces" but rather on the ``tangent space minus a hyperplane."
   
   \

Our definition of symmetry for $\pi$-connections  of degree $s$ 
with respect to the  torsion symbol $L_{ij}^{k(s)}$ in Equation (\ref{uf2})
is analogous to the classical differential geometric definition in Equation (\ref{llaa}) of symmetry for classical connections. One can take $L^{(s)}=(L_{ij}^{k(s)})$ to be equal to the Lie symbol $\ell^{(s)}=(\ell^{(s)}_{ij})$
or to a ``deformation" of the Lie symbol (as in the case of the multiplicative canonical torsion symbol $\ell^{(s)*}$). Then our $\delta$-Lie algebra of $\textup{GL}_n$ can be viewed as a (power of a) deformation of the classical Lie algebra of $\textup{GL}_n$; a similar construction can be made for an arbitrary group scheme in place of $\textup{GL}_n$.

\

Our invertible symmetric matrices $q$ with coefficients in $R_{\pi}$ are analogues of metrics in classical Riemannian geometry; 
as noticed in the paper it is sometimes natural to  further restrict attention to $q$'s all of whose 
diagonal entries are invertible in $R_{\pi}$, the latter condition being reminiscent of the condition that metrics in the classical case are positive definite. Our \textup{Ad}-invariance condition (cf. Definition \ref{adinv}) is analogous to the condition that a metric on a Lie group is bi-invariant. Bi-invariant 
metrics on Lie groups  depend ``linearly" on a number of parameters that is dictated by representation theoretic data (as seen, for instance, in the example of Equation (\ref{glgl})); this is also the case for \textup{Ad}-invariant metrics in the arithmetic case (as seen, for instance, in Example \ref{399}).

\

Our metric condition for a $\pi$-connection  of degree $s$ in \cite[Def. 3.25]{BM22} is analogous to the classical metric condition in Equation (\ref{metcom}). Indeed one can check that the condition in Equation (\ref{metcom}) is equivalent to the commutativity of the following diagrams
\begin{equation}
\label{dididi}
\begin{array}{rcl}
C^{\infty}(P) & \stackrel{X_i^{P,0}}{\longrightarrow} & C^{\infty}(P)\\
\mathcal H_q \downarrow & \ & \downarrow \mathcal H_q\\
C^{\infty}(P) & \stackrel{X_i^P}{\longrightarrow} & C^{\infty}(P)
\end{array}
\end{equation}
where 
$\mathcal H_q$ is the ring homomorphism induced by the map 
$$P\rightarrow P,\ \ (\xi,g)\mapsto (\xi,g^tqg),\ \ (\xi,g)\in P=M\times \textup{GL}_n.$$ The diagrams (\ref{dididi}) are analogous to the diagrams in \cite[Eq.3.22]{BM22}.
The formula giving the Christoffel symbols in Theorem \ref{LCC} is analogous to  formula (\ref{lilu}); correspondingly formula \cite[Eq. 3.28]{BM22} is analogous to the formula
(\ref{lilu2}). 

\

The Lie algebras of the form $\mathfrak d$ attached to involutive subsets 
$$\mathfrak D^{(c)}\subset \mathfrak F^{(c)}$$
in Definition \ref{ccaann}
are analogous to  involutive distributions (Equation \ref{tingg}) 
in classical differential geometry.
 The curvature $\mathcal R$ of a connection  on $M$ in Equation (\ref{nuc}) has, as analogue, the curvature $\mathcal R$ in Equation (\ref{nua}) and also the multiplicative curvature $\mathcal R^*$ in Definition \ref{eei}. 
 The classical Riemann curvature tensor in formula (\ref{nnuu}) has, as analogue, our arithmetic Riemann curvature tensor; cf. Definition \ref{arct}. Our formula for the arithmetic Riemann curvature tensor in Theorem \ref{tery} is similar to the ``order $2$ part" in the classical formula
 for the Riemann curvature tensor in Equation (\ref{zzoo}). The symmetries in Corollary
 \ref{antisymm} are analogous to the classical symmetries (\ref{summery}).
 The $I$-curvatures  $I(q,\overline{R})$ in Equation (\ref{scalar}) are analogues of the classical scalar curvature $S$ in Equation (\ref{rws}). 
 The Lie algebra of vector fields $\mathfrak X(M)$ has, as analogue, our canonical algebra $\mathfrak f$ (cf. Definition \ref{ccaann}).
 For a classical connection
 $\nabla$,  the smallest Lie algebra  $L_{\nabla}$
  containing the image of the natural map
 $$\nabla:\mathfrak X(M)\rightarrow \mathfrak X_{\textup{inv}}(P)=\{\textup{invariant vector fields on}\ P\},$$
  has, as analogue, the Atiyah extension $\mathfrak l_{\nabla}$ attached to a
 graded $\pi$-connection; cf. Definition \ref{liehulll}. 
 Recall that one can define the {\bf holonomy Lie algebra} of a classical connection $\nabla$ as the commutator of $L_{\nabla}$ above. This would make the commutator Lie ideal 
 $[\mathfrak l_{\nabla},\mathfrak l_{\nabla}]$
 of  $\mathfrak l_{\nabla}$ an arithmetic analogue of the holonomy Lie algebra.
 
 \

 The  classical Atiyah extension (\ref{reberluiza}) is analogous to the  extension  (\ref{nutziiii}). The obstruction to splitting the classical Atiyah extension (in the category of vector bundles) is what gives rise to the characteristic classes of $M$. On the other hand the  sequence (\ref{nutziiii}) is split in the category of graded Lie $\mathbb Z$-algebras so our analogy seems to break at this point. 

 \

 A way to restore the analogy is to interpret our Atiyah extension $\mathfrak l_{\nabla}$   as a refinement of (\ref{nutziiii}) and this refinement {\it does} lead to non-trivial classes in  Lie algebra cohomology; the latter classes can be viewed as arithmetic analogues of ``refinements"  of Chern classes (that depend on the metric but only in a mild manner, since they can be viewed as ``discrete invariants" of the metric). The analogue, in classical differential geometry, of (\ref{nutziiii}) would be the abelian extension of Lie algebras
 \begin{equation}
 \label{iouaz}
 0\rightarrow N^{\textup{ab}}\rightarrow L/[N,N]\rightarrow  \mathfrak X(M)\rightarrow 0\end{equation} 
 where $L=L_{\nabla}$, $N$ is the kernel of the natural homomorphism $L_{\nabla}\rightarrow \mathfrak X(M)$, and $N^{\textup{ab}}:= N/[N,N]$. 
  This construction can be extended, of course,  to the case when $M$ is not necessarily parallelizable.

  \
  
  Going back to the parallelizable case, 
in accordance with our analogy between the ``space" corresponding to $R_{\pi}$ and a Lie group, and in view of Cartan's isomorphism (\ref{ecartan}), it is tempting to take the Lie cohomology groups $H^2_{\textup{Lie}}(\mathfrak d,\mathfrak n/\mathfrak n^2)$
or  $H^2_{\textup{Hoch}}(\mathfrak d,\mathfrak n/\mathfrak n^2)$
as recipients for our cohomology classes
attached to our arithmetic objects and hence as analogues of the de Rham (or Betti) cohomology.
In fact, as already mentioned, the groups $H^2_{\textup{Lie}}(\mathfrak d,\mathfrak n_{\nabla}/\mathfrak n_{\nabla}^2)$ and 
$H^2_{\textup{Lie}}(\mathfrak d_{\nabla},\mathfrak n_{\nabla}/\mathfrak n_{\nabla}^2)$
have a discrete/combinatorial flavor which makes them somewhat analogous to a topologically defined object; and moreover the classes in these groups attached to $\nabla$ are images of classes in $H^2_{\textup{Lie}}(\mathfrak d,\mathfrak n_{\mathfrak D}/\mathfrak n_{\mathfrak D}^2)$ and 
$H^2_{\textup{Lie}}(\mathfrak d_{\mathfrak D},
\mathfrak n_{\mathfrak D}/\mathfrak n_{\mathfrak D}^2)$ that are independent of $\nabla$, a situation similar to the one in Chern-Weil theory.

\

Finally note the analogy between our ``torsors" $\mathcal A^u$ in (\ref{eqtor})
and the rings of smooth functions $C^{\infty}(P(E))$ on the principal bundles $P(E)$ attached to vector bundles $E$ on a manifold $M$. However, unlike in the case of classical differential geometry where the Chern classes of $E$ are ``remembered" by curvature, the $1$-cocycles $u$ used to construct our  ``torsors" $\mathcal A^u$ are  apparently not ``remembered" by our arithmetic curvature; indeed the class of $u$ in $H^1_{\textup{gr}}(\mathfrak T,\mathfrak W)$ and the Lie characteristic class of $\nabla$ in $H^2_{\textup{Lie}}(\mathfrak d,\mathfrak n_{\nabla}/\mathfrak n^2_{\nabla})$ are apparently not related in our formalism developed so far. It would be interesting, in fact, to find a formalism that relates generalizations of these two classes.
 
 \

We summarize our discussion above in the following table comparing classical and arithmetic differential geometry (DG):

\bigskip

\begin{center}
\begin{tabular}{||c|c||} \hline \hline
classical DG  & arithmetic DG \\ 
\hline
$C^{\infty}(M)$ & $R_{\pi}$ \\
\hline
 $\sum_i \lambda^i\delta^i_t f/i!$ & $P(a,\delta a,\ldots)$\\
\hline
$X_i$ & $\phi_i^{(s)},\delta_i^{(s)}$ \ \textup{on}\ $R_{\pi}$\\
\hline
$C^{\infty}(P)$ & $\mathcal A$\\
\hline
$E$ & $u\in H^1_{\textup{gr}}(\mathfrak T,\mathfrak W)$\\
\hline
$C^{\infty}(P(E))$ & $\mathcal A^u$\\
\hline
$X_i^{P,0}$ & $\phi_{i,0}^{(s)},\delta_{i,0}^{(s)}$\ \textup{on}\ $\mathcal A$ \\
\hline
$X_i^P, \nabla$ & $\phi_i^{(s)},\delta_i^{(s)}$\ \textup{on}\ $\mathcal A$ \\
\hline
$c:I\rightarrow M$ & $c\in R^n$, $c^*$\\
\hline
$Y_{c}$ & $(c,w)$\\
\hline
$x_{ij}$ & $x_{ij}$\\
\hline
$q_{ij}$ & $q^{(s)}_{ij}$\\
\hline
$q>0$  & $q_{hh}\in R_{\pi}^{\times}$,\\
\hline
$q$ \ \ \textup{bi-invariant} & $q$\ \ \textup{Ad-invariant}\\
\hline
$f_{ij}^k$ & $\ell_{ij}^{k(s)}, \ell^{k(s)*}_{ij}, L^{k(s)}_{ij}$\\
\hline
$\Gamma_{ij}^k$ & $\Gamma_{ij}^{k(s)}$\\
\hline
$\Gamma_{ijk}$ & $\Gamma_{ijk}^{(s)}$\\
\hline
$c,c^*$ & $c^*$\\
\hline
$\mathcal R$ & $\mathcal R,\mathcal R^*$\\
\hline
$R_{ijkl}$ & $R_{ijkl}$\\
\hline
$S$ & $I(q,\overline{R})$\\
\hline
$\mathfrak X(M)$ & $\mathfrak f$\\
\hline
$\mathfrak D, T$ & $\mathfrak D^{(c)},\mathfrak d$\\
\hline
$L=L_{\nabla}$ & $\mathfrak l=\mathfrak l_{\nabla}$ \\
\hline
$\textup{End}(T)\rightarrow  Q\rightarrow T$ & $\mathfrak r\rightarrow \mathfrak e\rightarrow \mathfrak d$\\
\hline
$N^{\textup{ab}}\rightarrow L/[N,N]\rightarrow \mathfrak X(M)$ & $\mathfrak n/\mathfrak n^2\rightarrow \mathfrak l/\mathfrak n^2\rightarrow \mathfrak d$\\
\hline 
$H^2_{DR}(G,\mathbb R)$ & $H^2_{\textup{Lie}}(\mathfrak d,\mathfrak n/\mathfrak n^2)$\\
\hline
\end{tabular}
\end{center}

\bigskip
\bigskip

\bibliographystyle{amsalpha}

\end{document}